\documentclass[english, hidelinks]{article}
\usepackage{geometry}
\usepackage{setspace}
\usepackage{authblk}
\usepackage[utf8]{inputenc}
\usepackage{color,soul,yhmath}
\usepackage{xcolor}
\usepackage{algpseudocode}
\usepackage{bm}
\usepackage{graphicx}
\graphicspath{ {.} }
\usepackage{hyperref}
\usepackage{natbib}
\usepackage{bbm}
\usepackage{amsmath}
\usepackage{amssymb}
\usepackage[makeroom]{cancel}

\DeclareMathOperator*{\argmin}{arg\,min}
\hypersetup{
    colorlinks=false,
    linkcolor=black,
    filecolor=magenta,
    urlcolor=cyan,
    }

\usepackage{booktabs}
\newcommand{\ra}[1]{\renewcommand{\arraystretch}{#1}}

\newtheorem{theorem}{Theorem}

\newtheorem{definition}{Definition}
\newtheorem{example}{Example}

\newtheorem{lemma}{Lemma}

\newtheorem{remark}{Remark}

\numberwithin{equation}{section}
\def\c#1{\mathcal{#1}}
\def\b#1{\mathbbm{#1}}
\def\bf#1{\mathbf{#1}}

\def\mat#1#2{\text{mat}_{#1}({#2})}

\def\wt#1{\widetilde{#1}}
\def\wh#1{\widehat{#1}}
\def\vec#1{\textnormal{\textbf{vec}}\big(#1\big)}

\def\cA{\mathcal{A}}
\def\cY{\mathcal{Y}}

\def\cC{\mathcal{C}}
\def\cF{\mathcal{F}}
\def\cE{\mathcal{E}}
\def\cK{\mathcal{K}}
\def\cX{\mathcal{X}}
\def\cN{\mathcal{N}}

\def\A{\mathbf{A}}
\def\B{\mathbf{B}}

\def\D{\mathbf{D}}
\def\E{\mathbf{E}}
\def\F{\mathbf{F}}

\def\H{\mathbf{H}}
\def\I{\mathbf{I}}
\def\Q{\mathbf{Q}}

\def\Y{\mathbf{Y}}
\def\M{\mathbf{M}}

\def\R{\mathbf{R}}

\def\X{\mathbf{X}}
\def\Z{\mathbf{Z}}
\def\0{\mathbf{0}}
\def\1{\mathbf{1}}

\def\diag{\textnormal{diag}}
\def\Fold{\textnormal{Fold}}
\def\Reshape{\textnormal{Reshape}}
\def\Amk{\A_{\text{-}k}}
\def\Qmk{\Q_{\text{-}k}}

\def\bSigma{\boldsymbol{\Sigma}}

\def\bgamma{\boldsymbol{\gamma}}
\def\bGamma{\boldsymbol{\Gamma}}

\def\bepsilon{\boldsymbol{\epsilon}}
\def\tr{\textnormal{tr}}

\def\dmk{d_{\text{-}k}}
\def\rmk{r_{\text{-}k}}
\def\Amk{\A_{\text{-}k}}

\def\a{\bf{a}}

\newcommand{\norm}[1]{\big\| #1 \big\|}

\allowdisplaybreaks[1]

\begin{document}
\setlength{\parindent}{18pt}
	\doublespacing
	\begin{titlepage}
		
		\title{On Testing Kronecker Product Structure in Tensor Factor Models}

		    \author{Zetai Cen\thanks{Zetai Cen is PhD student, Department of Statistics, London School of Economics. Email: Z.Cen@lse.ac.uk}}
            \author{Clifford Lam\thanks{Clifford Lam is Professor, Department of Statistics, London School of Economics. Email: C.Lam2@lse.ac.uk}}
		\affil{Department of Statistics, London School of Economics and Political Science}
		
		\date{}
		
		\maketitle

\begin{abstract}
We propose a test for testing the Kronecker product structure of a factor loading matrix implied by a tensor factor model with Tucker decomposition in the common component. Through defining a Kronecker product structure set, we define if a tensor time series response $\{\cY_t\}$ has a Kronecker product structure, equivalent to the ability to decompose $\{\cY_t\}$ according to a tensor factor model. Our test is built on analysing and comparing the residuals from fitting a full tensor factor model, and the residuals from fitting a (tensor) factor model on a reshaped version of the data. In the most extreme case, the reshaping is the vectorisation of the tensor data, and the factor loading matrix in such a case can be general if there is no Kronecker product structure present. Theoretical results are developed through asymptotic normality results on estimated residuals. Numerical experiments suggest that the size of the tests gets closer to the pre-set nominal value as the sample size or the order of the tensor gets larger, while the power increases with mode dimensions and the number of combined modes. We demonstrate out tests through a NYC taxi traffic data and a Fama-French matrix portfolio of returns.
\end{abstract}
		
		\bigskip
		\bigskip


\noindent

	\end{titlepage}
	
	\setcounter{page}{2}

\maketitle


\newpage
\section{Introduction}\label{sec: introduction}
With rapid advance in information technology, high dimensional time series data observed in tensor form are becoming more readily available for analysis in fields such as finance, economics, bioinformatics or computer science, to name but a few areas. In many cases, low-rank structures in the tensor time series observed can be exploited, facilitating analysis and interpretations. The most commonly used devices are the CP-decomposition and the multilinear/Tucker decomposition of a tensor, leading to CP-tensor factor models \citep{Changetal2023, Hanetal2024} and Tucker-decomposition tensor factor models \citep{Chenetal2022, Hanetal2024AOS, Chenetal2024,Chen_Lam} for tensor time series factor models, respectively. 
While tensor time series can be transformed back to vector time series through vectorisation and be analysed using traditional factor models for vector time series, the tensor structure of the data is lost and hence any corresponding interpretations from it. Moreover, vectorisation increases the dimension of the factor loading matrix significantly relative to the sample size, leading potentially to less accurate estimation and inferences \citep{Chen_Lam}.

However, a tensor factor model comes with its assumptions. For using the Tucker decomposition in particular, a tensor factor model assumes that the factor loading matrix for the vectorised data is the Kronecker product of smaller dimensional factor loading matrices. For instance, suppose at each $t\in[T]$, a mean-zero matrix $\Y_t\in \b{R}^{d_1\times d_2}$ is observed. Consider a matrix factor model (first studied in \cite{Wangetal2019}, and further extended/analysed in \cite{Yuetal2022a} or \cite{ChenFan2023} for example) of the form
\begin{equation}
\label{eqn: matrix_fm}
\Y_t = \A_1 \F_t \A_2' + \E_t,
\end{equation}
where $\F_t\in\b{R}^{r_1\times r_2}$ is the core factor, $\A_k\in\b{R}^{d_k\times r_k}$ is the mode-$k$ factor loading matrix, i.e., $\A_1$ and $\A_2$ are respectively the row and column loading matrices, and $\E_t$ is the noise. The vectorisation of \eqref{eqn: matrix_fm} is
\begin{equation}
\label{eqn: matrix_fm_vec}
\vec{\Y_t} = (\A_2 \otimes \A_1) \,\vec{\F_t} + \vec{\E_t}
\equiv \A_V \,\vec{\F_t} + \vec{\E_t},
\end{equation}
where $\A_V := \A_2 \otimes \A_1$, which is a vector factor model for the time series data $\{\vec{\Y_t}\}$ with factor loading matrix $\A_V$. Clearly, the implicit assumption of a Kronecker product structure for $\A_V$ when using a matrix factor model for matrix-valued time series data should be the first thing to check before such a factor model is applied.

Motivated by this simple example, we propose a test in this paper to test the Kronecker product structure of the factor loading matrix implied in the vectorised data when using a Tucker-decomposition tensor factor model (TFM), and extend it to higher order tensors. \cite{Heetal2023} has also noted this implicit assumption in a Tucker-decomposition matrix factor model, and proposes to test the ``boundary'' cases of each column (\text{resp.} row) of the data following a factor model with a common factor loading matrix, but with possibly distinct factors, or the whole matrix is just pure noise. Model \eqref{eqn: matrix_fm_vec} with a general $\A_V$ also implies a vector factor model with potentially different factor loading matrices for each column (\text{resp.} row) of the data, but they share the same factors. To explore the data as a matrix, connectedness through having a set of shared common factors rather than having the same factor loading matrix with all distinct factors is more meaningful. Practically, \eqref{eqn: matrix_fm_vec} is an alternative model easier to be satisfied by data than the ``boundary'' cases in  \cite{Heetal2023}, since the data still follows a more general factor model, just the implied Kronecker product structure in the factor loading matrix $\A_V$ is lost. This comes as no surprise then, that in all of the tests in \cite{Heetal2023} for their real data analyses, they cannot reject the null hypothesis of a matrix factor model. An easier alternative such as (\ref{eqn: matrix_fm_vec}) with just a general $\A_V$ can provide a more critical test for the null hypothesis of a matrix factor model. See our portfolio return example in Section \ref{subsec: realdataanalysis} for cases where our test can reject the null hypothesis of a matrix factor model, when \cite{Heetal2023} cannot.

We also stress that our model is fundamentally different from those used in testing for Kronecker product structure in the covariance matrix of the data. For example, \cite{Yuetal2022} and \cite{Guggenbergeretal2023} both propose tests for the Kronecker product structure of the covariance matrix of a vectorised matrix data. For model \eqref{eqn: matrix_fm}, even in the simplest hypothetical case of $\E_t$ and $\F_t$ being independent and $\F_t$ contains independent standard normal random variables, we have
\[\text{Cov}(\vec{\Y_t}) = \A_2\A_2'\otimes\A_1\A_1' + \text{Cov}(\vec{\E_t}),\]
so that the covariance matrix is never exactly of Kronecker product structure because of $\E_t$. Moreover, even with $\E_t = \0$, both $\A_1\A_1'$ and $\A_2\A_2'$ are of low rank, which is different from the full rank component matrices in the two papers mentioned above.

Our contributions in this paper are threefold. Firstly, as a first in the literature, we propose a test to test directly a Tucker-decomposition TFM against the alternative of a (tensor) factor model with Kronecker product structure lost in some of its factor loading matrices. As shown in Section~\ref{sec: model_identification}, for higher order tensors, testing against a tensor-decomposition TFM can be against a tensor factor model for the reshaped data, but not necessarily the vectorised data. This gives rise to flexibility and in fact statistical power in practical situations. Secondly, our analysis allows for weak factors, with our theoretical results developed to spell out rates of convergence explicitly. Last but not least, as a useful by-product, we developed tensor reshape theorems which can be useful in their own rights.

The rest of the paper is organised as follows. Section~\ref{sec: notations_tenreshape} introduces the notations used throughout this paper, and define the tensor reshape operation used for our tests. Section~\ref{sec: model_identification} introduces the Kronecker product structure set and pinpoints exactly through a theorem when a tensor time series $\{\cY_t\}$ follows a Tucker-decomposition TFM. This becomes the basis for the construction of our test statistics. Section~\ref{sec: Assumptions} lays down all the assumptions of the paper, and presents the main theoretical results for our test statistics to be valid. Section~\ref{sec: empirical} presents our simulation results and two sets of real data analyses. Finally, Section~\ref{sec: Appendix} provides links to the supplement of this paper and to the instructions of the R package \texttt{KOFM} for this paper.

\section{Notations and Tensor Reshape}\label{sec: notations_tenreshape}

\subsection{Notations}\label{subsec: notations}

Throughout this paper, we use the lower-case letter, bold lower-case letter, bold capital letter, and calligraphic letter, i.e., $x,\bf{x},\bf{X},\c{X}$, to denote a scalar, a vector, a matrix, and a tensor respectively.
We also use $x_i, X_{ij}, \bf{X}_{i\cdot}, \bf{X}_{\cdot i}$ to denote, respectively, the $i$-th element of $\bf{x}$, the $(i,j)$-th element of $\bf{X}$, the $i$-th row vector (as a column vector) of $\bf{X}$, and the $i$-th column vector of $\bf{X}$. We denote a column vector of $1$'s with length $a$ by $\1_a$. We use $\otimes$ to represent the Kronecker product, and $\circ$ the Hadamard product. By convention, the total Kronecker product for an index set is computed in descending index. We use $a\asymp b$ to denote $a=O(b)$ and $b=O(a)$. Hereafter, given a positive integer $m$, define $[m]:=\{1,2,\dots,m\}$. The $i$-th largest eigenvalue of a matrix $\X$ is denoted by $\lambda_i(\bf{X})$. The notation $\bf{X}\succcurlyeq 0$ (\text{resp.} $\bf{X} \succ 0$) means that $\bf{X}$ is positive semi-definite (\text{resp.} positive definite). We denote by $\tr(\X)$ the trace of $\X$, $\bf{X}'$ the transpose of $\bf{X}$, and $\diag(\X)$ a diagonal matrix with the diagonal elements of $\X$, while $\diag(\{x_1, \dots, x_n\})$ represents the diagonal matrix with $\{x_1, \dots, x_n\}$ on the diagonal. Define $d:=\prod_{k=1}^K d_k$, $\dmk :=d/d_k$, $r:= \prod_{k=1}^K r_k$ and $\rmk := r/r_k$.

\textbf{Norm notations}:
Sets are also denoted by calligraphic letters. For a given set $\cA$, we denote $|\cA|$ and $\cA_i$ as its cardinality and the $i$-th element respectively. We use $\|\bf{\cdot}\|$ to denote the spectral norm of a matrix or the $L_2$ norm of a vector, and $\|\bf{\cdot}\|_F$ to denote the Frobenius norm of a matrix. We use $\|\cdot\|_{\max}$ to denote the maximum absolute value of the elements in a vector, a matrix or a tensor. The notations $\|\cdot\|_1$ and $\|\cdot\|_{\infty}$ denote the $L_1$ and $L_{\infty}$-norm of a matrix respectively, defined by $\|\X\|_{1} := \max_{j}\sum_{i}|X_{ij}|$ and $\|\X\|_{\infty} := \max_{i}\sum_{j}|X_{ij}|$. Without loss of generality, we always assume the eigenvalues of a matrix are arranged by descending orders, and so are their corresponding eigenvectors.

\textbf{Tensor-related notations}: For the rest of this section, we briefly introduce the notations and operations for tensor data, which will be sufficient for this paper. For more details on tensor manipulations, readers are referred to \cite{KolderBader2009}. A multidimensional array with $K$ dimensions is an \textit{order}-$K$ tensor, with its $k$-th dimension termed as \textit{mode}-$k$. For an order-$K$ tensor $\c{X} = (X_{i_1,\ldots,i_K}) \in \mathbb{R}^{I_1\times\cdots\times I_K}$, a column vector $(X_{i_1, \ldots, i_{k-1}, i, i_{k+1}, \ldots, i_K})_{i\in[I_k]}$ represents a \textit{mode-$k$ fibre} for the tensor $\c{X}$. We denote by $\mat{k}{\c{X}} \in \mathbb{R}^{I_k\times I_{\text{-}k}}$ (or sometimes $\bf{X}_{(k)}$ for convenience, with $I_{\text{-}k} := (\prod_{j=1}^K I_j)/I_k$) the \textit{mode-$k$ unfolding/matricisation} of a tensor, defined by placing all mode-$k$ fibres into a matrix. We denote by $\c{X} \times_k \A$ the \textit{mode-$k$ product} of a tensor $\c{X}$ with a matrix $\A$, defined by
\[\mat{k}{\c{X}\times_k \A} := \A\, \mat{k}{\c{X}}.\]
We use the notation $\vec{\cdot}$ to denote the vectorisation of a matrix or the vectorisation of the mode-1 unfolding of a tensor. The \textit{refolding/tensorisation} of a vector $\a\in \b{R}^{I_1\ldots I_K}$ on $\{I_1,\dots, I_K\}$ is defined to be an order-$K$ tensor $\Fold\big(\a, \{I_1,\dots, I_K\}\big) \in\b{R}^{I_1\times \dots\times I_K}$ such that $\a =\vec{\Fold \big(\a, \{I_1,\dots, I_K\}\big)}$. The \textit{refolding/tensorisation} of a matrix $\A\in\b{R}^{I_k\times I_{\text{-}k}}$ on $\{I_1,\dots, I_K\}$ along mode-$k$ is defined to be $\Fold_k \big(\A, \{I_1,\dots, I_K\}\big) \in\b{R}^{I_1\times \dots\times I_K}$ such that $\A= \text{mat}_k\big(\Fold_k \big(\A, \{I_1,\dots, I_K\}\big)\big)$.

\subsection{Introduction to tensor reshape}\label{subsec: reshape}

In this subsection, we introduce tensor reshape. Given an order-$K$ tensor $\cX\in\b{R}^{I_1\times \dots \times I_K}$ and a set with ordered, strictly ascending elements $\{a_1, \dots, a_\ell\} \subseteq [K]$, the Reshape$(\cdot,\cdot)$ operator is defined as follows:
\begin{align*}
    \text{If } \ell=1, \;
    &\Reshape\big( \cX , \{a_1\} \big) := \Fold_K\big(\mat{a_1}{\cX}, \{I_1,\dots, I_{a_1-1}, I_{a_1+1}, \dots, I_K, I_{a_1}\}\big); \\
    \text{if } \ell=2, \;
    &\Reshape\big( \cX , \{a_1, a_2\} \big) := \Fold_{K-1}\big( \cX_{a_1\sim a_2}, \{I_1,\dots,I_{a_1-1}, I_{a_1+1},\dots, I_{a_2-1}, I_{a_2+1}, \dots, I_{K}, I_{a_1}I_{a_2}\} \big), \\
    &\text{where }
    \cX_{a_1\sim a_2} := \begin{pmatrix}
    \text{mat}_{a_1}\big\{\Fold\big( \mat{a_2}{\cX}_{1\cdot} , \{I_1, \dots, I_{a_2-1}, I_{a_2+1}, \dots, I_K\}\big) \big\} \\
    \ldots \\
    \text{mat}_{a_1}\big\{\Fold\big( \mat{a_2}{\cX}_{I_{a_2}\cdot} , \{I_1, \dots, I_{a_2-1}, I_{a_2+1}, \dots, I_K\}\big) \big\}
    \end{pmatrix} ; \\
    \text{if } \ell \geq 3, \;
    &\Reshape\big( \cX , \{a_1, \dots, a_\ell\} \big) :=
    \Reshape\big[ \Reshape\big( \cX , \{a_{\ell-1}, a_{\ell}\} \big), \{a_1, \dots, a_{\ell-2}, K-1\} \big].
\end{align*}
Hence, reshaping an order-$K$ tensor along $\{a_1, \dots, a_\ell\}$ results in an order-$(K-(\ell-1))$ tensor. A heuristic view of $\Reshape\big( \cX , \{a_1, \dots, a_\ell\} \big)$ is that all modes of $\cX$ with indices $\{a_1, \dots, a_\ell\}$ are ``merged'' into a single mode acting as the last mode as a result. Note that one may recover $\cX$ from $\Reshape\big( \cX , \{a_1, \dots, a_\ell\} \big)$ given the original dimension of $\cX$ and $\{a_1,\dots, a_\ell\}$.

As a simple example on tensor reshape, consider a matrix $\X\in\b{R}^{I_1 \times I_2}$. Trivially, we have
\[
\Reshape\big( \X , \{2\} \big) =\X, \;\;\;
\Reshape\big( \X , \{1\} \big) =\X'.
\]
Moreover, $\Reshape\big( \X , \{1,2\} \big) = \Fold_1(\X_{1\sim 2},\{I_1 I_2\}) =\X_{1\sim 2} =\vec{\X}$ since
\begin{align*}
    \X_{1\sim 2} =\begin{pmatrix}
    \text{mat}_{1}\big\{\Fold\big( \mat{2}{\X}_{1\cdot} , \{I_1\}\big) \big\} \\
    \ldots \\
    \text{mat}_{1}\big\{\Fold\big( \mat{2}{\X}_{I_2\cdot} , \{I_1\}\big) \big\}
    \end{pmatrix}
    =\begin{pmatrix}
    (\X')_{1\cdot} \\
    \ldots \\
    (\X')_{I_2\cdot}
    \end{pmatrix} = \vec{\X}.
\end{align*}
In fact, it holds for any order-$K$ tensor $\cX$ that $\Reshape\big( \cX , [K]\big) =\vec{\cX}$.

\textbf{Some useful algebra of tensor reshape}:
The reshape operator is linear in the first argument, i.e.,
\[
\Reshape\big( b_1\cX_1 + b_2\cX_2 , \{a_1, \dots, a_\ell\} \big) = b_1 \cdot \Reshape\big( \cX_1 , \{a_1, \dots, a_\ell\} \big) + b_2 \cdot \Reshape\big( \cX_2 , \{a_1, \dots, a_\ell\} \big) .
\]
Moreover, for two sets $\{a_1,\dots, a_\ell\}$, $\{b_1,\dots, b_g\}$ such that $a_\ell < b_1$ (i.e., all elements in the first set are less than those in the second), it holds that
\[
\Reshape\big( \cX , \{a_1, \dots, a_\ell, b_1,\dots, b_g\} \big)
= \Reshape\big[ \Reshape\big( \cX , \{b_1,\dots, b_g\} \big), \{a_1,\dots, a_\ell, K-g+1\} \big] .
\]
Note that $\{a_1,\dots, a_\ell, K-g+1\}$ is indeed strictly ascending since
\[
a_\ell \leq b_1-1 \leq b_g- (g-1) -1 = b_g -g\leq K-g .
\]

\section{Factor Model and Testing its Kronecker Product Structure}\label{sec: model_identification}

Section~\ref{subsec: model_KPS} introduces the concept of factor models with Kronecker product structure and lays down the technical details for the testing problem. For an integral reading experience, readers can go straight to Section~\ref{subsec: test_KPS} where equations and terms can be referred back to Section~\ref{subsec: model_KPS} whenever necessary.

\subsection{Factor model and Kronecker product structure}\label{subsec: model_KPS}

We begin by introducing the Kronecker product structure set which facilitates description of our models.

\begin{definition}\label{def: kron_structure_set}
\textnormal{(Kronecker product structure set)}
Given an ordered set of positive integers $\{b_1, \dots, b_{\kappa}\}$, the Kronecker product structure set is a set of full column rank matrices defined as
\[
\cK_{b_1\times \dots \times b_\kappa} := \big\{ \A \,|\, \A= \A_\kappa \otimes\dots \otimes \A_1 \text{ with } \A_j \in\b{R}^{b_j\times u_j} \text{ of rank } u_j\ll b_j,\,  \|\A_{j,\cdot i}\|^2 \asymp b_j^{\delta_{j,i}},\, \delta_{j,i} \in (0,1] \big\} .
\]
\end{definition}

The Kronecker product structure set defined by Definition~\ref{def: kron_structure_set} characterises the factor loading matrix, and requiring $\delta_{j,i}>0$ is to ensure certain factor strength in each loading matrix. See Assumptions (L1) and (L2) in Section~\ref{subsec: assumption} for the technical details. The form of factor models is depicted below, with the feature of Kronecker product structure.

\begin{definition}\label{def: kron_structure}
\textnormal{(Factor models and Kronecker product structure)}
Given a series of mean-zero order-$K$ tensors $\cY_t\in \b{R}^{d_1\times \dots\times d_K}$ for $t\in[T]$ and a set with ordered, ascending elements $\cA= \{a_1, \dots, a_\ell\} \subseteq[K]$, we say $\{\cY_t\}$ follows a factor model along $\cA$ if for $t\in[T]$,
\begin{equation}
\label{eqn: fm_wo_Kron}
\begin{split}
\Reshape(\cY_t, \cA) &= \cC_{\textnormal{reshape},t} + \cE_{\textnormal{reshape},t} =
\cF_{\textnormal{reshape},t} \times_{j=1}^{K-\ell+1} \A_{\textnormal{reshape},j} + \cE_{\textnormal{reshape},t},
\end{split}
\end{equation}
where $\Reshape(\cY_t, \cA) \in\b{R}^{I_1\times \dots \times I_{K-\ell+1}}$ is the order-$(K-\ell+1)$ tensor by reshaping $\cY_t$ along $\cA$, the common component $\cC_{\textnormal{reshape},t}$ consists of the core factor $\cF_{\textnormal{reshape},t} \in\b{R}^{\pi_1\times \dots\times \pi_{K-\ell+1}}$ and loading matrices $\A_{\textnormal{reshape},j} \in \b{R}^{I_j\times \pi_j}$ with rank $\pi_j\ll I_j$ for $j\in[K-\ell+1]$, and $\cE_{\textnormal{reshape},t}$ is the noise. We further make the following classifications.
\begin{itemize}
    \item [1.] $\{\cY_t\}$ has a Kronecker product structure if $\A_{\textnormal{reshape}, K-\ell+1} \in \cK_{d_{a_1} \times \dots \times d_{a_\ell}}$;
    \item [2.] $\{\cY_t\}$ has no Kronecker product structure along $\cA$ if $\A_{\textnormal{reshape}, K-\ell+1} \notin \cK_{d_{a_1} \times \dots \times d_{a_\ell}}$.
\end{itemize}
\end{definition}

Definition~\ref{def: kron_structure} formally defines the form of factor models considered in this paper. A key information lying in Definition~\ref{def: kron_structure}.1 is that if the Kronecker product structure holds along some $\cA$, the structure holds along any $\cA$; see the discussion below Theorem~\ref{thm: reshape} for details. Note that if $\ell=1$ in Definition~\ref{def: kron_structure}, i.e., $\cA$ contains only one element (representing the mode index), for each order-$K$ tensor $\cY_t$, $\Reshape( \cY_t , \{a_1\})$ is the order-$K$ tensor constructed from $\cY_t$ by treating mode-$a_1$ as mode-$K$. Hence, the factor model of $\cY_t$ along $\{a_1\}$ returns to a Tucker-decomposition TFM \citep{Chenetal2022, Heetal2022a} of $\cY_t$ with mode indices changed. For instance, we may read \eqref{eqn: fm_wo_Kron} along $\cA=\{K\}$ as
\[
\cY_t = \cC_{\textnormal{reshape},t} + \cE_{\textnormal{reshape},t} = \cF_{\textnormal{reshape},t} \times_1 \A_{\textnormal{reshape},1} \times_2 \dots \times_K \A_{\textnormal{reshape},K} + \cE_{\textnormal{reshape},t}.
\]
Hence Definition~\ref{def: kron_structure}.1 automatically describes $\{\cY_t\}$ if $\ell=1$, implying that Kronecker product structure is only non-trivial for $\ell\geq 2$ (hence $K\geq 2$).
To demystify Definition~\ref{def: kron_structure}.1, we next present Theorem~\ref{thm: reshape} which, as a first in the literature, spells out the equivalence of Tucker-decomposition TFM under tensor reshape.

\begin{theorem}\label{thm: reshape}
\textnormal{(Tensor Reshape Theorem I)}
With the notations in Definition~\ref{def: kron_structure}, $\{\cY_t\}$ following \eqref{eqn: fm_wo_Kron} along any given $\cA = \{a_1,\ldots,a_\ell\} \subseteq [K]$ with a Kronecker product structure is equivalent to $\{\cY_t\}$ following a Tucker-decomposition factor model such that
\begin{equation}
\label{eqn: fm_w_Kron}
\cY_t = \cC_t + \cE_t = \cF_t \times_1 \A_1 \times_2 \dots \times_K \A_K + \cE_t,
\end{equation}
where $\cC_t$ is the common component, $\cF_t\in \b{R}^{r_1\times \dots \times r_K}$ is the core factor, each $\A_k\in \b{R}^{d_k\times r_k}$ with $r_k\ll d_k$ is the mode-$k$ loading matrix, and $\cE_t$ is the noise. More importantly, with $\cA^\ast := [K]\setminus \cA$,
\begin{align*}
    &\cF_{\textnormal{reshape},t} = \Reshape(\cF_t, \cA), \;\;\;
    \cE_{\textnormal{reshape},t} = \Reshape(\cE_t, \cA), \\
    &\A_{\textnormal{reshape},K-\ell+1} = \otimes_{i\in \cA} \A_i, \;\;\;
    \A_{\textnormal{reshape},j} = \A_{\cA_j^\ast}
    \;\;\; \text{for } j\in[K-\ell].
\end{align*}
Moreover, \eqref{eqn: fm_w_Kron} uniquely determines \eqref{eqn: fm_wo_Kron}, and \eqref{eqn: fm_wo_Kron} determines \eqref{eqn: fm_w_Kron} up to an arbitrary set $\{\A_i\}_{i\in\cA}$.
\end{theorem}

Theorem \ref{thm: reshape} reveals that a factor model on $\{\cY_t\}$ with Kronecker product structure in Definition~\ref{def: kron_structure} is in fact a Tucker-decomposition TFM on $\{\cY_t\}$. We point out that Theorem~\ref{thm: reshape} enables us to reshape any order-$K$ tensor $\cY_t$ (under Tucker-decomposition TFM) along some selected $\cA$, so that $\{\Reshape(\cY_t, \cA)\}$ still follows some Tucker-decomposition TFM. This salvages the practical issue of insufficient mode dimensions and potentially boosts performances of tests and estimations, while interpretability is not necessarily lost much; see Remark~\ref{remark: reshape_thm}. The identification of \eqref{eqn: fm_wo_Kron} and \eqref{eqn: fm_w_Kron} are discussed in the supplement.

\begin{remark}\label{remark: other_decomp}
Both \eqref{eqn: fm_wo_Kron} and \eqref{eqn: fm_w_Kron} are based on a Tucker decomposition for the observed tensor. Other tensor decompositions are possible, such as the CP decomposition \citep{KolderBader2009} and the Low Separation Rank (LSR) decomposition \citep{Takietal2024}, etc. As CP decomposition is a special Tucker decomposition, our defined factor model is more general. The LSR decomposition is generalised further from Tucker decomposition, but the structure is less helpful here and brings in unnecessary complication due to the arbitrary separation rank.
\end{remark}

\begin{remark}\label{remark: reshape_thm}
We discuss here the practical significance of Theorem~\ref{thm: reshape}, which can be of independent interest. As mentioned in Section~\ref{subsec: reshape}, modes selected by $\cA$ are merged into a single one from $\cY_t$ to $\Reshape(\cY_t, \cA)$. Hence if $\{\cY_t\}$ is deemed to follow a Tucker-decomposition TFM, we may choose $\cA$ to include modes with small dimensions relative to others' and directly work on $\{\Reshape(\cY_t, \cA)\}$ which follows a Tucker-decomposition TFM guaranteed by Theorem~\ref{thm: reshape}. Consider a simple example that for $t\in[T]$, $\cY_t\in \b{R}^{d_1\times d_2 \times d_3}$ satisfies \eqref{eqn: fm_w_Kron} with all factors pervasive. If $T\asymp d_1 \asymp d_2^2\asymp d_3^2$, i.e., the observed tensor is unbalanced, then the squared error of each (PCA-based) estimated common component on $\{\cY_t\}$ is $O_P(1/T)$, that on $\{\textnormal{\textbf{vec}}(\cY_t)\}$ is $O_P(1/T)$ and that on $\{ \Reshape(\cY_t, \{2,3\} )\}$ is $O_P(1/T^2)$.
\end{remark}

\subsection{A test on Kronecker product structure}\label{subsec: test_KPS}

The testing problem on Kronecker product structure is formally defined in the following, with an example on an order-2 tensor (i.e., a matrix) given.

For each $t\in[T]$, we observe a mean-zero order-$K$ tensor $\cY_t\in \b{R}^{d_1\times \dots\times d_K}$ with $K\geq 2$ (otherwise the test is trivial as explained in Section~\ref{subsec: model_KPS}). Without loss of generality, let $v<K$ be a given positive integer and denote $\cA =\{v,v+1,\dots, K-1,K\}$ which contains the mode indices along which the Kronecker product structure might be lost; see the alternative hypothesis $H_1$ below. Suppose $\{\cY_t\}$ follows a factor model along $\cA$ as in Definition~\ref{def: kron_structure}, with notations therein except that we now read \eqref{eqn: fm_wo_Kron} as
\begin{equation}
\label{eqn: test_fm_wo_Kron}
\begin{split}
\Reshape(\cY_t, \cA) &= \cC_{\textnormal{reshape},t} + \cE_{\textnormal{reshape},t} =
\cF_{\textnormal{reshape},t} \times_1 \A_1 \times_2 \dots \times_{v-1} \A_{v-1} \times_v \A_V + \cE_{\textnormal{reshape},t},
\end{split}
\end{equation}
where $\A_j\in \b{R}^{d_j\times r_j}$ for $j\in[v-1]$ (if $v> 1$) and $\A_V\in \b{R}^{d_V \times r_V}$ with $d_V:= \prod_{i=v}^K d_i$. Essentially, the order-$v$ tensor $\Reshape(\cY_t, \cA)$ follows a Tucker-decomposition TFM. The set $\{r_1, \dots, r_{v-1}, r_V\}$ is assumed known and any consistent estimators \citep[e.g.][]{Hanetal2022, ChenLam2024a} can be used in practice. With $\cK_{d_v\times \dots \times d_K}$ defined in Definition~\ref{def: kron_structure_set}, we consider a hypothesis test as follows:
\begin{equation}
\label{eqn: test_kron_structure}
\begin{split}
    & H_0: \text{$\{\cY_t\}$ has a Kronecker product structure, i.e., $\A_V \in \cK_{d_v\times \dots \times d_K}$} ; \\
    & H_1: \text{$\{\cY_t\}$ has no Kronecker product structure along $\cA$, i.e., $\A_V \notin \cK_{d_v\times \dots \times d_K}$} .
\end{split}
\end{equation}

Besides the complexity of being a composite testing problem, the difficulty of \eqref{eqn: test_kron_structure} is elevated by the fact that $\cY_t$ under the alternative has no explicit form without reshaping along $\cA$. Fortunately, the factor structure in \eqref{eqn: test_fm_wo_Kron} is stable under both hypotheses. That is, the estimation of $\{\cF_{\textnormal{reshape},t}, \A_1, \dots, \A_{v-1}, \A_V\}$ is always feasible. In particular, thanks to Theorem~\ref{thm: reshape}, we have the following under $H_0$:
\begin{equation}
\label{eqn: test_fm_w_Kron}
\cY_t = \cC_t + \cE_t = \cF_t \times_1 \A_1 \times_2 \dots \times_{v-1} \A_{v-1} \times_v \A_v \times_{v+1} \dots \times_K \A_K + \cE_t,
\end{equation}
where $\A_k\in\b{R}^{d_k\times r_k}$ for $k\in[K]$ (hence the first $v-1$ loading matrices are exactly those in \eqref{eqn: test_fm_wo_Kron}), and that
\begin{align*}
    \Reshape(\cF_t, \cA) = \cF_{\textnormal{reshape},t}, \;\;\;
    \Reshape(\cE_t, \cA) = \cE_{\textnormal{reshape},t}, \;\;\;
    \A_K \otimes \A_{K-1} \otimes \dots \otimes \A_v = \A_V .
\end{align*}

\begin{example}\label{example: test_order2}
Let $\Y_t\in \b{R}^{d_1\times d_2}$ be a matrix-valued observation at $t\in[T]$. For the setup, we can only specify $\cA=\{1,2\}$ (which is the only non-trivial case here as discussed in Section~\ref{subsec: model_KPS}). The hypothesis test \eqref{eqn: test_kron_structure} is simplified as follows, with $\cA$ reflected by the vectorisation:
\begin{align*}
    & H_0: \Y_t = \A_1 \F_t \A_2' + \E_t ; \\
    & H_1: \vec{\Y_t} = \A_V \vec{\F_t} + \vec{\E_t}, \;\;\;
    \text{with $\A_V \notin \cK_{d_1\times d_2}$} .
\end{align*}
\end{example}

\subsection{Constructing the test statistic}\label{subsec: construct_test}

Despite the obscure $\cK_{d_v\times \dots \times d_K}$ in \eqref{eqn: test_kron_structure}, we may resort to the Tucker-decomposition TFM in \eqref{eqn: test_fm_w_Kron} under $H_0$. To construct the test, we first obtain estimators for the (standardised) loading matrices in \eqref{eqn: test_fm_wo_Kron}. For $j\in[v-1]$, $\wt\Q_{j}$ is defined as the eigenvector matrix corresponding to the $r_j$ largest eigenvalues of
\[
\frac{1}{T} \sum_{t=1}^T \Reshape(\cY_t, \cA)_{(j)} \,\Reshape(\cY_t, \cA)_{(j)}' ,
\]
where $\Reshape(\cY_t, \cA)_{(j)}$ represents the mode-$j$ unfolding matrix of $\Reshape(\cY_t, \cA)$ (see Section~\ref{subsec: notations}). Similarly, $\wt\Q_{V}$ is the eigenvector matrix corresponding to the $r_V$ largest eigenvalues of
\[
\frac{1}{T} \sum_{t=1}^T \Reshape(\cY_t, \cA)_{(v)} \,\Reshape(\cY_t, \cA)_{(v)}' .
\]
Then $\cC_{\textnormal{reshape},t}$ and $\cE_{\textnormal{reshape},t}$ are respectively estimated by
\begin{align}
    \wt\cC_{\textnormal{reshape},t} &:= \Reshape(\cY_t, \cA) \times_{j=1}^{v-1} (\wt\Q_j \wt\Q_j') \times_v (\wt\Q_V \wt\Q_V') ,
    \label{eqn: tilde_cCt} \\
    \wt\cE_{\textnormal{reshape},t} &:= \Reshape(\cY_t, \cA) - \wt\cC_{\textnormal{reshape},t} .
    \label{eqn: tilde_cEt}
\end{align}

For \eqref{eqn: test_fm_w_Kron}, $\wh\Q_{j}$ for $j\in[v-1]$ is defined as the eigenvector matrix corresponding to the $r_j$ largest eigenvalues of $T^{-1} \sum_{t=1}^T \Y_{t,(j)} \Y_{t,(j)}'$. Next, denote $\c{R}$ as the set of all divisor combinations of $r_V$, i.e.,
\begin{equation}
\label{eqn: set_R}
\c{R} := \Big\{ \big(\pi_1, \pi_2, \dots, \pi_{K-v+1} \big) \,\Big|\, \prod_{j=1}^{K-v+1} \pi_j=r_V \text{ and each }  \pi_j \in \b{Z}^{+} \text{ with } \pi_j \leq d_{j+v-1} \Big\}.
\end{equation}
Let the $m$-th element of $\c{R}$ be $(\pi_{m,1}, \dots, \pi_{m,K-v+1})$. Then for $i\in\{v, v+1,\dots, K\}$, we obtain $\wh\Q_{m,i}$ as the eigenvector matrix corresponding to the $\pi_{m,i-v+1}$ largest eigenvalues of $T^{-1} \sum_{t=1}^T \Y_{t,(i)} \Y_{t,(i)}'$. The common component and residual estimators are hence obtained as
\begin{align}
    \wh\cC_{m,t} &:= \cY_t \times_{j=1}^{v-1} (\wh\Q_j \wh\Q_j') \times_{i=v}^{K} (\wh\Q_{m,i} \wh\Q_{m,i}') ,
    \label{eqn: hat_cCt} \\
    \wh\cE_{m,t} &:= \cY_t - \wh\cC_{m,t} .
    \label{eqn: hat_cEt}
\end{align}

Let $\wt\cE_t$ be the order-$K$ tensor with the same dimension as $\cY_t$ such that $\Reshape(\wt\cE_t, \cA)= \wt\cE_{\textnormal{reshape},t}$. Define $k^\ast:= \argmin_{k\in[K]} \{d_k\}$ and denote the mode-$k^\ast$ unfolding of $\wt\cE_t$ and $\wh\cE_{m,t}$ as $\wt\E_{t,(k^\ast)}$ and $\wh\E_{m,t,(k^\ast)}$, respectively. Theorem~\ref{thm: noise_aggregate_dist} (in Section~\ref{subsec: main_theorem}) tells us that there exists $m\in[|\c{R}|]$ such that for each $t\in[T]$, $j\in [d/d_{k^\ast}]$, both
\begin{align}
    x_{j,t} := \frac{1}{d_{k^\ast}} \sum_{i=1}^{d_{k^\ast}} \wt{E}_{t,(k^\ast),ij}^2
    ,\;\;\;
    y_{m,j,t} := \frac{1}{d_{k^\ast}} \sum_{i=1}^{d_{k^\ast}} \wh{E}_{m,t,(k^\ast),ij}^2 ,
    \notag
\end{align}
are asymptotically distributed the same under $H_0$, and $x_{j,t}$ in particular is distributed the same under either $H_0$ or $H_1$. Let $\b{P}_{x,j}$ and $\b{P}_{y,m,j}$ respectively denote the empirical probability measures induced by the empirical cumulative distribution functions for $\{x_{j,t}\}_{t\in[T]}$ and $\{y_{m,j,t}\}_{t\in[T]}$:
\begin{equation}\label{eqn: CDF}
\begin{split}
&\b{F}_{x,j}(c) := \frac{1}{T} \sum_{t=1}^T\b{1}\{x_{j,t} \leq c\}, \;\;\;
\b{F}_{y,m,j}(c) := \frac{1}{T} \sum_{t=1}^T\b{1}\{y_{m,j,t} \leq c\}.
\end{split}
\end{equation}
Let $\wh{q}_{x,j}(\alpha) := \inf\big\{ c\;|\;\b{F}_{x,j}(c) \geq 1-\alpha \big\}$. The intuition here is that if $H_0$ is satisfied, then over different $j\in[d/d_{k^\ast}]$, the cumulative distribution functions $\b{F}_{x,j}(\cdot)$ and $\b{F}_{y,m,j}(\cdot)$ should be similar. However, if $H_1$ is true, then we expect the residuals in $\wh{\E}_{m,t,(k^\ast)}$ to be inflated, so that the empirical probability of $\{y_{m,j,t} \geq \wh{q}_{x,j}(\alpha)\}$ is expected to be larger than $\alpha$; see the theoretical statement on this in Theorem \ref{thm: test}. To incorporate this across different  $j\in[d/d_{k^\ast}]$, we compare the 5\% quantile of $T^{-1}\sum_{t=1}^T \b{1}\{y_{m,j,t} \geq \wh{q}_{x,j}(\alpha)\}$ over $j\in[d/d_{k^\ast}]$ to $\alpha$, and expect it to be larger than $\alpha$ under $H_1$. 

Since with the wrong number of factors, a particular $m\in[|\c{R}|]$ will in general inflate the residuals $y_{m,j,t}$ further,
 in practice, to be on the conservative side, we reject $H_0$ if
\begin{equation}
\label{eqn: reject_rule}
\min_{m\in [|\c{R}|]} \Bigg\{ \text{5\% quantile of } \frac{1}{T} \sum_{t=1}^T \b{1}\{y_{m,j,t} \geq \wh{q}_{x,j}(\alpha)\} \text{ over $ j\in [d/d_{k^\ast}]$} \Bigg\} >\alpha ,
\end{equation}
noting that exactly one element in $\c{R}$ represents the true number of factors on the modes with indices in $\cA$. We also point out that there are other possible ways to aggregate the information from each $j$, but \eqref{eqn: reject_rule} empirically works well and circumvents possible issues such as heavy-tailed noise, under- or over-estimation on the number of factors, and insufficient data dimensions; see Section~\ref{subsec: simulation}.

\begin{remark}\label{remark: divisor_combination_rV}
(Explanation of $\c{R}$ in \eqref{eqn: set_R})
It is possible to perform the test directly using the number of factors for modes in $\cA$, i.e., $\c{R}$ only contains the number of factors $r_j$, $j=v, \dots, K$, in \eqref{eqn: test_fm_w_Kron}. This is guaranteed by the Tucker-decomposition TFM under $H_0$ in \eqref{eqn: test_kron_structure}. However, usually in practice we need to estimate the number of factors which are invalid under $H_1$ in \eqref{eqn: test_kron_structure}. This leads to unstable estimated number of factors and hence unstable test statistic, which we address by introducing $\c{R}$ in \eqref{eqn: set_R}.
\end{remark}

\section{Assumptions and Theoretical Results}\label{sec: Assumptions}

\subsection{Assumptions}\label{subsec: assumption}

This section presents all the assumptions for testing $H_0$ against $H_1$ in \eqref{eqn: test_kron_structure}. Another version (with only different notations) of Assumptions (L1) and (L2) for the identification of \eqref{eqn: fm_wo_Kron} and \eqref{eqn: fm_w_Kron} is included in the supplement, with identification theorem presented and proved there.

\begin{itemize}
\item[(L1)]
\em
For each $j\in[v-1]$, we assume that $\A_{j}$ in \eqref{eqn: test_fm_wo_Kron} is of full rank and as $d_j\to\infty$,
\begin{equation}
\label{eqn: L1}
\Z_{j}^{-1/2} \A_{j}' \A_{j} \Z_{j}^{-1/2}
\to \bSigma_{A,j},
\end{equation}
where $\bSigma_{A,j}$ is positive definite with all eigenvalues bounded away from 0 and infinity, and $\Z_{j}$ is a diagonal matrix with $(\Z_{j})_{hh}\asymp I_j^{\delta_{j,h}}$ for $h\in[r_j]$ and the ordered factor strengths $1/2<\delta_{j,r_j}\leq \dots\leq \delta_{j,2}\leq \delta_{j,1}\leq 1$.

We assume that $\A_V$ also has the above form with $\Z_V$ and $\bSigma_{A,V}$, except that only the maximum and minimum factor strengths are ordered, i.e.,  $1/2<\delta_{V,r_V}\leq \delta_{V,h} \leq \delta_{V,1}\leq 1$ for any $h\in[r_V]$.
\end{itemize}

\begin{itemize}
\item[(L2)]
\em
With $\cA= \{v,v+1,\dots,K\}$, we assume that for each $i\in\cA$, $\A_i$ in \eqref{eqn: test_fm_w_Kron} is of full rank and as $d_i\to\infty$,
\begin{equation}
\label{eqn: L2}
\Z_i^{-1/2} \A_{i}' \A_i \Z_i^{-1/2}
\to \bSigma_{A,i},
\end{equation}
where $\bSigma_{A,i}$ is positive definite with all eigenvalues bounded away from 0 and infinity, and $\Z_i$ is a diagonal matrix with $(\Z_i)_{hh}\asymp d_i^{\delta_{i,h}}$ for $h\in[r_i]$ and the ordered factor strengths $1/2<\delta_{i,r_i}\leq \dots\leq \delta_{i,2}\leq \delta_{i,1}\leq 1$.
\end{itemize}

\begin{itemize}
\item[(F1)] (Time series in $\cF_{\textnormal{reshape},t}$)
\em
There is $\cX_{\textnormal{reshape},f,t}$ the same dimension as $\cF_{\textnormal{reshape},t}$ such that $\cF_{\textnormal{reshape},t} = \sum_{w\geq 0} a_{f,w} \cX_{\textnormal{reshape},f,t-w}$. The time series $\{ \cX_{\textnormal{reshape},f,t}\}$ has \text{i.i.d.} elements with mean $0$, variance $1$ and uniformly bounded fourth order moments. The coefficients $a_{f,w}$ satisfy $\sum_{w\geq 0} a_{f,w}^2=1$ and $\sum_{w\geq 0}|a_{f,w}| \leq c$ for some constant $c$.
\end{itemize}

\begin{itemize}
\item[(E1)] (Decomposition of $\c{E}_t$)
\em
The noise $\c{E}_t$ (such that $\cE_{\textnormal{reshape},t} = \Reshape(\cE_t, \cA)$) can be decomposed as
\begin{equation}
\label{eqn: E1}
    \c{E}_t = \c{F}_{e,t}\times_1 \A_{e,1}\times_2 \cdots \times_K\A_{e,K}
    + \bSigma_{\epsilon}\circ\bepsilon_t,
\end{equation}
where order-$K$ tensors $\c{F}_{e,t}\in \b{R}^{r_{e,1} \times \cdots \times r_{e,K}}$ and $\bepsilon_t\in \b{R}^{d_1 \times \cdots \times d_K}$ contain independent mean zero elements with unit variance, with the two time series $\{\bepsilon_t\}$ and $\{\c{F}_{e,t}\}$ being independent. The order-$K$ tensor $\bSigma_{\epsilon}$ contains the standard deviations of the corresponding elements in $\bepsilon_t$, and has elements uniformly bounded.

Moreover, for each $k\in[K]$, $\A_{e,k} \in \mathbb{R}^{d_k\times r_{e,k}}$ is (approximately) sparse such that $\norm{\A_{e,k}}_1=O(1)$.
\end{itemize}

\begin{itemize}
\item[(E2)] (Time series in $\cE_t$)
\em
There is $\c{X}_{e,t}$ the same dimension as $\c{F}_{e,t}$, and $\c{X}_{\epsilon,t}$ the same dimension as $\bm\epsilon_t$, such that $\c{F}_{e,t}=\sum_{q\geq 0}a_{e,q}\c{X}_{e,t-q}$ and $\bm\epsilon_t=\sum_{q\geq 0}a_{\epsilon,q}\c{X}_{\epsilon, t-q}$, with $\{\c{X}_{e,t}\}$ and $\{\c{X}_{\epsilon, t}\}$ independent of each other. $\{\c{X}_{e,t}\}$ has independent elements while $\{\c{X}_{\epsilon, t}\}$ has \text{i.i.d.} elements, and all elements have mean zero with unit variance and uniformly bounded fourth order moments. Both $\{\c{X}_{e,t}\}$ and $\{\c{X}_{\epsilon, t}\}$ are independent of $\{ \cX_{\textnormal{reshape},f,t}\}$ from (F1).

The coefficients $a_{e,q}$ and $a_{\epsilon, t}$ are such that $\sum_{q\geq 0} a_{e,q}^2=\sum_{q\geq 0} a_{\epsilon,q}^2=1$ and $\sum_{q\geq 0}|a_{e,q}|, \sum_{q\geq 0}|a_{\epsilon,q}|\leq c$ for some constant $c$.
\end{itemize}

\begin{itemize}
\item[(R1)] (Rate assumptions)
\em
With $g_s := \prod_{k=1}^K d_k^{\delta_{k,1}}$ and $\gamma_s := d_V^{\delta_{V,1}} \prod_{j=1}^{v-1} d_j^{\delta_{j,1}}$, we assume that
\begin{align*}
d g_s^{-2} T^{-1} d_k^{2(\delta_{k,1} - \delta_{k,r_k})+1} , \;\;\;
d g_s^{-1} d_k^{\delta_{k,1} -\delta_{k,r_k}-1/2} &, \\
d \gamma_s^{-2} T^{-1} d_V^{2(\delta_{V,1} - \delta_{V,r_V})+1} , \;\;\;
d \gamma_s^{-1} d_V^{\delta_{V,1} -\delta_{V,r_V}-1/2} &= o(1) .
\end{align*}
\end{itemize}

\begin{itemize}
\item[(R2)] (Further rate assumptions)
\em
With $g_w := \prod_{k=1}^K d_k^{\delta_{k,r_k}}$ and $\gamma_w := d_V^{\delta_{V,r_V}} \prod_{j=1}^{v-1} d_j^{\delta_{j,r_j}}$, we assume that
\begin{align*}
    \max_{k\in[K]}\Bigg\{ d_k^{2(\delta_{k,1} - \delta_{k,r_k})} \Bigg( \frac{1}{T\dmk d_k^{1- \delta_{k,1}}} +\frac{1}{ d_k^{1+ \delta_{k, r_k}}} \Bigg) \frac{d^2}{g_s g_w} \Bigg\} &, \\
    \max_{j\in[v-1]}\Bigg\{ d_j^{2(\delta_{j,1} - \delta_{j,r_j})} \Bigg( \frac{1}{T\dmk d_j^{1- \delta_{j,1}}} +\frac{1}{ d_j^{1+ \delta_{j, r_j}}} \Bigg) \frac{d^2}{\gamma_s \gamma_w} \Bigg\} &, \\
    d_V^{2(\delta_{V,1} - \delta_{V,r_V})} \Bigg( \frac{1}{Td d_V^{- \delta_{V,1}}} ,\;\;\;
    \frac{1}{ d_V^{1+ \delta_{V, r_V}}} \Bigg) \frac{d^2}{\gamma_s \gamma_w} ,\;\;\;
    \frac{d}{\gamma_w^2} ,\;\;\;
    \frac{d}{g_w^2}
    &= o\Big(\max_{k\in[K]} \{d_k^{-1}\} \Big) .
\end{align*}
\end{itemize}

With Assumption (L1), the standardised loading matrix $\Q_j:= \A_j \Z_j^{-1/2}$ satisfies $\Q_j' \Q_j \to\bSigma_{A,j}$ for $j\in[v-1]$, and $\Q_V:= \A_V \Z_V^{-1/2}$ satisfies $\Q_V' \Q_V \to\bSigma_{A,V}$. Similar implication holds for (L2), except that (L2) is only valid under the null. Hence with (L2), $\Z_V$ and $\bSigma_{A,V}$ in (L1) satisfy
\begin{equation}
\label{eqn: L1L2_equivalence}
\Z_V = \Z_K\otimes \dots\otimes \Z_v, \;\;\;
\bSigma_{A,V} = \bSigma_{A,K} \otimes \dots\otimes \bSigma_{A,v}.
\end{equation}
Note that the factor strength requirement for $\A_V$ in (L1) is satisfied by (L2), since from \eqref{eqn: L1L2_equivalence} $(\Z_V)_{r_V r_V} \asymp \prod_{i=v}^K d_i^{\delta_{i,r_i}} \geq d_V^{\min_{i=v}^K \delta_{i, r_i}} > d_V^{1/2}$. Assumption (L1) characterises the loading matrix behaviour generally for \eqref{eqn: test_kron_structure}, and the additional (L2) is specific for the null. Both assumptions allow for weak factors which are common feature in the literature \citep{LamYao2012, Onatski2012, CenLam2024}. When all factors are pervasive, for instance, \eqref{eqn: L1} can be interpreted as $d_j^{-1} \A_j' \A_j \to \bSigma_{A,j}$ if all factors are pervasive, which coincides with Assumption 3 of \cite{ChenFan2023} for matrix time series.

Assumption (F1) assumes that $\cF_{\textnormal{reshape},t}$ is a general linear process with weakly serial dependence. Theorem~\ref{thm: reshape} ensures that the core factor in \eqref{eqn: test_fm_w_Kron} (under $H_0$) reserves its structure of (F1) such that
\begin{equation}
\label{eqn: F1_under_H0}
\cF_{t} = \sum_{w\geq 0} a_{f,w} \cX_{f,t-w},
\;\;\; \text{with }
\Reshape(\cX_{f,t}, \cA) =\cX_{\textnormal{reshape},f,t} .
\end{equation}
Note that it holds for each $k\in[K]$, as $T\to\infty$,
\begin{equation}
\label{eqn: cFt_limit}
\frac{1}{T} \sum_{t=1}^T \mat{k}{\cF_{t}} \, \mat{k}{\cF_{t}}' \xrightarrow{p} \rmk\, \I_{r_k} ,
\end{equation}
which is direct from Proposition 1.3 in the supplement of \cite{CenLam2024}. In comparison, \cite{Heetal2022a} and \cite{Barigozzietal2023b} assume the form of \eqref{eqn: cFt_limit} with $\rmk\, \I_{r_k}$ replaced by a positive definite matrix. This does not imply (F1) is particularly stronger as our factor loading matrices already incorporate some positive definite matrices by (L1) and (L2).

Assumptions (E1) and (E2) depict a general noise time series on the factor models \eqref{eqn: test_fm_wo_Kron} and \eqref{eqn: test_fm_w_Kron}. It is worth noting that the noise tensor $\cE_t$ is allowed to be (weakly) dependent across modes and time, regardless of the existence of Kronecker product structure. From \eqref{eqn: E1}, we have
\begin{equation}
\label{eqn: cEt_reshape}
\begin{split}
    \Reshape(\c{E}_t, \cA) &=
    \Reshape(\c{F}_{e,t}, \cA) \times_1 \A_{e,1}\times_2 \cdots \times_{v-1} \A_{e,v-1} \times_v (\A_K \otimes \dots \otimes \A_v) \\
    &\hspace{12pt}
    + \Reshape(\bSigma_{\epsilon}, \cA) \circ \Reshape(\bepsilon_t, \cA),
\end{split}
\end{equation}
so that the structure of (E1) and (E2) are preserved by $\Reshape(\c{E}_t, \cA)$. Assumption (R1) details the rate assumptions on factor strengths and is hence satisfied automatically when all factors are pervasive. (R2) also concerns factor strength and would hold for all strong factors if $v> 1$; for $v=1$, (R2) holds when $\min_{k\in[K]} d_k=o(T)$ in addition to strong factors.

\begin{remark}\label{remark: R2}
When $v=1$ and all factors are strong, Assumption (R2) requires $d_{k*} = \min_{k\in[K]} d_k = o(T)$ which seems restricted. This is to ensure the asymptotic normality when we aggregate $d_{k*}$ number of estimated residuals in $x_{j,t}$ and $y_{m,j,t}$ in Section~\ref{subsec: construct_test}. However, from the proof of Theorem~\ref{thm: noise_aggregate_dist} in the supplement, it is feasible to aggregate $d_{k^*}^\beta$ for any $0< \beta < 2$ such that $d_{k*}^\beta = o(T)$. Therefore, (R2) is in fact arguably as mild as Assumption B5 in \cite{Heetal2023}. We do not pursue such an aggregation scheme to keep the practical procedure as simple as possible.
\end{remark}

\subsection{Main results and practical test design}\label{subsec: main_theorem}

We first present below the results for our residual estimators in \eqref{eqn: tilde_cEt} and \eqref{eqn: hat_cEt}, which inspire the testing procedure in Section~\ref{subsec: construct_test}. Following Theorem~\ref{thm: noise_aggregate_dist}, the theoretical guarantee of the test is also provided.

\begin{theorem}\label{thm: noise_aggregate_dist}
Let Assumptions (F1), (L1), (L2), (E1), (E2), (R1) and (R2) hold. With the notations in Section~\ref{subsec: construct_test}, under $H_0$, there exists $m\in [|\c{R}|]$ such that for each $t\in[T]$, $j\in [d/d_k^\ast]$,
\begin{align*}
    &\frac{\sum_{i=1}^{d_{k^\ast}} \big(\wh{E}_{m,t,(k^\ast),ij}^2 - \Sigma_{\epsilon, (k^\ast), ij}^2 \big)}{\sqrt{\sum_{i=1}^{d_{k^\ast}} \textnormal{Var}(\epsilon_{t, (k^\ast), ij}^2) \Sigma_{\epsilon, (k^\ast), ij}^4 }}, \; \frac{\sum_{i=1}^{d_{k^\ast}} \big(\wt{E}_{t,(k^\ast),ij}^2 - \Sigma_{\epsilon, (k^\ast), ij}^2 \big)}{\sqrt{\sum_{i=1}^{d_{k^\ast}} \textnormal{Var}(\epsilon_{t, (k^\ast), ij}^2) \Sigma_{\epsilon, (k^\ast), ij}^4 }}
    \xrightarrow{p} Z_{j},
\end{align*}
where $Z_{j} \xrightarrow{\c{D}} \c{N}(0,1)$ and $Z_{h}$ is independent of $Z_{\ell}$ for $h\neq \ell$. Under $H_1$, the asymptotic result for $\wt{E}_{t,(k^\ast),ij}$  above still holds true.
\end{theorem}

\begin{theorem}\label{thm: test}
Let all the assumptions in Theorem~\ref{thm: noise_aggregate_dist} hold. With the notations in Section~\ref{subsec: construct_test}, we have the following for any $j\in [d/d_{k^\ast}]$ under $H_0$. There exists $m\in [|\c{R}|]$ such that, as $T,d_1,\dots,d_K \to \infty$,
\[
\b{P}_{y,m,j}\big( y_{m,j,t} \geq \wh{q}_{x,j}(\alpha) \big) \leq \alpha .
\]
In other words, as $T,d_1,\dots,d_K \to \infty$, it holds with probability going to 1 that
\begin{equation}
\label{eqn: min_test_probability}
\min_{m\in [|\c{R}|]} \Big\{ \frac{1}{T} \sum_{t=1}^T \b{1}\{y_{m,j,t} \geq \wh{q}_{x,j}(\alpha)\} \Big\}
\leq \alpha .
\end{equation}
\end{theorem}

Theorem~\ref{thm: test} presents the grounds for our construction of the test statistic in (\ref{eqn: reject_rule}).
For related explanations, see the discussions immediately after (\ref{eqn: CDF}), and before Remark \ref{remark: divisor_combination_rV}.


The setup of the problem \eqref{eqn: test_kron_structure} specifies the set $\cA$ which is only needed in $H_1$ due to \eqref{eqn: test_fm_w_Kron} under $H_0$. It is direct to specify $\cA$ for a series of matrix-valued observations (i.e., order-2 tensor), see Example~\ref{example: test_order2}. However, for a general order-$K$ tensor with $K\geq 3$, $\cA$ might be misspecified without any prior knowledge. To resolve this, we present the second theorem on tensor reshape as follows.

\begin{theorem}\label{thm: reshape2}
\textnormal{(Tensor Reshape Theorem II)}
Consider a tensor time series $\{\cY_t\}$ and a set of mode indices $\cA$. With Definition~\ref{def: kron_structure}, the time series $\big\{ \Reshape(\cY_t, \cA) \big\}$ has a Kronecker product structure if and only if $\{\cY_t\}$ either has a Kronecker product structure or has no Kronecker product structure along a subset of $\cA$.
\end{theorem}

Suppose now $\{\cY_t\}$ has no Kronecker product structure along some $\cA^\ast$. Theorem~\ref{thm: reshape2} tells us that testing the Kronecker product structure of the reshaped series $\big\{ \Reshape(\cY_t, \cA) \big\}$ effectively tests if $\cA^\ast \subseteq \cA$. In light of this, a testing design is feasible when $\cA$ is unspecified, with a minimal assumption that $\Reshape(\cY_t,[K]) = \vec{\cY_t}$ has a factor structure, i.e., the vectorised $\cY_t$ follows a standard vector factor model. For illustration, consider $\Reshape(\cY_t, [K] \setminus \{1\}) = \Reshape(\cY_t, \{2,\dots, K\})$ which is an order-2 tensor. Using the last property of Reshape$(\cdot,\cdot)$ in Section~\ref{subsec: reshape}, we have
\begin{align*}
    \Reshape\Big( \Reshape(\cY_t, \{2,\dots, K\}), \{1,2\} \Big) = \Reshape(\cY_t, \{1, 2,\dots, K\}) = \vec{\cY_t}.
\end{align*}
According to Definition~\ref{def: kron_structure}, $\big\{\Reshape( \cY_t, \{2,\dots, K\}) \big\}$ follows a factor model along $\{1,2\}$. This is always correctly specified since $\big\{\vec{\cY_t} \big\}$ follows a factor model (which also implies $\cA^\ast \subseteq [K]$). By Theorem~\ref{thm: reshape2}, $\Reshape(\cY_t, \{2,\dots, K\})$ has no Kronecker product structure if and only if $1\in \cA^\ast$. Hence on testing \eqref{eqn: test_kron_structure} with $\cY_t$ replaced by $\{\Reshape(\cY_t, \{2,\dots, K\})\}$ and $\cA=\{1,2\}$, rejection of the null implies $1\in \cA^\ast$.

By the fact that $\big\{\vec{\cY_t} \big\}$ with any permutation on $\vec{\cY_t}$ also follows a factor model, the above scheme is in fact valid on $\Reshape(\cY_t, [K] \setminus \{k\})$ for any $k\in[K]$. Eventually, $\cA^\ast$ can be identified, and the above procedure is summarised into the following algorithm.

\bigskip
\noindent \ul{Practical testing algorithm}
\begin{itemize}
    \item[1.] Given an order-$K$ tensor time series $\{\cY_t\}$  with $K\geq 2$ and $\vec{\cY_t}$ following a factor model with $r_\text{vec}$ number of factors, initialise $\wh\cA^\ast =\phi$, the empty set.
    \item[2.] Initialise $k=1$. Define a test as \eqref{eqn: test_kron_structure} with $\{\cY_t\}$ replaced by $\big\{ \Reshape(\cY_t, [K] \setminus \{k\}) \big\}$ and $\cA$ by $\{1,2\}$.
    \item[3.] Follow the steps in Section~\ref{subsec: construct_test} to test the problem in step 2, with $r_V$ replaced by $r_\text{vec}$. If the null is rejected, include $k$ in the set $\wh\cA^\ast$.
    \item[4.] Repeat from step 2 to step 3 with $k=2,3,\dots,K$.
\end{itemize}
\bigskip

With the resulted $\wh\cA^\ast$ from the algorithm, we interpret that $\{\cY_t\}$ has no Kronecker product structure along $\wh\cA^\ast$. In practice, $\wh\cA^\ast$ being an empty set implies $\{\cY_t\}$ has a Kronecker product structure.

\begin{remark}\label{remark: wild_test}
Definition~\ref{def: kron_structure} considers the absence of Kronecker product structure over a single set $\cA$ only, which does not fully characterise all scenarios for $\cY_t$ with order at least 4. However, we do not pursue this complication here, albeit our practical design can be readily adapted.
\end{remark}

\section{Numerical Studies}\label{sec: empirical}

\subsection{Simulations}\label{subsec: simulation}

In this section, we demonstrate the empirical performance of our test by Monte Carlo simulations. As discussed in Section~\ref{subsec: model_KPS}, the test is only non-trivial when the data order $K$ is at least 2. We hence consider from $K=2$ to $K=4$.

The data generating processes adapt Assumptions (F1), (E1) and (E2). Specifically, we set the number of factors as $r_k=2$ for any $k\in[K]$, and first generate $\cF_t$ in \eqref{eqn: test_fm_w_Kron} with each element being independent standardised AR(2) with AR coefficient 0.7 and -0.3. The elements in $\cF_{e,t}$ and $\bm\epsilon_t$ are generated similarly, but their AR coefficients are (-0.5, 0.5) and (0.4, 0.4) respectively. The standard deviation of each element in $\bm\epsilon_t$ is generated by \text{i.i.d.} $|\cN(0,1)|$. Unless otherwise specified, all innovation processes in constructing $\cF_t$, $\cF_{e,t}$ and $\bepsilon_t$ are \text{i.i.d.} standard normal. For each $j\in[v-1]$, each factor loading matrix $\A_j$ is generated independently with $\A_j=\bf{U}_j\B_j$, where each entry of $\bf{U}_j\in \b{R}^{d_j\times r_j}$ is \text{i.i.d.} $\cN(0,1)$, and $\B_j\in \b{R}^{r_j\times r_j}$ is diagonal with the $h$-th diagonal entry being $d_j^{-\zeta_{j,h}}$, $0\leq \zeta_{j,h}\leq 0.5$. Pervasive factors have $\zeta_{j,h}=0$, while weak factors have $0<\zeta_{j,h}\leq 0.5$. Each entry of $\A_{e,j} \in\b{R}^{d_j\times r_{e,j}}$ is \text{i.i.d.} $\cN(0,1)$, but has independent probability of 0.95 being set exactly to 0. We set $r_{e,k}=2$ for all $j\in[v-1]$ throughout all experiments. For any $\cA$ (specified later), we obtain
\[
\Reshape(\cF_t, \cA) = \cF_{\textnormal{reshape},t}, \;\;\;
\Reshape(\cE_t, \cA) = \cE_{\textnormal{reshape},t}.
\]
Lastly, similar to $\{\A_j\}_{j\in[v-1]}$, we generate $\{\A_v,\dots,\A_K\}$ and let $\A_V = \A_K \otimes \A_{K-1} \otimes \dots \otimes \A_v$ under $H_0$, or $\A_V$ directly under $H_1$. Whenever $r_V$ is required, it is computed as $\prod_{j\in\cA} r_j $, which has $2^{|\cA|}$ combinations since we set two factors for each mode. According to \eqref{eqn: test_fm_wo_Kron} and \eqref{eqn: test_fm_w_Kron}, we then respectively construct $\Reshape(\cY_t, \cA)$ (and hence the corresponding $\cY_t$) and  $\cY_t$ directly.

We mainly consider a series of performance indicators and each simulation setting is repeated 500 times. With notations in Section~\ref{subsec: construct_test}, we calculate the following with $\alpha\in \{0.01, 0.05\}$:
\begin{equation}
\label{eqn: hat_alpha_p_def}
\begin{split}
    \wh\alpha &:= \min_{m\in [|\c{R}|]} \Big\{ \frac{1}{T} \sum_{t=1}^T \b{1}\{y_{m,1,t} \geq \wh{q}_{x,1}(\alpha)\} \Big\} , \\
    \wh{p} &:= \b{1}\{ \wh{q}_\alpha 
    \leq\alpha \} , \;\;\;
    \text{where }
    \wh{q}_\alpha := \min_{m\in [|\c{R}|]} \Bigg\{ \text{5\% quantile of } \frac{1}{T} \sum_{t=1}^T \b{1}\{y_{m,j,t} \geq \wh{q}_{x,j}(\alpha)\} \text{ over $ j\in [d/d_{k^\ast}]$} \Bigg\},
\end{split}
\end{equation}
where $\wh\alpha$ is the significance level under the measure $\b{P}_{y,m,1}$ taken minimum over $m\in [|\c{R}|]$, and $\wh{p}$ is an indicator function of the decision rule \eqref{eqn: reject_rule} leading to retaining $H_0$. Under $H_0$, we expect $\wh\alpha$ to be close to $\alpha$ and $\wh{p}$ to be 1 according to Theorem~\ref{thm: test}.

\subsubsection*{Test size and power}

Consider first $H_0$ with $\cA$ containing the last two modes of $\cY_t$, i.e., $\cA=\{1,2\}$ for $K=2$, $\cA=\{2,3\}$ for $K=3$ and $\cA=\{3,4\}$ for $K=4$. We experiment on all pervasive factors. Table~\ref{tab: simulation_test_size} presents the simulation results under various settings for $K=2,3,4$, and all of them well align with Theorem~\ref{thm: test}. 
Note that for $K =
3, 4$, all $\wh{p}$’s are 1, and for K = 2, the proportion of repetitions with $\wh{p}$ = 1 is increasing with dimensions and time in general.
The results under $H_1$ are presented in Table~\ref{tab: simulation_test_power} which confirms the power of our test. While larger dimensions generally improve the test performance, it is unsurprising from Table~\ref{tab: simulation_test_power} that under the same $(T,d_k)$ setting, testing the Kronecker product structure along two modes on $\cY_t$ is harder for higher-order $\cY_t$. This is reasonable since the testing problem \eqref{eqn: test_kron_structure} is genuinely harder when $\A_V$ plays a less significant role in a higher order data. To demonstrate this, suppose $K=3$, $(T,d_1,d_2,d_3) = (360,10,15,20)$, and all factors are pervasive. We experiment through $\cA= \{1,2\}, \{1,3\}, \{2,3\}, \{1,2,3\}$. The results reported in Table~\ref{tab: simulation_test_different_cA} indeed shows that when the tested loading matrix $\A_V$ has a larger size, the test has larger power in general. The setting with $\cA=\{2,3\}$ is an exception, suggesting a potential issue of unbalanced spatial dimensions.

\begin{table}[htp!]\centering
\ra{1.3}
\begin{tabular}{@{}crrrrrrrrrrrrrrr@{}}\toprule
& \multicolumn{4}{c}{$K=2$} &  & \multicolumn{4}{c}{$K=3$} &  & \multicolumn{4}{c}{$K=4$} \\
\toprule
$T=120$
& \multicolumn{2}{c}{$d_k=15$} & \multicolumn{2}{c}{$d_k=30$} &  & \multicolumn{2}{c}{$d_k=15$} & \multicolumn{2}{c}{$d_k=30$} &  & \multicolumn{2}{c}{$d_k=10$} & \multicolumn{2}{c}{$d_k=15$} \\
\cmidrule{2-5} \cmidrule{7-10} \cmidrule{12-15}
$\alpha$
&  $1\%$ & $5\%$ & $1\%$ & $5\%$ && $1\%$ & $5\%$ & $1\%$ & $5\%$ && $1\%$ & $5\%$ & $1\%$ & $5\%$ \\ 
\cmidrule{2-5} \cmidrule{7-10} \cmidrule{12-15}
$\wh\alpha$ & $.020$ & $.071$ & $.020$ & $.078$ && $.013$ & $.055$ & $.013$ & $.055$ && $.012$ & $.054$ & $.012$ & $.053$ \\
$\wh{p}$ &  $.974$ & $.836$ & $.996$ & $.860$ && $1$ & $1$ & $1$ & $1$ && $1$ & $1$ & $1$ & $1$ \\
\midrule
$T=360$
& \multicolumn{2}{c}{$d_k=15$} & \multicolumn{2}{c}{$d_k=30$} &  & \multicolumn{2}{c}{$d_k=15$} & \multicolumn{2}{c}{$d_k=30$} &  & \multicolumn{2}{c}{$d_k=10$} & \multicolumn{2}{c}{$d_k=15$} \\
\cmidrule{2-5} \cmidrule{7-10} \cmidrule{12-15}
$\alpha$
&  $1\%$ & $5\%$ & $1\%$ & $5\%$ && $1\%$ & $5\%$ & $1\%$ & $5\%$ && $1\%$ & $5\%$ & $1\%$ & $5\%$ \\ 
\cmidrule{2-5} \cmidrule{7-10} \cmidrule{12-15}
$\wh\alpha$ &  $.011$ & $.057$ & $.012$ & $.059$ && $.010$ & $.051$ & $.010$ & $.052$ && $.010$ & $.051$ & $.010$ & $.051$ \\
$\wh{p}$ &  $.988$ & $.862$ & $1$ & $.842$ && $1$ & $1$ & $1$ & $1$ && $1$ & $1$ & $1$ & $1$ \\
\midrule
$T=720$
& \multicolumn{2}{c}{$d_k=15$} & \multicolumn{2}{c}{$d_k=30$} &  & \multicolumn{2}{c}{$d_k=15$} & \multicolumn{2}{c}{$d_k=30$} &  & \multicolumn{2}{c}{$d_k=10$} & \multicolumn{2}{c}{$d_k=15$} \\
\cmidrule{2-5} \cmidrule{7-10} \cmidrule{12-15}
$\alpha$
&  $1\%$ & $5\%$ & $1\%$ & $5\%$ && $1\%$ & $5\%$ & $1\%$ & $5\%$ && $1\%$ & $5\%$ & $1\%$ & $5\%$ \\ 
\cmidrule{2-5} \cmidrule{7-10} \cmidrule{12-15}
$\wh\alpha$ &  $.011$ & $.053$ & $.012$ & $.054$ && $.010$ & $.051$ & $.010$ & $.051$ && $.010$ & $.050$ & $.010$ & $.051$ \\
$\wh{p}$ & $.994$ & $.916$ & $1$ & $.920$ && $1$ & $1$ & $1$ & $1$ && $1$ & $1$ & $1$ & $1$ \\
\bottomrule
\end{tabular}
\caption {Results of $\wh\alpha$ and $\wh{p}$ under $H_0$ in \eqref{eqn: test_kron_structure} for $K=2,3,4$. For each setting, $d_k$ is the same for all $k\in[K]$. Each cell is the average of $\wh\alpha$ or $\wh{p}$ computed under the corresponding setting over 500 runs.}
\label{tab: simulation_test_size}
\end{table}

\begin{table}[htp!]\centering
\ra{1.3}
\begin{tabular}{@{}crrrrrrrrrrrrrrr@{}}\toprule
& \multicolumn{4}{c}{$K=2$} &  & \multicolumn{4}{c}{$K=3$} &  & \multicolumn{4}{c}{$K=4$} \\
\toprule
$T=120$
& \multicolumn{2}{c}{$d_k=15$} & \multicolumn{2}{c}{$d_k=30$} &  & \multicolumn{2}{c}{$d_k=15$} & \multicolumn{2}{c}{$d_k=30$} &  & \multicolumn{2}{c}{$d_k=10$} & \multicolumn{2}{c}{$d_k=15$} \\
\cmidrule{2-5} \cmidrule{7-10} \cmidrule{12-15}
$\alpha$
&  $1\%$ & $5\%$ & $1\%$ & $5\%$ && $1\%$ & $5\%$ & $1\%$ & $5\%$ && $1\%$ & $5\%$ & $1\%$ & $5\%$ \\ 
\cmidrule{2-5} \cmidrule{7-10} \cmidrule{12-15}
$\wh\alpha$ & $.839$ & $.898$ & $.928$ & $.956$ && $.674$ & $.742$ & $.790$ & $.834$ && $.583$ & $.655$ & $.649$ & $.712$ \\
$\wh{p}$ &  $0$ & $0$ & $0$ & $0$ && $0$ & $0$ & $0$ & $0$ && $.012$ & $.002$ & $0$ & $0$ \\
\midrule
$T=360$
& \multicolumn{2}{c}{$d_k=15$} & \multicolumn{2}{c}{$d_k=30$} &  & \multicolumn{2}{c}{$d_k=15$} & \multicolumn{2}{c}{$d_k=30$} &  & \multicolumn{2}{c}{$d_k=10$} & \multicolumn{2}{c}{$d_k=15$} \\
\cmidrule{2-5} \cmidrule{7-10} \cmidrule{12-15}
$\alpha$
&  $1\%$ & $5\%$ & $1\%$ & $5\%$ && $1\%$ & $5\%$ & $1\%$ & $5\%$ && $1\%$ & $5\%$ & $1\%$ & $5\%$ \\ 
\cmidrule{2-5} \cmidrule{7-10} \cmidrule{12-15}
$\wh\alpha$ &  $.818$ & $.888$ & $.917$ & $.951$ && $.659$ & $.738$ & $.776$ & $.832$ && $.571$ & $.653$ & $.636$ & $.709$ \\
$\wh{p}$ &  $0$ & $0$ & $0$ & $0$ && $0$ & $0$ & $0$ & $0$ && $.002$ & $0$ & $0$ & $0$ \\
\midrule
$T=720$
& \multicolumn{2}{c}{$d_k=15$} & \multicolumn{2}{c}{$d_k=30$} &  & \multicolumn{2}{c}{$d_k=15$} & \multicolumn{2}{c}{$d_k=30$} &  & \multicolumn{2}{c}{$d_k=10$} & \multicolumn{2}{c}{$d_k=15$} \\
\cmidrule{2-5} \cmidrule{7-10} \cmidrule{12-15}
$\alpha$
&  $1\%$ & $5\%$ & $1\%$ & $5\%$ && $1\%$ & $5\%$ & $1\%$ & $5\%$ && $1\%$ & $5\%$ & $1\%$ & $5\%$ \\ 
\cmidrule{2-5} \cmidrule{7-10} \cmidrule{12-15}
$\wh\alpha$ &  $.817$ & $.885$ & $.918$ & $.951$ && $.652$ & $.731$ & $.787$ & $.837$ && $.559$ & $.640$ & $.629$ & $.701$ \\
$\wh{p}$ & $0$ & $0$ & $0$ & $0$ && $0$ & $0$ & $0$ & $0$ && $.002$ & $0$ & $0$ & $0$ \\
\bottomrule
\end{tabular}
\caption {Results of $\wh\alpha$ and $\wh{p}$ under $H_1$ in \eqref{eqn: test_kron_structure} for $K=2,3,4$. Refer to Table~\ref{tab: simulation_test_size} for the explanation of each cell.}
\label{tab: simulation_test_power}
\end{table}

\begin{table}[htp!]\centering
\ra{1.3}
\begin{tabular}{@{}crrrrrrrrrrrrrrr@{}}\toprule
$H_0$ & \multicolumn{2}{c}{$\cA=\{1,2\}$} &  & \multicolumn{2}{c}{$\cA=\{1,3\}$} &  & \multicolumn{2}{c}{$\cA=\{2,3\}$} & & \multicolumn{2}{c}{$\cA=\{1,2,3\}$} \\
\cmidrule{2-12}
$\alpha$
&  $1\%$ & $5\%$ && $1\%$ & $5\%$ && $1\%$ & $5\%$ && $1\%$ & $5\%$ \\ 
\cmidrule{2-3} \cmidrule{5-6} \cmidrule{8-9} \cmidrule{11-12}
$\wh\alpha$ & $.010$ & $.051$ && $.010$ & $.052$ && $.010$ & $.052$ && $.014$ & $.063$  \\
$\wh{p}$ & $1$ & $1$ && $1$ & $1$ && $1$ & $1$ && $1$ & $.956$  \\
\midrule
$H_1$ & \multicolumn{2}{c}{$\cA=\{1,2\}$} &  & \multicolumn{2}{c}{$\cA=\{1,3\}$} &  & \multicolumn{2}{c}{$\cA=\{2,3\}$} & & \multicolumn{2}{c}{$\cA=\{1,2,3\}$} \\
\cmidrule{2-12}
$\alpha$
&  $1\%$ & $5\%$ && $1\%$ & $5\%$ && $1\%$ & $5\%$ && $1\%$ & $5\%$ \\ 
\cmidrule{2-3} \cmidrule{5-6} \cmidrule{8-9} \cmidrule{11-12}
$\wh\alpha$ & $.702$ & $.775$ && $.705$ & $.779$ && $.673$ & $.748$ && $.927$ & $.959$  \\
$\wh{p}$ & $0$ & $0$ && $0$ & $0$ && $0$ & $0$ && $0$ & $0$  \\
\bottomrule
\end{tabular}
\caption {Results of $\wh\alpha$ and $\wh{p}$ over different $\cA$'s in \eqref{eqn: test_kron_structure} for $(T,d_1,d_2,d_3) = (360,15,20,25)$. Refer to Table~\ref{tab: simulation_test_size} for the explanation of each cell. For each $\cA= \{1,2\}, \{1,3\}, \{2,3\}, \{1,2,3\}$, the number of rows of $\A_V$ in \eqref{eqn: test_kron_structure} is respectively $300, 375, 500, 7500$.}
\label{tab: simulation_test_different_cA}
\end{table}

\subsubsection*{Robustness for weak factor, heavy-tailed noise and misspecified number of factors}

In the following, we fix $K=3$ and $\cA=\{2,3\}$ to investigate the robustness of our test. Consider Setting I and II, each with four sub-settings:
\begin{itemize}
\item[(Ia)] $T=180, \; d_1=d_2=d_3=15$. All factors are pervasive with $\zeta_{j,h} =0$.
\item[(Ib)] Same as (Ia), but one factor is weak with $\zeta_{j,1}=0.1$.
\item[(Ic)] Same as (Ia), but both factors are weak with $\zeta_{j,1}= \zeta_{j,2}=0.1$.
\item[(Id)] Same as (Ia), but all innovation processes in constructing $\cF_t$, $\cF_{e,t}$ and $\bepsilon_t$ are \text{i.i.d.} $t_3$.
\item[(IIa-d)] Same as (Ia) to (Id) respectively, except that $r_V$ is randomly specified from $\{2,3,4,5,6\}$ with equal probability.
\end{itemize}

Setting (Ia) is our benchmark and all other settings feature some defects from weak factors, heavy-tailed noise, or misspecified number of factors. Table~\ref{tab: simulation_robust} reports the results for both $H_0$ and $H_1$. In contrast to (Ia), all other settings have lower test power to various extents. However, the size of the test is hardly influenced by weak factors or heavy-tailed noise from the results of (Ib), (Ic) and (Id). Although number-of-factor misspecification is detrimental, our decision rule $\wh{p}$ still has satisfying performance.

\begin{table}[htp!]\centering
\ra{0.95}
\begin{tabular}{@{}crrrrrrrrrrrrrrrrrrrrrrrr@{}}\toprule
& \multicolumn{11}{c}{Setting I} \\
\toprule
$H_0$
& \multicolumn{2}{c}{(Ia)} && \multicolumn{2}{c}{(Ib)} && \multicolumn{2}{c}{(Ic)} && \multicolumn{2}{c}{(Id)} \\
\cmidrule{2-3} \cmidrule{5-6} \cmidrule{8-9} \cmidrule{11-12}
$\alpha$
&  $1\%$ & $5\%$ && $1\%$ & $5\%$ && $1\%$ & $5\%$ && $1\%$ & $5\%$ \\ 
\cmidrule{2-3} \cmidrule{5-6} \cmidrule{8-9} \cmidrule{11-12}
$\wh\alpha$ &  $.008$ & $.054$ && $.008$ & $.053$ && $.008$ & $.053$ && $.008$ & $.053$ \\
$\wh{p}$ &  $1$ & $1$ && $1$ & $1$ && $1$ & $1$ && $1$ & $1$  \\
\midrule
$H_1$
& \multicolumn{2}{c}{(Ia)} && \multicolumn{2}{c}{(Ib)} && \multicolumn{2}{c}{(Ic)} && \multicolumn{2}{c}{(Id)} \\
\cmidrule{2-3} \cmidrule{5-6} \cmidrule{8-9} \cmidrule{11-12}
$\alpha$
&  $1\%$ & $5\%$ && $1\%$ & $5\%$ && $1\%$ & $5\%$ && $1\%$ & $5\%$ \\ 
\cmidrule{2-3} \cmidrule{5-6} \cmidrule{8-9} \cmidrule{11-12}
$\wh\alpha$ &  $.691$ & $.765$ && $.593$ & $.684$ && $.441$ & $.553$ && $.519$ & $.693$  \\
$\wh{p}$ &  $0$ & $0$ && $.014$ & $0$ && $.034$ & $.004$ && $.070$ & $.002$  \\
\midrule
& \multicolumn{11}{c}{Setting II} \\
\midrule
$H_0$
& \multicolumn{2}{c}{(IIa)} && \multicolumn{2}{c}{(IIb)} && \multicolumn{2}{c}{(IIc)} && \multicolumn{2}{c}{(IId)} \\
\cmidrule{2-3} \cmidrule{5-6} \cmidrule{8-9} \cmidrule{11-12}
$\alpha$
&  $1\%$ & $5\%$ && $1\%$ & $5\%$ && $1\%$ & $5\%$ && $1\%$ & $5\%$ \\ 
\cmidrule{2-3} \cmidrule{5-6} \cmidrule{8-9} \cmidrule{11-12}
$\wh\alpha$ &  $.054$ & $.113$ && $.035$ & $.091$ && $.021$ & $.074$ && $.030$ & $.092$ \\
$\wh{p}$ &  $.972$ & $.932$ && $.988$ & $.966$ && $.998$ & $.996$ && $.996$ & $.964$  \\
\midrule
$H_1$
& \multicolumn{2}{c}{(IIa)} && \multicolumn{2}{c}{(IIb)} && \multicolumn{2}{c}{(IIc)} && \multicolumn{2}{c}{(IId)} \\
\cmidrule{2-3} \cmidrule{5-6} \cmidrule{8-9} \cmidrule{11-12}
$\alpha$
&  $1\%$ & $5\%$ && $1\%$ & $5\%$ && $1\%$ & $5\%$ && $1\%$ & $5\%$ \\ 
\cmidrule{2-3} \cmidrule{5-6} \cmidrule{8-9} \cmidrule{11-12}
$\wh\alpha$ &  $.553$ & $.652$ && $.504$ & $.620$ && $.378$ & $.509$ && $.424$ & $.596$  \\
$\wh{p}$ &  $.034$ & $0$ && $.018$ & $.002$ && $.036$ & $.006$ && $.128$ & $.004$  \\
\bottomrule
\end{tabular}
\caption {Results of $\wh\alpha$ and $\wh{p}$ under $H_0$ and $H_1$ in \eqref{eqn: test_kron_structure} over Setting I and II. Refer to Table~\ref{tab: simulation_test_size} for the explanation of each cell.}
\label{tab: simulation_robust}
\end{table}

\subsubsection*{Numerical performance of the practical testing algorithm}

On the practical testing algorithm which does not require $\cA$ to be specified, we consider Setting III and IV with $K=3$, and each has three sub-settings:
\begin{itemize}
\item[(IIIa)] $T=360, \; d_1=d_2=d_3=10$. All factors are strong and the data has a Kronecker product structure.
\item[(IIIb)] Same as (Ia), but the data has no Kronecker product structure along $\{2,3\}$.
\item[(IIIc)] Same as (Ia), but the data has no Kronecker product structure along $\{1,2,3\}$.
\item[(IVa-c)] Same as (IIIa) to (IIIc) respectively, except that $T=720$.
\end{itemize}

Table~\ref{tab: simulation_practical} verifies that our algorithm is able to test the Kronecker product structure of a given data without pre-specifying $\cA$. The performance is improved with more observations, and the level of $\alpha=0.01$ works particularly well.

\begin{table}[htp!]\centering
\ra{0.95}
\begin{tabular}{@{}crrrrrrrrrrrrrrrrrrrrrrrrr@{}}\toprule
& \multicolumn{8}{c}{Setting III}
&& \multicolumn{8}{c}{Setting IV} \\
\cmidrule{2-9} \cmidrule{11-18}
& \multicolumn{2}{c}{(IIIa)} && \multicolumn{2}{c}{(IIIb)} && \multicolumn{2}{c}{(IIIc)} && \multicolumn{2}{c}{(IVa)} && \multicolumn{2}{c}{(IVb)} && \multicolumn{2}{c}{(IVc)}  \\
\cmidrule{2-3} \cmidrule{5-6} \cmidrule{8-9} \cmidrule{11-12} \cmidrule{14-15} \cmidrule{17-18}
$\alpha$
&  $1\%$ & $5\%$ && $1\%$ & $5\%$ && $1\%$ & $5\%$
&& $1\%$ & $5\%$ && $1\%$ & $5\%$ && $1\%$ & $5\%$ \\ 
\cmidrule{2-3} \cmidrule{5-6} \cmidrule{8-9} \cmidrule{11-12} \cmidrule{14-15} \cmidrule{17-18}
Mode 1 &  $0$ & $.024$ && $.030$ & $.202$  && $1$ & $1$ && $0$ & $0$ && $.004$ & $.228$  && $1$ & $1$  \\
Mode 2 &  $0$ & $.036$ && $1$ & $1$  && $1$ & $1$ && $0$ & $.004$ && $1$ & $1$  && $1$ & $1$  \\
Mode 3 &  $0$ & $.034$ && $.998$ & $1$  && $1$ & $1$ && $0$ & $.002$ && $1$ & $1$  && $1$ & $1$  \\
\bottomrule
\end{tabular}
\caption {Results of Setting III and IV for the practical testing algorithm. Each cell is the fraction of the corresponding mode identified over 500 runs for the corresponding sub-settings.}
\label{tab: simulation_practical}
\end{table}

\subsection{Real data analysis}\label{subsec: realdataanalysis}

We apply our test on two real data examples described as follows.
\begin{itemize}
    \item [1.] New York city taxi traffic. The data considered includes all individual taxi rides operated by Yellow Taxi within Manhattan Island of New York City, published at 
        
        \url{https://www1.nyc.gov/site/tlc/about/tlc-trip-record-data.page}.
        
    The dataset contains trip records within the period of January 1, 2018 to December 31, 2022. We focus on the pick-up and drop-off dates/times, and pick-up and drop-off locations which are coded according to 69 predefined zones in the dataset. Moreover, each day is divided into 24 hourly periods to represent the pick-up and drop-off times, with the first hourly period from 0 \text{a.m.} to 1 a.m. Hence each day we have $\cY_t\in\b{R}^{69\times 69\times 24}$, where $y_{i_1,i_2,i_3,t}$ is the number of trips from zone $i_1$ to zone $i_2$ and the pick-up time is within the $i_3$-th hourly period on day $t$. We consider business days and non-business days separately, so that we will analyse two tensor time series. The business-day series and the non-business-day series are 1,260 and 566 days long, respectively.
    \item [2.] Fama-French portfolio returns. This is a set of portfolio returns data, where stocks are respectively categorised into ten levels of market equity and book-to-equity ratio which is the book equity for the last fiscal year divided by the end-of-year market equity; both criteria use NYSE deciles as breakpoints at the end of June each year. See details in 
        
        \url{https://mba.tuck.dartmouth.edu/pages/faculty/ken.french/Data_Library/det_100_port_sz.html}.
        
    The stocks in each of the $10\times 10$ categories form exactly two portfolios, one being value-weighted, and the other of equal-weight. That is, we will study two sets of 10 by 10 portfolios with their time series. We use monthly data from January 2010 to June 2021, and hence for both value-weighted and equal-weighted portfolios we have each of our data set as an order-2 tensor $\c{X}_t\in \b{R}^{10\times 10}$ for $t\in[138]$.
\end{itemize}

The two taxi series are order-3 tensor time series, and we only test their Kronecker product structure along $\cA= \{1,2\}$, i.e., we speculate that there is a merged ``location'' factor instead of ``pick-up'' and ``drop-off'' factors along mode-1 and -2 respectively. On the other hand, the two portfolio series are order-2 tensor time series, hence naturally we test along $\cA=\{1,2\}$. Furthermore, we remove the market effect via the capital asset pricing model (CAPM) as
\[
\vec{\c{X}_t} = \vec{ \bar{\c{X}}} + (r_t -\bar{r}) \bm{\beta} + \vec{\c{Y}_t},
\]
where $\vec{\c{X}_t} \in\b{R}^{100}$ is the vectorised returns at time $t$, $\vec{ \bar{\c{X}}}$ is the sample mean of $\vec{\c{X}_t}$, $\bm{\beta}$ is the coefficient vector, $r_t$ is the return of the NYSE composite index at time $t$, $\bar{r}$ is the sample mean of $r_t$, and $\vec{\c{Y}_t}$ is the CAPM residual. The least squares solution is
\[
\wh{\bm{\beta}} = \frac{\sum_{t=1}^{138} (r_t -\bar{r}) \{\vec{\c{X}_t} -\vec{ \bar{\c{X}}} \}}{\sum_{t=1}^{138} (r_t -\bar{r})^2} ,
\]
so that the estimated residual series $\{ \wh{\c{Y}}_t\}_{t \in[138]}$ with $\wh{\c{Y}}_t \in\b{R}^{10\times 10}$ is constructed as $\{\vec{\c{X}_t} -\vec{ \bar{\c{X}}} - (r_t -\bar{r}) \wh{\bm{\beta}} \}_{t \in[138]}$.

Hence, we study six time series in total: business-day taxi series, non-business-day taxi series, value-weighted portfolio series, equal-weighted portfolio series, value-weighted residual series and equal-weighted residual series. For each series, we perform the test described in Section~\ref{subsec: construct_test}. To estimate the rank, we use BCorTh by \cite{Chen_Lam}, iTIP-ER by \cite{Hanetal2022} and RTFA-ER by \cite{Heetal2022b} directly on each time series due to their large dimensions. Each mode of the six series has one or two estimated number of factors. Since the test results are similar for those rank settings, we present the results with two factors each mode and hence $\wh{r}_V=4$.

\begin{table}[htp!]\centering
\ra{1.3}
\begin{tabular}{@{}lrrrrrrccc@{}}\toprule
& \multicolumn{2}{c}{$\wh\alpha$} && \multicolumn{2}{c}{$\wh{q}_\alpha$} && \multicolumn{3}{c}{Tests in \cite{Heetal2023}} \\
\cmidrule{2-3} \cmidrule{5-6} \cmidrule{8-10}
& $1\%$ & $5\%$ && $1\%$ & $5\%$ && $H_0$ versus $H_{1,\text{row}}$ && $H_0$ versus $H_{1,\text{col}}$ \\ 
\cmidrule{2-3} \cmidrule{5-6}  \cmidrule{8-8} \cmidrule{10-10}
Business-day taxi &  $.020$ & $.093$  && $.002$ & $.003$ && - && - \\
Non-business-day taxi &  $.018$ & $.095$  && $.004$ & $.011$ && - && - \\
Value-weighted portfolio &  $.058$ & $.087$  && \textbf{.011} & \textbf{.053} && Not reject && Not reject \\
Equal-weighted portfolio &  $.036$ & $.051$  && \textbf{.018} & $.039$ && Not reject && Not reject \\
Value-weighted residual &  $.022$ & $.065$  && \textbf{.011} & $.047$ && Not reject && Not reject \\
Equal-weighted residual &  $.014$ & $.051$  && \textbf{.011} & $.047$ && Not reject && Not reject \\
\bottomrule
\end{tabular}
\caption {Test results for the studied series. The first two columns report the results for our hypothesis of interest \eqref{eqn: test_kron_structure} with $\cA=\{1,2\}$; $\wh{q}_\alpha$ larger than the corresponding $\alpha$ level is in bold. The last two columns report the results according to \cite{Heetal2023}.}
\label{tab: real_data_result}
\end{table}

In addition, we also conduct the hypotheses tests in \cite{Heetal2023} on our matrix time series data sets. To explain their hypotheses, for a matrix time series $\{\Y_t\}$ with $\Y_t\in \b{R}^{d_1\times d_2}$, under the null we have \eqref{eqn: test_fm_w_Kron} for $K=2$:
\[
H_0: \Y_t = \A_1 \F_t \A_2' +\E_t,
\]
where $\F_t\in \b{R}^{r_1\times r_2}$. However, under their two alternatives we test
\begin{align*}
    H_{1,\text{row}} : r_2 = 0, \;\;\; H_{1,\text{col}} : r_1 = 0,
\end{align*}
where according to \cite{Heetal2023}, $r_1 > 0, r_2 = 0$ (\text{resp.} $r_2 > 0, r_1 = 0$) denotes a one-way factor model along the row dimension, so that $\Y_t = \A_1 \F_{1,t} +\E_t$ with $\F_{1,t}\in \b{R}^{r_1\times d_2}$ (\text{resp.} the column dimension, so that $\Y_t = \F_{2,t} \A_2' +\E_t$ with $\F_{2,t}\in \b{R}^{d_1\times r_2}$), and $r_1 = r_2 = 0$ denotes the absence of any factor structure, so  that $\Y_t=\E_t$. All hyperparameter setups in Table 8 and 9 in \cite{Heetal2023} are experimented and all conclusions are the same.

Table~\ref{tab: real_data_result} reports $\wh\alpha$ and $\wh{q}_\alpha$ defined in \eqref{eqn: hat_alpha_p_def}, with $\alpha=0.01, \, 0.05$, together with the corresponding tests by \cite{Heetal2023}. For our hypothesis of interest, there is no evidence to reject the null for the two taxi series, but there is mild evidence (especially at $1\%$ level, with $\wh\alpha$ observed to be mildly larger than 1\%) to conclude that for the Fama-French time series, there is no Kronecker product structure along $\{1,2\}$. In other words, there is evidence to suggest that the portfolio return series
have structures deviating from 
the low-rank structure along its respective categorisations by market equity and book-to-equity ratio, meaning the vectorised data may have a more distinct factor structure. The comparisons between the portfolio and residual series justifies the removal of the market effect, which is intuitive since the market effect should be pervasive in financial return data and is irrelevant of our categorisations. On the other hand, we cannot reject the null by considering those alternative hypotheses considered in \cite{Heetal2023}.

\section{Appendix}\label{sec: Appendix}
Model identification, and proofs of all the theorems in this paper can be found in the supplement of this paper at \href{http://stats.lse.ac.uk/lam/Supp-KOFM.pdf}{http://stats.lse.ac.uk/lam/Supp-KOFM.pdf}.
Instruction in using our R package \texttt{KOFM} can be found \href{http://stats.lse.ac.uk/lam/intro-to-KOFM.html}{here}.

\newpage

\appendix

\begin{center}
{\large \textbf{Supplementary materials for the paper ``On Testing Kronecker Product Structure in Tensor Factor Models'' - Model identification, and proof of all theorems and auxiliary results}}
\end{center}

\section{Appendix: Identification}\label{sec: identification}

This section concerns the identification of the model in Definition~\ref{def: kron_structure}, following a discussion on Definition~\ref{def: kron_structure_set}. First, consider Definition~\ref{def: kron_structure}.1 so that in \eqref{eqn: fm_wo_Kron}, we have $\A_{\textnormal{reshape},K-\ell+1} \in \cK_{d_{a_1}\times \dots \times d_{a_\ell}}$. Given a general $\cA=\{a_1,\dots, a_\ell\}$, the ordered set of matrices $\{\A_j\}_{j \in[\kappa]}$ decomposing (as in Definition~\ref{def: kron_structure_set}) $\A_{\textnormal{reshape},K-\ell+1}$ might not be unique. For instance, suppose $K=2$, $d_1= d_2$, $\cA=\{1, 2\}$ (hence $\A_{\textnormal{reshape},K-\ell+1} \equiv \A_{\textnormal{reshape},1}$ has $d_1^2$ rows) and let $\A_{\textnormal{reshape},1}= d_1^{-1/4} \1_{d_1^2}$, then we have
\begin{equation}
\label{eqn: identification_example}
\A_{\textnormal{reshape},1} = \big( \underbrace{d_1^{-1/4} \1_{d_1}}_{\A_1} \big) \otimes \underbrace{\1_{d_1}}_{\A_2} = \underbrace{\1_{d_1}}_{\ddot{\A}_1} \otimes\, \big( \underbrace{d_1^{-1/4} \1_{d_1}}_{\ddot{\A}_2} \big),
\end{equation}
where it is clear that $\|\A_{1,\cdot 1}\|^2 \asymp d_1^{1/2}$ and $\|\ddot{\A}_{1,\cdot 1} \|^2 \asymp d_1$. Such defeat can be rectified by allocating the ``factor strength'' in $\A_{\textnormal{reshape},K-\ell+1}$ to each mode in $\cA$; see the following Assumption (S1) as one example.

\begin{itemize}
\item[(S1)]
\em
For \eqref{eqn: fm_w_Kron} such that $\A_{\textnormal{reshape},K-\ell+1} = \otimes_{i\in\cA} \A_i$ for a given $\cA$, we assume that for any $i\in\cA$,
\[
\frac{\|\A_i\|_F^2}{r_i d_i} = \Bigg( \frac{\|\A_{\textnormal{reshape},K-\ell+1} \|_F^2}{\prod_{i\in\cA} r_i d_i} \Bigg)^{1/|\cA|} .
\]
\end{itemize}
The issue of indeterminacy in the factor strength is fixed by (S1) which has same spirit of Assumption (IC) in \cite{ChenLam2024a}. Its heuristic is to allocate the factor strength in $\A_{\textnormal{reshape},K-\ell+1}$ to each mode according to their number of factors and dimensions. Note that Assumption (S1) holds automatically if all factors are pervasive. Now recall the example~\eqref{eqn: identification_example}, $\A_1= \A_2= d_1^{-1/8} \1_{d_1}$ are identified by (S1). Note that the discussion above is only for completeness and would not influence our testing problem.

Next, we identify models in Definition~\ref{def: kron_structure} and Theorem~\ref{thm: reshape}. For \eqref{eqn: fm_wo_Kron} and \eqref{eqn: fm_w_Kron}, we state below the two assumptions (L1') and (L2'), which are (notation-wise) general versions of (L1) and (L2), respectively.

\begin{itemize}
\item[(L1')] (Factor strength in $\A_{\textnormal{reshape},j}$)
\em
For each $j\in[K-\ell]$, we assume that $\A_{\textnormal{reshape},j}$ in \eqref{eqn: fm_wo_Kron} is of full rank and as $I_j\to\infty$,
\begin{equation}
\label{eqn: L1'}
\Z_{\textnormal{reshape},j}^{-1/2} \A_{\textnormal{reshape},j}' \A_{\textnormal{reshape},j} \Z_{\textnormal{reshape},j}^{-1/2}
\to \bSigma_{\textnormal{reshape},A,j},
\end{equation}
where $\bSigma_{\textnormal{reshape},A,j}$ is positive definite with all eigenvalues bounded away from 0 and infinity, and $\Z_{\textnormal{reshape},j}$ is a diagonal matrix with $(\Z_{\textnormal{reshape},j})_{hh}\asymp I_j^{\delta_{\textnormal{reshape},j,h}}$ for $h\in[\pi_j]$ and the ordered factor strengths $1/2<\delta_{\textnormal{reshape},j,\pi_j}\leq \dots\leq \delta_{\textnormal{reshape},j,2}\leq \delta_{\textnormal{reshape},j,1}\leq 1$.

We assume that $\A_{\textnormal{reshape}, j}$ for $j=K-\ell+1$ also has the above form, except that only the maximum and minimum factor strengths are ordered, i.e.,  $1/2<\delta_{\textnormal{reshape},j,\pi_j}\leq \delta_{\textnormal{reshape},j,h} \leq \delta_{\textnormal{reshape},j,1}\leq 1$ for any $h\in[\pi_j]$.
\end{itemize}

\begin{itemize}
\item[(L2')] (Factor strength in $\A_i$)
\em
For \eqref{eqn: fm_w_Kron} with a given $\cA$, we assume that for each $i\in\cA$, $\A_i$ is of full rank and as $d_i\to\infty$,
\begin{equation}
\label{eqn: L2'}
\Z_i^{-1/2} \A_{i}' \A_i \Z_i^{-1/2}
\to \bSigma_{A,i},
\end{equation}
where $\bSigma_{A,i}$ is positive definite with all eigenvalues bounded away from 0 and infinity, and $\Z_i$ is a diagonal matrix with $(\Z_i)_{hh}\asymp d_i^{\delta_{i,h}}$ for $h\in[r_i]$ and the ordered factor strengths $1/2<\delta_{i,r_i}\leq \dots\leq \delta_{i,2}\leq \delta_{i,1}\leq 1$.
\end{itemize}

Theorem~\ref{thm: identification} presents the identification of the model \eqref{eqn: fm_wo_Kron} both in general and with Kronecker product structure (equivalently \eqref{eqn: fm_w_Kron} according to Theorem~\ref{thm: reshape}). Its proof is given directly after the statement.

\begin{theorem}\label{thm: identification}
\textnormal{(Identification)}
Let Assumption (F1) and (L1') hold, and $\cA$ is given. Then the factor structure in \eqref{eqn: fm_wo_Kron} is asymptotically identified up to some invertible matrix $\M_j\in \b{R}^{\pi_j\times \pi_j}$ such that the following sets of factor structure are equivalent,
\[
\Big(\cF_{\textnormal{reshape},t}, \big\{\A_{\textnormal{reshape},j} \big\}_{j\in [K-\ell+1]} \Big) = \Big(\cF_{\textnormal{reshape},t} \times_{j=1}^{K-\ell+1} \M_j^{-1}, \big\{\A_{\textnormal{reshape},j} \M_j \big\}_{j\in [K-\ell+1]} \Big).
\]

Let Assumption (L2') further holds. With a Kronecker product structure on \eqref{eqn: fm_wo_Kron}, we have \eqref{eqn: fm_w_Kron} where for each $k\in[K]$, $\A_k$ has unique rank and the factor structure in \eqref{eqn: fm_w_Kron} is asymptotically identified up to some invertible matrices.
\end{theorem}

\textbf{\textit{Proof of Theorem~\ref{thm: identification}.}}
Consider first \eqref{eqn: fm_wo_Kron}. Let $\big(\ddot\cF_{\textnormal{reshape},t}, \big\{\ddot\A_{\textnormal{reshape},j} \big\}_{j\in [K-\ell+1]} \big)$ be another set of parameters such that $\ddot\cF_{\textnormal{reshape},t} \times_{j=1}^{K-\ell+1} \ddot\A_{\textnormal{reshape},j} = \cF_{\textnormal{reshape},t} \times_{j=1}^{K-\ell+1} \A_{\textnormal{reshape},j}$. Define
\begin{equation}
\label{eqn: def_A_reshape_mj}
\A_{\textnormal{reshape},\text{-}j} := \A_{\textnormal{reshape}, K-\ell+1} \otimes \dots \otimes \A_{\textnormal{reshape},j+1} \otimes \A_{\textnormal{reshape},j-1} \otimes \dots \otimes \A_{\textnormal{reshape},1}.
\end{equation}
Define $\ddot\A_{\textnormal{reshape},\text{-}j}$ similarly. Without loss of generality, for any $j\in[K-\ell+1]$ we write
\begin{equation}
\label{eqn: tilde_A_reshape_rewrite}
\ddot\A_{\textnormal{reshape},j} = \A_{\textnormal{reshape},j} \M_{\textnormal{reshape},j} + \bGamma_{\textnormal{reshape},j}, \;\;\;
\text{where } \bGamma_{\textnormal{reshape},j}' \A_{\textnormal{reshape},j} = \0,
\end{equation}
with $\M_{\textnormal{reshape},j} \in\b{R}^{\pi_j \times \pi_j}$ and $\bGamma_{\textnormal{reshape},j} \in\b{R}^{I_j\times \pi_j}$, but can have zero columns. Then
\begin{align*}
    \0 &= \bGamma_{\textnormal{reshape},j}' \A_{\textnormal{reshape},j} \F_{\textnormal{reshape},t,(j)} \A_{\textnormal{reshape},\text{-}j}' \ddot\A_{\textnormal{reshape},\text{-}j} \\
    &= \bGamma_{\textnormal{reshape},j}' \ddot\A_{\textnormal{reshape},j} \ddot\F_{\textnormal{reshape},t,(j)} \ddot\A_{\textnormal{reshape},\text{-}j}' \ddot\A_{\textnormal{reshape},\text{-}j} \\
    &= \bGamma_{\textnormal{reshape},j}' \bGamma_{\textnormal{reshape},j} \ddot\F_{\textnormal{reshape},t,(j)} \ddot\A_{\textnormal{reshape},\text{-}j}' \ddot\A_{\textnormal{reshape},\text{-}j} ,
\end{align*}
which can only be true in general if $\bGamma_{\textnormal{reshape},j} =\0$ since $\ddot\F_{\textnormal{reshape},t,(j)} $ is random by Assumption (F1) and $\ddot\A_{\textnormal{reshape},\text{-}j}' \ddot\A_{\textnormal{reshape},\text{-}j}$ converges to some full rank matrix by (L1'). Hence $\M_{\textnormal{reshape},j}$ has full rank, and $\ddot\A_{\textnormal{reshape},j}$ and $\A_{\textnormal{reshape},j}$ share the same column space. $\cF_{\textnormal{reshape},t}$ is identified once $\big\{\A_{\textnormal{reshape},j} \big\}_{j\in [K-\ell+1]}$ is given correspondingly.

Suppose \eqref{eqn: fm_wo_Kron} has a Kronecker product structure and consider \eqref{eqn: fm_w_Kron}. By an argument similar to \eqref{eqn: tilde_A_reshape_rewrite} but over $k\in[K]$ (omitted), each matrix $\A_k$ for $k\in[K]$ has a unique rank and is identified up to some invertible matrix using (L2'). Hence $\cF_t$ is also identified, which completes the proof of the theorem. $\square$

\section{Appendix: Proof of theorems and lemmas}\label{appendix: proof}

{\large \textbf{A high-level summary of proofs:}}\\
The design of \eqref{eqn: hat_cCt} and \eqref{eqn: hat_cEt} over all divisor combinations of $r_V$ is due to pseudo-ranks from mode-$v$ to mode-$K$ of $\cY_t$ under $H_1$ in \eqref{eqn: test_kron_structure}. On the other hand, under $H_0$, there must be one $m\in[|\c{R}|]$ such that $(\pi_{m,1}, \dots, \pi_{m,K-v+1}) = (r_v, \dots, r_K)$. Hence we consider throughout the proof that $\{r_k\}_{k\in[K]}$ in \eqref{eqn: test_fm_w_Kron} (hence under $H_0$) is correctly specified. Hence, we simplify the notations $\wh\Q_{m,i}$, $\wh\cC_{m,t}$ and $\wh\cE_{m,t}$ as $\wh\Q_{i}$, $\wh\cC_{t}$ and $\wh\cE_{t}$, respectively. Whenever (L2) is assumed, we are implicitly considering $H_0$ in \eqref{eqn: test_kron_structure}, where we also define
\begin{align*}
    \Qmk &:= \Q_K\otimes \dots \otimes \Q_{k+1} \otimes \Q_{k-1}\otimes \dots \otimes \Q_1 \;\;\; \text{for $k\in[K]$}, \\
    \cF_{Z,t} &:= \cF_t \times_{k=1}^K \Z_k^{1/2} \;\;\; \text{with $\cF_t$ from \eqref{eqn: F1_under_H0}}.
\end{align*}

Theorem~\ref{thm: reshape} reveals the importance of the reshape operator and is the key to formalise the testing problem. Lemma~\ref{lemma: correlation_Et_Ft} to Lemma~\ref{lemma: norm_hat_D} serve as technical steps. The steps of all other proofs are as follows.
\begin{itemize}
    \item [1.] Under $H_0$, consider \eqref{eqn: test_fm_w_Kron}. We first derive the rate of convergence for $\wh\Q_k$ as an estimator of $\Q_k$ for each $k\in[K]$ (Lemma~\ref{lemma: Q_hat_rate}). Then the rates for the corresponding core factor and hence the common component can also be obtained (Lemma~\ref{lemma: hat_rate_cFcC}).
    \item [2.] Under $H_0$, consider \eqref{eqn: test_fm_wo_Kron}. We derive the rates of convergence for $\{\wt\Q_j\}_{j\in[v-1]}$ and $\wt\Q_V$ as estimators of $\{\Q_j\}_{j\in[v-1]}$ and $\Q_V$, respectively. Similar to step 1, the rate for the common component is obtained. See Lemma~\ref{lemma: Q_tilde_rate}.
    \item [3.] With steps 1 and 2, we then show that $\sum_{i=1}^{d_{k^\ast}} \wh{E}_{m,t,(k^\ast),ij}^2$ and $\sum_{i=1}^{d_{k^\ast}} \wt{E}_{t,(k^\ast),ij}^2$ under $H_0$ have the same distribution asymptotically (Theorem~\ref{thm: noise_aggregate_dist}).
    \item [4.] With step 3, the test statistic can be constructed with theoretical support (Theorem~\ref{thm: test}).
\end{itemize}

\subsection{Lemmas and proofs}

By the interplay between the core factor and the noise series in Assumptions (F1), (E1) and (E2), we state below Lemma~\ref{lemma: correlation_Et_Ft} which is direct from Proposition 1.1 and 1.2 of \cite{CenLam2024}.

\begin{lemma}\label{lemma: correlation_Et_Ft}
Let Assumptions (F1), (E1) and (E2) hold. Then
\begin{itemize}
    \item [1.]
    (Weak correlation of noise $\cE_t$ across different modes and time) there exists some positive constant $C<\infty$ so that for any $t\in[T]$, $k\in[K]$, $i_k, j\in[d_k]$, $h\in [\dmk]$, we have $\b{E}\cE_{t,i_1, \dots,i_K}=0$, $\b{E}\cE_{t,i_1, \dots,i_K}^4 \leq C$, and
    \begin{align*}
    & \sum_{j=1}^{d_k} \sum_{l=1}^{\dmk}
    \Big|\b{E}[E_{t,(k), i_k h}
    E_{t,(k), jl}]\Big| \leq C, \\
    & \sum_{l=1}^{\dmk} \sum_{s=1}^T \Big| \textnormal{Cov}\big( E_{t,(k),i_k h} E_{t,(k),jh}, E_{s,(k),i_k l} E_{s,(k),jl} \big) \Big| \leq C.
    \end{align*}
    \item [2.]
    (Weak dependence between factor $\cF_t$ and noise $\cE_t$) there exists some positive constant $C < \infty$ so that for any $k\in[K]$, $j\in[d_k]$, and any deterministic vectors $\bf{u}\in\b{R}^{r_k}$ and $\bf{v}\in\b{R}^{\rmk}$ with constant magnitudes, it holds for $\cF_t$ in \eqref{eqn: test_fm_w_Kron} that
    \[
    \b{E}\Bigg(
        \frac{1}{(\dmk T)^{1/2}}\sum_{h=1}^{\dmk}
        \sum_{t=1}^T E_{t,(k),jh}
        \bf{u}' \F_{t,(k)} \bf{v}
    \Bigg)^2 \leq C.
    \]
    \item [3.]
    Statement 2 holds similarly for $\Reshape(\cF_t,\cA)$ and $\Reshape(\cE_t,\cA)$.
\end{itemize}
\end{lemma}

\begin{lemma}\label{lemma: ft_lim}
Under Assumption (F1), with $\gamma_{v}:= \prod_{j=1}^{v-1} r_j$, we have as $T\to\infty$,
\begin{align}
    \frac{1}{T} \sum_{t=1}^T \F_{t,(k)} \M \F_{t,(k)}' \xrightarrow{p} \tr(\M) \cdot \I_{r_k} ,&\label{eqn: Ftk_lim} \\
    \frac{1}{T} \sum_{t=1}^T \F_{\textnormal{reshape},t,(j)} \bf{N} \F_{\textnormal{reshape},t,(j)}' \xrightarrow{p} \tr(\bf{N}) \cdot \I_{r_j} ,&\label{eqn: reshape_Ftk_lim} \\
    \frac{1}{T} \sum_{t=1}^T \F_{\textnormal{reshape},t,(v)} \bf{W} \F_{\textnormal{reshape},t,(v)}' \xrightarrow{p} \tr(\bf{W}) \cdot \I_{r_V} ,&\label{eqn: reshape_Ftv_lim}
\end{align}
where $\F_{t,(k)} \in\b{R}^{r_k\times \rmk}$ for $k\in[K]$ represents the mode-$k$ unfolding of $\cF_t$ in \eqref{eqn: test_fm_w_Kron}, and $\M$ is any $\rmk \times \rmk$ matrix independent of $\{\cF_t \}_{t\in[T]}$ with $\|\M\|_F$ bounded in probability; similarly for $\F_{\textnormal{reshape},t,(j)} \in\b{R}^{r_k\times (r_V \gamma_{v}/r_j)}$ for $j\in[v-1]$ in \eqref{eqn: test_fm_wo_Kron}; similarly for $\F_{\textnormal{reshape},t,(v)} \in\b{R}^{r_V\times \gamma_{v}}$ in \eqref{eqn: test_fm_wo_Kron}.
\end{lemma}

\textbf{\textit{Proof of Lemma~\ref{lemma: ft_lim}.}}
It suffices to show \eqref{eqn: Ftk_lim} and the other two follow similarly. With (F1), consider first
\begin{equation}
\label{eqn: F1_coef_bound}
\begin{split}
    &\hspace{12pt}
    \sum_{s=1}^T \Big(\sum_{w\geq 0} a_{f,w} a_{f,w-(t-s)} \Big) \Big(\sum_{q\geq 0} a_{f,q} a_{f,q-(t-s)} \Big) \\
    &\leq
    \Big(\sum_{s=1}^T \sum_{w\geq 0} |a_{f,w}| |a_{f,w-(t-s)}| \Big) \Big(\sum_{q\geq 0} |a_{f,q}| \Big) \cdot \max_q |a_{f,q}| =
    O(1) \cdot \Big( \sum_{w\geq 0} |a_{f,w}| \Big)^2 =O(1),
\end{split}
\end{equation}
where the last two equality used Assumption (F1). Similarly, it holds that
\begin{equation}
\label{eqn: F1_coef_bound2}
\sum_{s=1}^T \sum_{w\geq 0} a_{f,w}^2 a_{f,w-(t-s)}^2 \leq \Big(\sum_{w\geq 0} a_{f,w}^2 \Big)^2 =O(1).
\end{equation}

Now for $i\neq j \in[r_k]$, by Assumption (F1) we have
\begin{align*}
    &\hspace{12pt}
    \b{E}\Big(\frac{1}{T} \sum_{t=1}^T \F_{t,(k),i\cdot}' \M \F_{t,(k),j\cdot} \Big) = \b{E}\Big(\frac{1}{T} \sum_{t=1}^T \F_{t,(k),i\cdot}' \Big) \b{E}\Big(\M \F_{t,(k),j\cdot} \Big) = 0.
\end{align*}
For $i=j \in[r_k]$, with $\cX_{f,t}$ from \eqref{eqn: F1_under_H0}, we have
\begin{align*}
    &\hspace{12pt}
    \b{E}\Big(\frac{1}{T} \sum_{t=1}^T \F_{t,(k),i\cdot}' \M \F_{t,(k),i\cdot} \Big) = \frac{1}{T} \sum_{t=1}^T \sum_{w\geq 0} \sum_{q\geq 0} a_{f,w} a_{f,q} \b{E}\Big(\X_{f,t-w,(k),i\cdot}' \M \X_{f,t-q,(k),i\cdot} \Big) \\
    &=
    \frac{1}{T} \sum_{t=1}^T \sum_{w\geq 0} \sum_{q\geq 0} a_{f,w} a_{f,q} \tr(\M) \b{1}\{w=q\} = \tr(\M) \sum_{w\geq 0} a_{f,w}^2 = \tr(\M),
\end{align*}
so that $\b{E}\Big( T^{-1} \sum_{t=1}^T \F_{t,(k)} \M \F_{t,(k)}' \Big) =\tr(\M) \cdot \I_{r_k}$. Finally, consider any $i,j\in [r_k]$,
\begin{align*}
    &\hspace{12pt}
    \text{Var}\Big( \frac{1}{T} \sum_{t=1}^T \F_{t,(k),i\cdot}' \M \F_{t,(k),j\cdot} \Big) =
    \text{Var}\Big( \frac{1}{T} \sum_{t=1}^T \sum_{l=1}^{\rmk} \sum_{v=1}^{\rmk} F_{t,(k),il} M_{lv} F_{t,(k),jv} \Big) \\
    &=
    \frac{1}{T^2} \text{Cov}\Big( \sum_{t=1}^T \sum_{l=1}^{\rmk} \sum_{v=1}^{\rmk} \sum_{w\geq 0} \sum_{q\geq 0} a_{f,w} a_{f,q} X_{f,t-w,(k),il} M_{lv} X_{f,t-q,(k),jv}, \\
    &\hspace{48pt}
    \sum_{s=1}^T \sum_{g=1}^{\rmk} \sum_{u=1}^{\rmk} \sum_{h\geq 0} \sum_{m\geq 0} a_{f,h} a_{f,m} X_{f,s-h,(k),ig} M_{gu} X_{f,s-m,(k),ju} \Big) \\
    &=
    \frac{1}{T^2} \sum_{t=1}^T \sum_{s=1}^T \Big(\sum_{w\geq 0} a_{f,w} a_{f,w-(t-s)} \Big) \Big(\sum_{q\geq 0} a_{f,q} a_{f,q-(t-s)} \Big) \Big(\sum_{l=1}^{\rmk} \sum_{v=1}^{\rmk} M_{lv}^2\Big) \\
    &\hspace{12pt}
    + \frac{1}{T^2} \sum_{t=1}^T \sum_{w\geq 0} \sum_{q\geq 0} \sum_{s=1}^T \sum_{h\geq 0} \sum_{m\geq 0} a_{f,w} a_{f,q} a_{f,h} a_{f,m} \Big(\sum_{l=1}^{\rmk} M_{ll}^2\Big) \\
    &\hspace{36pt}
    \cdot \text{Cov}\Big( X_{f,t-w,(k),il} X_{f,t-q,(k),il}, X_{f,s-h,(k),il} X_{f,s-m,(k),il} \Big) \\
    &=
    O\Big( \frac{1}{T}\Big) \cdot \Big(\sum_{l=1}^{\rmk} \sum_{v=1}^{\rmk} M_{lv}^2\Big) + O\Big( \frac{1}{T}\Big) \cdot \Big(\sum_{l=1}^{\rmk} M_{ll}^2\Big) = O\Big( \frac{1}{T}\Big) \cdot \|\M \|_F^2 =o(1),
\end{align*}
where the third equality considered $i=j$ and $i\neq j$ separately, and the fourth used \eqref{eqn: F1_coef_bound} and \eqref{eqn: F1_coef_bound2}. This completes the proof of \eqref{eqn: Ftk_lim} and hence Lemma~\ref{lemma: ft_lim}. $\square$

\begin{lemma}\label{lemma: rate_R_kt}
(Bounding $\sum_{t=1}^T\R_{k,t}$) Under Assumptions (F1), (L1), (L2), (E1) and (E2), it holds that
\begin{align}
    & \Big\| \sum_{t=1}^T \Q_k \F_{Z,t,(k)} \Qmk' \E_{t,(k)}' \Big\|_F^2 =  O_P\Big(T d_k^{1+\delta_{k,1}} \dmk \Big),
    \label{eqn: CE_bound} \\
    & \Big\| \sum_{t=1}^T \E_{t,(k)} \E_{t,(k)}' \Big\|_F^2 = O_P\Big(T d_k^2 \dmk + T^2 d_k \dmk^2 \Big).
    \label{eqn: EE_bound}
\end{align}
Hence, with $\R_{k,t}$ defined in \eqref{eqn: sample_cov_unfolding}, we have
\[
\Big\| \sum_{t=1}^T \R_{k,t} \Big\|_F^2 = O_P\Big(T d_k^2 \dmk + T^2 d_k \dmk^2 \Big) .
\]
\end{lemma}

\textbf{\textit{Proof of Lemma~\ref{lemma: rate_R_kt}.}}
It is not hard to see \eqref{eqn: CE_bound} holds as follows,
\begin{equation*}
\begin{split}
    &\hspace{12pt}
    \Big\| \sum_{t=1}^T \Q_k \F_{Z,t,(k)} \Qmk' \E_{t,(k)}' \Big\|_F^2 =
    \sum_{i=1}^{d_k} \sum_{l=1}^{d_k}
    \Bigg( \sum_{t=1}^T \A_{k, i\cdot}' \F_{t,(k)} \Amk' \E_{t,(k), l\cdot} \Bigg)^2 \\
    &=
    \sum_{i=1}^{d_k} \|\A_{k, i\cdot}\|^2 \cdot \sum_{l=1}^{d_k} \Bigg( \sum_{h=1}^{\dmk} \sum_{t=1}^T E_{t,(k), lh} \frac{1}{\|\A_{k, i\cdot}\|} \A_{k, i\cdot}' \F_{t,(k)} \A_{\text{-}k, h\cdot} \Bigg)^2
    = O_P\Big(T d_k^{1+\delta_{k,1}} \dmk\Big) ,
\end{split}
\end{equation*}
where the last equality is from Assumptions (L1), (L2) and Lemma~\ref{lemma: correlation_Et_Ft}.

Consider now \eqref{eqn: EE_bound}. First note from Assumption (E1) that, for any $k\in[K]$, $i\in[d_k]$, $j\in[\dmk]$,
\[
E_{t,(k),ij} = \A_{e,k, i\cdot}' \F_{e,t,(k)} \A_{e,\text{-}k, j\cdot} + \Sigma_{\epsilon,(k), ij} \epsilon_{t,(k), ij},
\]
where $\A_{e,\text{-}k}:= \A_{e,K}\otimes \dots \otimes \A_{e,k+1} \otimes \A_{e,k-1} \otimes \dots\otimes \A_{e,1}$. Then with Assumption (E2), we have
\[
\text{Cov}(E_{t,(k),ij}, E_{t,(k),lj}) = \A_{e,k, i\cdot}'
\A_{e,k, l\cdot} \|\A_{e,\text{-}k, j\cdot}\|^2 + \Sigma_{\epsilon,(k), ij}^2 \b{1}_{\{i=l\}} ,
\]
and together with Lemma~\ref{lemma: correlation_Et_Ft},
\begin{equation*}
\begin{split}
    &\hspace{12pt}
    \b{E}\Big\{ \Big\| \sum_{t=1}^T \E_{t,(k)} \E_{t,(k)}' \Big\|_F^2 \Big\} = \sum_{i=1}^{d_k} \sum_{l=1}^{d_k} \b{E}\Big\{ \Big(\sum_{t=1}^T \sum_{j=1}^{\dmk} E_{t,(k),ij} E_{t,(k),lj} \Big)^2 \Big\} \\
    &=
    \sum_{i=1}^{d_k} \sum_{l=1}^{d_k} \Big\{\sum_{t=1}^T \sum_{j=1}^{\dmk} \sum_{s=1}^T \sum_{h=1}^{\dmk} \text{Cov}(E_{t,(k),ij} E_{t,(k),lj}, E_{s,(k),ih} E_{s,(k),lh}) + \Big( \sum_{t=1}^T \sum_{j=1}^{\dmk} \b{E}[E_{t,(k),ij} E_{t,(k),lj}] \Big)^2 \Big\} \\
    &=
    O\big( Td_k^2 \dmk\big) + \sum_{i=1}^{d_k} \sum_{l=1}^{d_k} O\Big( T\cdot \A_{e,k, i\cdot}' \A_{e,k, l\cdot} \|\A_{e,\text{-}k} \|_F^2 + T\dmk \cdot\b{1}_{\{i=l\}} \Big)^2
    = O\big(T d_k^2 \dmk + T^2 d_k \dmk^2 \big).
\end{split}
\end{equation*}
Recall from \eqref{eqn: sample_cov_unfolding} it is defined that
\[
\R_{k,t} = \Q_k \F_{Z,t,(k)} \Qmk' \E_{t,(k)}' + \E_{t,(k)} \Qmk \F_{Z,t,(k)}' \Q_k' + \E_{t,(k)} \E_{t,(k)}' ,
\]
the rate on $\big\| \sum_{t=1}^T \R_{k,t} \big\|_F^2$ is direct from \eqref{eqn: CE_bound} and \eqref{eqn: EE_bound}. Hence the proof of Lemma~\ref{lemma: rate_R_kt} is complete. $\square$

\begin{lemma}\label{lemma: norm_hat_D}
Let Assumptions (F1), (L1), (L2), (E1), (E2) and (R1) hold. For $k\in[K]$, define $\wh\D_k$ as the $r_k \times r_k$ diagonal matrix with the first largest $r_k$ eigenvalues of $T^{-1} \sum_{t=1}^T \Y_{t,(k)} \Y_{t,(k)}'$ on the main diagonal, such that $\wh\D_k = \wh\Q_k' \big(T^{-1} \sum_{t=1}^T \Y_{t,(k)} \Y_{t,(k)}' \big) \wh\Q_k$. Define $\omega_k:= d_k^{\delta_{k,r_k} -\delta_{k,1}} g_s$. It holds that
\begin{equation*}
     \big\| \wh\D_k^{-1} \big\|_F = O_P\big( \omega_k^{-1} \big) .
\end{equation*}
\end{lemma}

\textbf{\textit{Proof of Lemma~\ref{lemma: norm_hat_D}.}}
 Observe that $\wh\D_k$ has size $r_k\times r_k$, it suffices to find the lower bound of $\lambda_{r_k}(\wh\D_k)$. To do this, consider the decomposition
\begin{align}
\label{eqn: CC_decomposition}
  \frac{1}{T}\sum_{t=1}^T \Y_{t,(k)} \Y_{t,(k)}' = \frac{1}{T}\sum_{t=1}^T \Q_k \F_{Z,t,(k)} \Qmk' \Qmk \F_{Z,t,(k)}' \Q_k' + \frac{1}{T}\sum_{t=1}^T \R_{k,t},
\end{align}
which is direct from \eqref{eqn: sample_cov_unfolding}. Then for a unit vector $\bgamma \in\b{R}^{d_k}$, we can define
\begin{align}
    S_k(\bgamma) &:=
    \frac{1}{\omega_k}\bm{\gamma}' \Big(\frac{1}{T}\sum_{t=1}^T \Y_{t,(k)} \Y_{t,(k)}' \Big) \bm{\gamma}
    =: S_k^*(\bgamma) + \wt{S}_k(\bgamma), \; \text{ with} \notag\\
    S_k^*(\bgamma) &:=
    \frac{1}{\omega_k}\bm{\gamma}' \Big(\frac{1}{T}\sum_{t=1}^T \Q_k \F_{Z,t,(k)} \Qmk' \Qmk \F_{Z,t,(k)}' \Q_k' \Big) \bm{\gamma}
    , \;\;\;
    \wt{S}_k(\bgamma) :=
    \frac{1}{\omega_k}\bm{\gamma}' \Big(\frac{1}{T}\sum_{t=1}^T \R_{k,t} \Big) \bm{\gamma}
    . \notag
\end{align}
Since $\|\bgamma\|=1$, we have by Lemma~\ref{lemma: rate_R_kt},
\begin{align*}
    &\hspace{12pt}
    |\wt{S}_k(\bgamma)|^2 \leq \frac{1}{\omega_k^2 T^2} \Big\|\sum_{t=1}^T \R_{k,t} \Big\|_F^2
    = O_P\Big( T^{-1} d_k^{1+ 2(\delta_{k,1} -\delta_{k, r_k})} d g_s^{-2} + d_k^{2(\delta_{k,1} -\delta_{k, r_k} -1/2)} d^2 g_s^{-2} \Big) = o_P(1),
\end{align*}
where the last equality used Assumption (R1). Next, with Assumption (F1) and Lemma~\ref{lemma: ft_lim}, consider
\begin{equation}
\label{eqn: CC_lambda_rk}
\begin{split}
    &\hspace{12pt}
    \lambda_{r_k} \Big(\frac{1}{T}\sum_{t=1}^T \Q_k \F_{Z,t,(k)} \Qmk' \Qmk \F_{Z,t,(k)}' \Q_k' \Big)
    = \lambda_{r_k}\Big( \frac{1}{T}\sum_{t=1}^T \A_k \F_{t,(k)} \Amk' \Amk \F_{t,(k)}' \A_k' \Big) \\
    &\geq
    \lambda_{r_k}(\A_k' \A_k)\cdot \lambda_{r_k}\Big( \frac{1}{T} \sum_{t=1}^T \F_{t,(k)} \Amk' \Amk \F_{t,(k)}' \Big)
    \asymp_P
    d_k^{\delta_{k,r_k}} \cdot \lambda_{r_k}(\tr(\Amk' \Amk) \I_{r_k}) \asymp_P \omega_k.
\end{split}
\end{equation}
With this, going back to the decomposition \eqref{eqn: CC_decomposition},
\begin{align*}
    &\hspace{12pt}
    \omega_k^{-1} \lambda_{r_k}(\wh\D_k) = \omega_k^{-1} \lambda_{r_k} \Big(\frac{1}{T} \sum_{t=1}^T \Y_{t,(k)} \Y_{t,(k)}' \Big) \\
    &\geq
    \omega_k^{-1}\lambda_{r_k} \Big(\frac{1}{T} \sum_{t=1}^T \Q_k \F_{Z,t,(k)} \Qmk' \Qmk \F_{Z,t,(k)}' \Q_k' \Big) - \sup_{\|\bgamma\|=1} |\wt{S}_k(\bgamma)|\asymp_P 1,
\end{align*}
which implies $\norm{\wh\D_k^{-1}}_F = O_P\big( \lambda_{r_k}^{-1}(\wh\D_k) \big) = O_P\big( \omega_k^{-1} \big)$, which completes the proof of Lemma~\ref{lemma: norm_hat_D}. $\square$

\begin{lemma}\label{lemma: Q_hat_rate}
(Consistency of $\{\wh\Q_k\}_{k\in[K]}$)
Let Assumptions (F1), (L1), (L2), (E1), (E2) and (R1) hold. For any $k\in[K]$, define an $r_k\times r_k$ matrix
\[
\H_k := T^{-1}\wh\D_k^{-1} \wh\Q_k' \Q_k \sum_{t=1}^T \big( \F_{Z,t,(k)} \Qmk' \Qmk \F_{Z,t,(k)}' \big).
\]
As $T,d_1,\dots, d_K \to\infty$ we have $\H_k$ invertible with $\|\H_k\|_F =O_P(1)$ and
\begin{align*}
    \big\| \wh\Q_k - \Q_k\H_k' \big\|_F^2 =
    O_P\Bigg\{ d_k^{2(\delta_{k,1} - \delta_{k,r_k})} \Bigg( \frac{1}{T\dmk d_k^{1- \delta_{k,1}}} +\frac{1}{ d_k^{1+ \delta_{k, r_k}}} \Bigg) \frac{d^2}{g_s^2} \Bigg\} .
\end{align*}
\end{lemma}

\textbf{\textit{Proof of Lemma~\ref{lemma: Q_hat_rate}.}}
First, we may write \eqref{eqn: test_fm_w_Kron} as
\begin{equation}
\label{eqn: fm_w_Kron_unfolding}
\cY_t = \cF_{Z,t} \times_1 \Q_1 \times_2 \dots \times_K \Q_K + \cE_t .
\end{equation}
For any $k\in[K]$, taking the mode-$k$ unfolding on \eqref{eqn: fm_w_Kron_unfolding}, we have
\[
\Y_{t,(k)} = \Q_k \F_{Z,t,(k)} \Qmk' + \E_{t,(k)},
\]
where $\F_{Z,t,(k)}$ denotes the mode-$k$ unfolding of $\cF_{Z,t}$. Hence,
\begin{equation}
\label{eqn: sample_cov_unfolding}
\begin{split}
    \Y_{t,(k)} \Y_{t,(k)}' &=
    \Q_k \F_{Z,t,(k)} \Qmk' \Qmk \F_{Z,t,(k)}' \Q_k' + \Q_k \F_{Z,t,(k)} \Qmk' \E_{t,(k)}' + \E_{t,(k)} \Qmk \F_{Z,t,(k)}' \Q_k' + \E_{t,(k)} \E_{t,(k)}' \\
    &=: \Q_k \F_{Z,t,(k)} \Qmk' \Qmk \F_{Z,t,(k)}' \Q_k' + \R_{k,t}.
\end{split}
\end{equation}
Recall from Lemma~\ref{lemma: norm_hat_D} that $\wh\D_k$ is the $r_k \times r_k$ diagonal matrix with the first largest $r_k$ eigenvalues of $T^{-1} \sum_{t=1}^T \Y_{t,(k)} \Y_{t,(k)}'$ on the main diagonal, and since $\wh\Q_k$ consists of the corresponding eigenvectors, we have
\begin{equation}
\label{eqn: sample_cov_unfolding_svd}
    \wh\Q_k \wh\D_k = \frac{1}{T} \sum_{t=1}^T \Y_{t,(k)} \Y_{t,(k)}' \wh\Q_k .
\end{equation}
With (\ref{eqn: sample_cov_unfolding}), we can write the $j$-th row of estimated mode-$k$ factor loading as
\begin{equation*}
\begin{split}
    \wh\Q_{k,j\cdot} &= \frac{1}{T} \wh\D_k^{-1} \sum_{i=1}^{d_k} \wh\Q_{k,i\cdot} \sum_{t=1}^T \big( \Y_{t,(k)} \Y_{t,(k)}' \big)_{ij} \\
    &=
    \frac{1}{T} \wh\D_k^{-1} \sum_{i=1}^{d_k} \wh\Q_{k,i\cdot} \Q_{k,i\cdot}' \sum_{t=1}^T \big( \F_{Z,t,(k)} \Qmk' \Qmk \F_{Z,t,(k)}' \big) \Q_{k,j\cdot}
    + \frac{1}{T} \wh\D_k^{-1} \sum_{i=1}^{d_k} \wh\Q_{k,i\cdot} \sum_{t=1}^T ( \R_{k,t} )_{ij}.
\end{split}
\end{equation*}
Hence with the definition $\H_k = T^{-1}\wh\D_k^{-1} \wh\Q_k' \Q_k \sum_{t=1}^T \big( \F_{Z,t,(k)} \Qmk' \Qmk \F_{Z,t,(k)}' \big)$, we may decompose
\begin{align}
    &\hspace{12pt}
    \wh\Q_{k,j\cdot} - \H_k\Q_{k,j\cdot} =
    \frac{1}{T} \wh\D_k^{-1} \sum_{i=1}^{d_k} \wh\Q_{k,i\cdot} \sum_{t=1}^T ( \R_{k,t} )_{ij}
    \label{eqn: Q_norm_decomp1} \\
    &=
    \frac{1}{T} \wh\D_k^{-1} \sum_{i=1}^{d_k} \big( \wh\Q_{k, i\cdot} - \H_k \Q_{k,i\cdot} \big) \sum_{t=1}^T (\R_{k,t})_{ij} + \frac{1}{T} \wh\D_k^{-1} \sum_{i=1}^{d_k} \H_k \Q_{k,i\cdot} \sum_{t=1}^T (\R_{k,t})_{ij} .
    \label{eqn: Q_norm_decomp2}
\end{align}
With the decomposition \eqref{eqn: Q_norm_decomp1}, it holds that
\begin{equation}
\label{eqn: Q_norm_rate1}
\begin{split}
    &\hspace{12pt}
    \big\| \wh\Q_k - \Q_k\H_k' \big\|_F^2 = \sum_{j=1}^{d_k} \big\| \wh\Q_{k,j\cdot} - \H_k \Q_{k,j\cdot} \big\|^2 = \sum_{j=1}^{d_k} \Big\| \frac{1}{T} \wh\D_k^{-1} \sum_{i=1}^{d_k} \wh\Q_{k,i\cdot} \sum_{t=1}^T ( \R_{k,t} )_{ij} \Big\|^2 \\
    &=
    \sum_{j=1}^{d_k} \Big\| \frac{1}{T} \wh\D_k^{-1} \wh\Q_{k}' \Big(\sum_{t=1}^T \R_{k,t} \Big)_{j\cdot} \Big\|^2
    \leq \frac{1}{T^2} \cdot \|\wh\D_k^{-1}\|_F^2 \cdot \|\wh\Q_k\|_F^2 \cdot \Big\| \sum_{t=1}^T \R_{k,t} \Big\|_F^2 \\
    &=
    O_P\Bigg( d_k^{2 (\delta_{k,1} - \delta_{k,r_k})} \bigg(\frac{1}{T\dmk} +\frac{1}{d_k} \bigg) \frac{d^2}{g_s^2} \Bigg) = o_P(1),
\end{split}
\end{equation}
where the second last equality used Lemmas \ref{lemma: rate_R_kt} and \ref{lemma: norm_hat_D}, and the last used Assumption (R1).

Before improving the rate of \eqref{eqn: Q_norm_rate1}, we now use it to show $\H_k$ has full rank and $\|\H_k \|_F =O_P(1)$ asymptotically. To this end, it is sufficient to observe
\begin{align*}
    \I_{r_k} &= \wh\Q_k' \wh\Q_k = \wh\Q_k' \big( \wh\Q_k -\Q_k \H_k' \big) + \wh\Q_k' \Q_k \H_k' = \Q_k' \wh\Q_k \H_k' + o_P(1) \\
    &=
    \H_k \Q_k' \Q_k \H_k' + o_P(1) = \H_k \bSigma_{A,k} \H_k' + o_P(1),
\end{align*}
where the last equality used Assumption (L1) or (L2) and it is immediate that $\H_k$ has full rank asymptotically. Let $\sigma_i(\X)$ denote the $i$-th largest singular value for any give matrix $\X$, we have
\[
\sigma_1(\H_k) \cdot \sigma_{r_k}(\bSigma_{A,k}) \cdot \sigma_{r_k}(\H_k^T) \leq \sigma_1(\H_k) \cdot \sigma_{r_k}(\bSigma_{A,k} \H_k^T) \leq \sigma_1(\H_k \bSigma_{A,k} \H_k') = O_P(\sigma_1(\I_{r_k})) = O_P(1),
\]
which implies $\|\H_k \|_F =O_P(1)$ by Assumption (L1) or (L2).

Consider the decomposition \eqref{eqn: Q_norm_decomp2}, we have
\begin{equation}
\label{eqn: Q_norm_rate2}
\begin{split}
    &\hspace{12pt}
    \big\| \wh\Q_k - \Q_k\H_k' \big\|_F^2 = \sum_{j=1}^{d_k} \big\| \wh\Q_{k,j\cdot} - \H_k \Q_{k,j\cdot} \big\|^2 \\
    &=
    O_P\Bigg( \sum_{j=1}^{d_k} \Big\|\frac{1}{T} \wh\D_k^{-1} \sum_{i=1}^{d_k} \H_k \Q_{k,i\cdot} \sum_{t=1}^T (\R_{k,t})_{ij} \Big\|^2 + \sum_{j=1}^{d_k} \Big\|\frac{1}{T} \wh\D_k^{-1} \sum_{i=1}^{d_k} \big( \wh\Q_{k, i\cdot} - \H_k \Q_{k,i\cdot} \big) \sum_{t=1}^T (\R_{k,t})_{ij} \Big\|^2 \Bigg) \\
    &=
    O_P\Bigg( \sum_{j=1}^{d_k} \Big\|\frac{1}{T} \wh\D_k^{-1} \H_k \sum_{i=1}^{d_k} \Q_{k,i\cdot} \sum_{t=1}^T (\R_{k,t})_{ij} \Big\|^2\Bigg) + O_P\Bigg( \frac{1}{T^2} \big\|\wh\D_k^{-1} \big\|_F^2 \cdot \big\| \wh\Q_{k} -\Q_{k} \H_k' \big\|_F^2 \cdot \Big\| \sum_{t=1}^T \R_{k,t} \Big\|_F^2 \Bigg) \\
    &=:
    O_P\Big( T^{-2} d_k^{2(\delta_{k,1} -\delta_{k,r_k})} g_s^{-2} \Big) \cdot \sum_{j=1}^{d_k} \big(\c{I}_{1,j} +\c{I}_{2,j} +\c{I}_{3,j} \big) + o_P\Big( \big\|\wh\Q_{k} -\Q_{k} \H_k' \big\|_F^2 \Big),
\end{split}
\end{equation}
where the last equality used $\|\H_k \|_F =O_P(1)$, Lemmas~\ref{lemma: rate_R_kt} and \ref{lemma: norm_hat_D}, and the definitions
\begin{align*}
    & \c{I}_{1,j} := \Big\| \sum_{i=1}^{d_k} \Q_{k,i\cdot} \sum_{t=1}^T \Q_{k,i\cdot}' \F_{Z,t,(k)} \Qmk' \E_{t,(k), j\cdot} \Big\|^2 , \;\;\;
    \c{I}_{2,j} := \Big\| \sum_{i=1}^{d_k} \Q_{k,i\cdot} \sum_{t=1}^T \E_{t,(k), i\cdot}' \Qmk \F_{Z,t,(k)}' \Q_{k,j\cdot} \Big\|^2 ,\\
    & \c{I}_{3,j} := \Big\| \sum_{i=1}^{d_k} \Q_{k,i\cdot} \sum_{t=1}^T \E_{t,(k), i\cdot}' \E_{t,(k), j\cdot} \Big\|^2 .
\end{align*}

By the Cauchy-Schwarz inequality,
\begin{align*}
    \c{I}_{1,j} &= \Big\| \sum_{i=1}^{d_k} \Q_{k,i\cdot} \sum_{t=1}^T \A_{k,i\cdot}' \F_{t,(k)} \Amk' \E_{t,(k), j\cdot} \Big\|^2
    \leq \Big( \sum_{i=1}^{d_k} \big\|\Q_{k,i\cdot} \big\|^2\Big) \Big\{\sum_{i=1}^{d_k} \Big( \sum_{t=1}^T \A_{k,i\cdot}' \F_{t,(k)} \Amk' \E_{t,(k), j\cdot} \Big)^2 \Big\} \\
    &=
    \big\|\Q_k \big\|_F^2 \cdot \Big\{\sum_{i=1}^{d_k} \|\A_{k,i\cdot} \|^2 \Big( \sum_{h=1}^{\dmk} \sum_{t=1}^T \frac{\A_{k,i\cdot}'}{ \|\A_{k,i\cdot} \| } \F_{t,(k)} \A_{\text{-}k, h\cdot} \E_{t,(k), jh} \Big)^2 \Big\}
    =
    O_P\Big( T d_k^{\delta_{k,1}} \dmk \Big) ,
\end{align*}
where the last equality used Assumptions (L1), (L2) and Lemma~\ref{lemma: correlation_Et_Ft}.

Consider $\c{I}_{2,j}$. From Assumptions (E1) and (E2), we have for any $t\in[T], k\in[K], i\in[j_k], h \in[\dmk]$,
\begin{equation*}
    E_{t,(k), ih} = \sum_{w\geq 0} a_{e,w}
    \A_{e,k,i\cdot}' \X_{e,t-w, (k)} \A_{e,\text{-}k, h\cdot} + \Sigma_{\epsilon, (k), ih} \sum_{w \geq 0} a_{\epsilon,w} X_{\epsilon, t-w,(k), ih} .
\end{equation*}
By Assumptions (F1), (E1) and (E2), we first have
\begin{equation}
\begin{split}
\label{eqn: I2_step1}
    &\hspace{12pt}
    \b{E}\Big\{ \Big( \sum_{t=1}^T \sum_{h=1}^{\dmk} (\sum_{w \geq 0} a_{e,w} \A_{e,k, i\cdot}' \X_{e,t-w, (k)} \A_{e,\text{-}k, h\cdot}) \A_{\text{-}k, h\cdot}' \F_{t, (k)}' \A_{k,j\cdot} \Big)^2 \Big\} \\
    &=
    \text{Cov}\Big( \sum_{t=1}^T \sum_{h=1}^{\dmk} \A_{\text{-}k, h\cdot}' (\sum_{w\geq 0} a_{f,w} \X_{f,t-w, (k)}' ) \A_{k,j\cdot} (\sum_{w \geq 0} a_{e,w} \A_{e,k, i\cdot}' \X_{e,t-w, (k)} \A_{e,\text{-}k, h\cdot}) , \\
    &\hspace{36pt}
    \sum_{t=1}^T \sum_{h=1}^{\dmk} \A_{\text{-}k, h\cdot}' (\sum_{w\geq 0} a_{f,w} \X_{f,t-w, (k)}' ) \A_{k,j\cdot} (\sum_{w \geq 0} a_{e,w} \A_{e,k, i\cdot}' \X_{e,t-w, (k)} \A_{e,\text{-}k, h\cdot}) \Big) \\
    &=
    \sum_{h=1}^{\dmk} \sum_{l=1}^{\dmk} \sum_{t=1}^T
    \sum_{w\geq 0} a_{f,w}^2 a_{e,w}^2 \cdot
    \|\A_{k, j\cdot}\|^2 \cdot
    \| \A_{\text{-}k, h\cdot}\|\cdot
    \| \A_{\text{-}k, l\cdot}\|\cdot
    \|\A_{e,\text{-}k, h\cdot}\|\cdot
    \|\A_{e,\text{-}k, l\cdot}\|\cdot
    \|\A_{e,k, i\cdot}\|^2 \\
    &=
    O(T)\cdot \|\A_{k,j\cdot}\|^2 \cdot \|\A_{e,k, i\cdot}\|^2 .
\end{split}
\end{equation}
Similarly, it holds that
\begin{equation}
\begin{split}
\label{eqn: I2_step2}
    &\hspace{12pt}
    \b{E}\Big\{ \Big\| \sum_{i=1}^{d_k} \sum_{h=1}^{\dmk} \sum_{t=1}^T \Q_{k,i\cdot} \Big( \Sigma_{\epsilon, (k), ih} \sum_{w\geq 0} a_{\epsilon,w} X_{\epsilon, t-w, (k), ih} \Big) \A_{\text{-}k, h\cdot}' \F_{t, (k)}' \A_{k,j\cdot} \Big\|^2 \Big\} \\
    &=
    \text{Cov}\Big( \sum_{i=1}^{d_k} \sum_{h=1}^{\dmk} \sum_{t=1}^T \Q_{k,i\cdot} \Big( \Sigma_{\epsilon, (k), ih} \sum_{w\geq 0} a_{\epsilon,w} X_{\epsilon, t-w, (k), ih} \Big) \A_{\text{-}k, h\cdot}' \F_{t, (k)}' \A_{k,j\cdot} , \\
    &\hspace{36pt}
    \sum_{i=1}^{d_k} \sum_{h=1}^{\dmk} \sum_{t=1}^T \Q_{k,i\cdot} \Big( \Sigma_{\epsilon, (k), ih} \sum_{w\geq 0} a_{\epsilon,w} X_{\epsilon, t-w, (k), ih} \Big) \A_{\text{-}k, h\cdot}' \F_{t, (k)}' \A_{k,j\cdot} \Big) \\
    &=
    \sum_{i=1}^{d_k} \sum_{h=1}^{\dmk} \sum_{t=1}^T \sum_{w\geq 0} a_{f, w}^2 a_{\epsilon, w}^2 \cdot \|\A_{k,j \cdot}\|^2 \cdot \|\A_{\text{-}k, h\cdot}\|^2 \cdot \Sigma_{\epsilon, (k), ih}^2
    \cdot \|\Q_{k,i\cdot}\|^2
    = O(T) \cdot \|\A_{k, j\cdot}\|^2 \cdot \|\Amk \|^2 .
\end{split}
\end{equation}
Hence using Lemma~\ref{lemma: norm_hat_D}, it holds that
\begin{align*}
    \c{I}_{2,j} &= \Big\| \sum_{i=1}^{d_k} \Q_{k,i\cdot} \sum_{t=1}^T \E_{t,(k), i\cdot}' \Amk \F_{t,(k)}' \A_{k,j\cdot} \Big\|^2
    = \Big\| \sum_{i=1}^{d_k} \Q_{k,i\cdot} \sum_{h=1}^{\dmk} \sum_{t=1}^T E_{t,(k), ih} \A_{\text{-}k, h\cdot}' \F_{t,(k)}' \A_{k,j\cdot} \Big\|^2 \\
    &\leq
    2\, \Big\{ \Big\| \sum_{i=1}^{d_k} \sum_{h=1}^{\dmk} \sum_{t=1}^T \Q_{k,i\cdot} \Big( \Sigma_{\epsilon, (k), ih} \sum_{w\geq 0} a_{\epsilon,w} X_{\epsilon, t-w, (k), ih} \Big) \A_{\text{-}k, h\cdot}' \F_{t, (k)}' \A_{k,j\cdot} \Big\|^2 \Big\} \\
    &\hspace{12pt}
    + 2\, \Big\{ \big\|\Q_{k} \big\|^2 \cdot \sum_{i=1}^{d_k} \Big( \sum_{t=1}^T \sum_{h=1}^{\dmk} (\sum_{w \geq 0} a_{e,w} \A_{e,k, i\cdot}' \X_{e,t-w, (k)} \A_{e,\text{-}k, h\cdot}) \A_{\text{-}k, h\cdot}' \F_{t, (k)}' \A_{k,j\cdot} \Big)^2 \Big\} \\
    &=
    O_P\Big( \|\A_{k, j\cdot}\|^2 \cdot T d_k^{-\delta_{k, 1}} g_s \Big)
\end{align*}
where we used \eqref{eqn: I2_step1} and \eqref{eqn: I2_step2} in the last equality.

For $\c{I}_{3,j}$, let $r_{e, \text{-}k}:= \prod_{p\neq k} r_{e,p}$. Then from the noise structure in Assumptions (E1) and (E2),
\begin{equation*}
\begin{split}
    &\hspace{12pt}
    \text{Var}\Big( \sum_{i=1}^{d_k} \sum_{h=1}^{\dmk} \sum_{t=1}^T \Q_{k,i\cdot} E_{t,(k), ih} E_{t,(k), jh} \Big) \\
    &=
    O\Big( \sum_{i=1}^{d_k} \sum_{u=1}^{d_k} \sum_{h=1}^{\dmk} \sum_{l=1}^{\dmk} \sum_{t=1}^T \sum_{n=1}^{r_{e,k}} \sum_{m=1}^{r_{e, \text{-}k}} \sum_{w \geq 0}
    a_{e,w}^4 A_{e,k,in} A_{e,k,un} A_{e,k,jn}^2 A_{e,\text{-}k, hm}^2 A_{e,\text{-}k, lm}^2 \\
    &\hspace{12pt}
    \cdot \|\Q_{k,i\cdot}\|
    \cdot \|\Q_{k,u\cdot}\|
    \cdot \text{Var}(X_{e,t-w, (k), nm}^2) \Big) \\
    &\hspace{12pt}
    + O\Big( \sum_{i=1}^{d_k} \sum_{h=1}^{\dmk} \sum_{t=1}^T \sum_{w\geq 0} a_{\epsilon, w}^4 \Sigma_{\epsilon, (k), ih}^2 \Sigma_{\epsilon, (k), jh}^2
    \cdot \|\Q_{k,i\cdot}\|^2
    \cdot \text{Var}( X_{e,t-w, (k), ih} X_{e,t-w, (k), jh}) \\
    &
    = O(T + T \dmk)=O(T \dmk).
\end{split}
\end{equation*}
Moreover, it also holds that
\begin{align*}
    &\hspace{12pt}
    \b{E}\Big( \sum_{i=1}^{d_k} \Big| \sum_{h=1}^{\dmk} \sum_{t=1}^T E_{t,(k), ih} E_{t,(k), jh} \Big|\Big) \\
    &=
    \sum_{i=1}^{d_k} \Big| \sum_{h=1}^{\dmk} \sum_{t=1}^T \Big( \|\A_{e,\text{-}k, h\cdot}\|^2 \cdot \|\A_{e,k,i\cdot}\| \cdot \|\A_{e,k,j\cdot}\| + \Sigma_{\epsilon, (k), ih} \b{1}_{\{i=j\}} \Big) \Big|
    = O(T \dmk),
\end{align*}
together with $\max_i\|\Q_{k, i\cdot}\|^2 \leq \|\A_{k,j\cdot}\|^2 \cdot \|\Z_k^{-1/2}\|^2 = O_P\big( d_k^{-\delta_{k,r_k}} \big)$, we arrive at
\[
\c{I}_{3,j} = \Big\| \sum_{i=1}^{d_k} \sum_{h=1}^{\dmk} \Q_{k, i\cdot} \sum_{t=1}^T E_{t,(k), ih} E_{t,(k), jh} \Big\|^2
= O_P\Big(T\dmk + T^2 d_k^{-\delta_{k,r_k}} \dmk^2 \Big).
\]

Finally, for \eqref{eqn: Q_norm_rate2} we have
\begin{align*}
    \big\| \wh\Q_k - \Q_k\H_k' \big\|_F^2 &=
    O_P\Big\{ T^{-2} d_k^{2(\delta_{k,1} -\delta_{k,r_k})} g_s^{-2} \cdot \Big( Td d_k^{\delta_{k,1}} + T^2 d_k^{1-\delta_{k,r_k}} \dmk^2 \Big) \Big\} + o_P\Big( \big\|\wh\Q_{k} -\Q_{k} \H_k' \big\|_F^2 \Big) \\
    &=
    O_P\Bigg\{ d_k^{2(\delta_{k,1} - \delta_{k,r_k})} \Bigg( \frac{1}{T\dmk d_k^{1- \delta_{k,1}}} +\frac{1}{ d_k^{1+ \delta_{k, r_k}}} \Bigg) \frac{d^2}{g_s^2} \Bigg\}.
\end{align*}
This completes the proof of Lemma~\ref{lemma: Q_hat_rate}. $\square$

\begin{lemma}\label{lemma: hat_rate_cFcC}
(Consistency of $\wh\cC_{t}$)
Let Assumptions in Lemma~\ref{lemma: Q_hat_rate} hold. With $\{\Z_k\}_{k\in[K]}$ from Assumptions (L1) and (L2) and $\{\H_k\}_{k\in[K]}$ from the statement of Lemma~\ref{lemma: Q_hat_rate}, define
\begin{align*}
    \Z_\otimes &:= \Z_K\otimes \dots \otimes \Z_1 , \;\;\;
    \H_\otimes := \H_K\otimes \dots \otimes \H_1 , \;\;\;
    \vec{\wh\cF_{t}} := \big(\wh\Q_{K} \otimes \dots \otimes \wh\Q_{1} \big)' \vec{\cY_t} .
\end{align*}
Then the estimators of the vectorised (renormalised) core factor and the $i$-th entry of the vectorised common component for \eqref{eqn: test_fm_w_Kron} are consistent such that
\begin{align}
    &\hspace{12pt}
    \big\|\vec{\wh\cF_{t}} - (\H_\otimes')^{-1} \vec{\cF_{Z,t}} \big\|^2 =\big\|\vec{\wh\cF_{t}} - (\H_\otimes')^{-1} \Z_\otimes^{1/2} \vec{\cF_t} \big\|^2 \notag \\
    &= O_P\Bigg( \max_{k\in[K]}\Bigg\{ d_k^{2(\delta_{k,1} - \delta_{k,r_k})} \Bigg( \frac{1}{T\dmk d_k^{1- \delta_{k,1}}} +\frac{1}{ d_k^{1+ \delta_{k, r_k}}} \Bigg) \frac{d^2}{g_s} \Bigg\} +\frac{d}{g_w} \Bigg) ,
    \label{eqn: hat_rate_cF} \\
    &\hspace{12pt}
    \Big\{ \big(\vec{\wh\cC_{t}}\big)_i - \big( \vec{\cC_{t}} \big)_i \Big\}^2 \notag \\
    &= O_P\Bigg( \max_{k\in[K]}\Bigg\{ d_k^{2(\delta_{k,1} - \delta_{k,r_k})} \Bigg( \frac{1}{T\dmk d_k^{1- \delta_{k,1}}} +\frac{1}{ d_k^{1+ \delta_{k, r_k}}} \Bigg) \frac{d^2}{g_s g_w} \Bigg\} +\frac{d}{g_w^2} \Bigg) .
    \label{eqn: hat_rate_cC}
\end{align}
\end{lemma}

\textbf{\textit{Proof of Lemma~\ref{lemma: hat_rate_cFcC}.}}
By \eqref{eqn: test_fm_w_Kron}, we have
\begin{equation}
\label{eqn: cFt_hat_decomp}
\begin{split}
    &\hspace{12pt}
    \vec{\wh\cF_{t}} - (\H_\otimes')^{-1} \Z_\otimes^{1/2} \vec{\cF_t}
    = \big(\wh\Q_{K} \otimes \dots \otimes \wh\Q_{1} \big)' \vec{\cY_t} - (\H_\otimes')^{-1} \Z_\otimes^{1/2} \vec{\cF_t} \\
    &=
    \big(\wh\Q_{K} \otimes \dots \otimes \wh\Q_{1} \big)' \vec{\cF_t \times_{k=1}^K \A_k + \cE_t} - (\H_\otimes')^{-1} \Z_\otimes^{1/2} \vec{\cF_t} \\
    &=
    \big(\wh\Q_{K} \otimes \dots \otimes \wh\Q_{1} \big)' \Big\{ \big(\Q_{K}\H_K' \otimes \dots \otimes \Q_{1}\H_1' \big) - \big(\wh\Q_{K} \otimes \dots \otimes \wh\Q_{1} \big) \Big\} (\H_\otimes')^{-1} \Z_\otimes^{1/2} \vec{\cF_t} \\
    &\hspace{12pt}
    + \big\{\big(\wh\Q_{K} \otimes \dots \otimes \wh\Q_{1} \big) - \big(\Q_{K}\H_K' \otimes \dots \otimes \Q_{1}\H_1' \big) \big\}' \vec{\cE_t} \\
    &\hspace{12pt}
    + \big(\Q_{K}\H_K' \otimes \dots \otimes \Q_{1}\H_1' \big)' \vec{\cE_t}
    =: \c{I}_{f,1} + \c{I}_{f,2} + \c{I}_{f,3}.
\end{split}
\end{equation}

We first show by an induction argument that for any positive integer $K$,
\begin{equation}
\label{eqn: Q_kron_induction}
\big\|\big( \wh\Q_{K} \otimes \dots \otimes \wh\Q_{1} \big) - \big(\Q_{K}\H_K' \otimes \dots \otimes \Q_{1}\H_1' \big) \big\|_F = O_P\Big( \max_{k\in[K]} \big\| \wh\Q_{k} - \Q_{k}\H_k' \big\|_F\Big).
\end{equation}
The initial case for $K=1$ is trivial. Suppose \eqref{eqn: Q_kron_induction} holds for $K-1$, then
\begin{align*}
    &\hspace{13pt}
    \big\|\big( \wh\Q_{K} \otimes \dots \otimes \wh\Q_{1} \big) - \big(\Q_{K}\H_K' \otimes \dots \otimes \Q_{1}\H_1' \big) \big\|_F \\
    &=
    \big\| \big( \wh\Q_{K} -\Q_{K}\H_K'\big) \otimes \big( \wh\Q_{K-1} \otimes \dots \otimes \wh\Q_{1} \big) \\
    &\hspace{12pt}
    + \Q_{K}\H_K' \otimes \big\{\big( \wh\Q_{K-1} \otimes \dots \otimes \wh\Q_{1} \big) - \big(\Q_{K-1} \H_{K-1}' \otimes \dots \otimes \Q_{1}\H_1' \big) \big\} \big\|_F \\
    &=
    O_P\Big( \big\| \wh\Q_{K} - \Q_{K}\H_K' \big\|_F + \big\|\big( \wh\Q_{K-1} \otimes \dots \otimes \wh\Q_{1} \big) - \big(\Q_{K-1}\H_{K-1}' \otimes \dots \otimes \Q_{1}\H_1' \big) \big\|_F \Big),
\end{align*}
which concludes \eqref{eqn: Q_kron_induction}. Hence for $\c{I}_{f,1}$, with Lemma~\ref{lemma: Q_hat_rate} we immediately have
\begin{equation}
\label{eqn: I_f1}
\begin{split}
    \|\c{I}_{f,1} \|^2 &= \Big\| \big(\wh\Q_{K} \otimes \dots \otimes \wh\Q_{1} \big)' \Big\{ \big(\Q_{K}\H_K' \otimes \dots \otimes \Q_{1}\H_1' \big) - \big(\wh\Q_{K} \otimes \dots \otimes \wh\Q_{1} \big) \Big\} (\H_\otimes')^{-1} \Z_\otimes^{1/2} \vec{\cF_t} \Big\|^2 \\
    &=
    O_P\Big( \big\|\big( \wh\Q_{K} \otimes \dots \otimes \wh\Q_{1} \big) - \big(\Q_{K}\H_K' \otimes \dots \otimes \Q_{1}\H_1' \big) \big\|_F^2 \cdot \big\| \Z_\otimes^{1/2} \vec{\cF_t} \big\|^2 \Big) \\
    &=
    O_P\Bigg( \max_{k\in[K]} \Bigg\{ d_k^{2(\delta_{k,1} - \delta_{k,r_k})} \Bigg( \frac{1}{T\dmk d_k^{1- \delta_{k,1}}} +\frac{1}{ d_k^{1+ \delta_{k, r_k}}} \Bigg) \frac{d^2}{g_s} \Bigg\} \Bigg),
\end{split}
\end{equation}
where the last equality used \eqref{eqn: Q_kron_induction}, Assumptions (F1), (L1) and (L2).

For $\c{I}_{f,2}$, observe first throughout the proof of Lemma~\ref{lemma: Q_hat_rate}, the consistency of $\wh\Q_{k,j\cdot}$ for $k\in[K], j\in[d_k]$ can be shown with the same argument (omitted), before eventually being aggregated over all $d_k$ rows. That is,
\begin{equation}
\label{eqn: Q_hat_rate_row}
\big\| \wh\Q_{k,j\cdot} - \H_k\Q_{k,j\cdot} \big\|^2 = O_P\Bigg\{ d_k^{2(\delta_{k,1} - \delta_{k,r_k})-1} \Bigg( \frac{1}{T\dmk d_k^{1- \delta_{k,1}}} +\frac{1}{ d_k^{1+ \delta_{k, r_k}}} \Bigg) \frac{d^2}{g_s^2} \Bigg\}.
\end{equation}
Then we have $\norm{\wh\Q_{k,j\cdot}}^2 = O_P\big(\norm{\wh\Q_{k,j\cdot} - \H_k\Q_{k,j\cdot}}^2 + \norm{\H_k\Z_{k}^{-1/2} \A_{k,j\cdot}}^2 \big) = O_P\big( d_k^{-\delta_{k,r_k}} \big)$ which used \eqref{eqn: Q_hat_rate_row}, Assumptions (L1) (or (L2)) and (R1). Note that this rate $d_k^{-\delta_{k,r_k}}$ is the same as the one for $\norm{\Q_{k,j\cdot}}^2$, shown in the proof of Lemma~\ref{lemma: Q_hat_rate}. Moreover, it holds that for any positive integer $k$ that
\begin{equation}
\label{eqn: Q_kron_row_induction}
\begin{split}
    &\hspace{12pt}
    \max_{\ell\in[d]} \Big\| \big[\big( \wh\Q_{K} \otimes \dots \otimes \wh\Q_{1} \big) - \big(\Q_{K}\H_K' \otimes \dots \otimes \Q_{1}\H_1' \big) \big]_{\ell \cdot} \Big\|^2 \\
    &=
    O\Bigg( \max_{k\in[K]} \Bigg\{ \max_{j\in[d_k]} \big\| \wh\Q_{k,j\cdot} - \H_k\Q_{k,j\cdot} \big\|^2 \cdot \prod_{ j\in[K] \setminus \{k\} } \max_{i\in[d_j]} \norm{\wh\Q_{j,i\cdot}}^2 \Bigg\} \Bigg),
\end{split}
\end{equation}
which can be shown by an induction argument for which the initial case for $K=1$ is trivial, and the induction step is seen by
\begin{align*}
    &\hspace{13pt}
    \max_{\ell\in[d]} \Big\| \big[\big( \wh\Q_{K} \otimes \dots \otimes \wh\Q_{1} \big) - \big(\Q_{K}\H_K' \otimes \dots \otimes \Q_{1}\H_1' \big) \big]_{\ell \cdot} \Big\|^2 \\
    &=
    O\Bigg(\max_{j\in[d_K]} \big\| \wh\Q_{K,j\cdot} - \H_K\Q_{K,j\cdot} \big\|^2 \cdot \prod_{k=1}^{K-1} \max_{i\in[d_k]} \norm{\wh\Q_{k,i\cdot}}^2 \\
    &\hspace{12pt}
    + \max_{i\in[d_K]} \norm{\Q_{K,i\cdot}}^2 \max_{i\in \big[\prod_{k=1}^{K-1} d_k \big]}\Big\| \big[\big( \wh\Q_{K-1} \otimes \dots \otimes \wh\Q_{1} \big) - \big(\Q_{K-1}\H_{K-1}' \otimes \dots \otimes \Q_{1}\H_1' \big) \big]_{i \cdot} \Big\|^2 \Bigg) .
\end{align*}
Hence, for $\c{I}_{f,2}$ we have
\begin{equation}
\label{eqn: I_f2}
\begin{split}
    \|\c{I}_{f,2} \|^2 &= \Big\| \big\{\big(\wh\Q_{K} \otimes \dots \otimes \wh\Q_{1} \big) - \big(\Q_{K}\H_K' \otimes \dots \otimes \Q_{1}\H_1' \big) \big\}' \vec{\cE_t} \Big\|^2 \\
    &=
    \max_{\ell\in[d]} \Big\| \big[\big( \wh\Q_{K} \otimes \dots \otimes \wh\Q_{1} \big) - \big(\Q_{K}\H_K' \otimes \dots \otimes \Q_{1}\H_1' \big) \big]_{\ell \cdot} \Big\|^2 \sum_{j,l=1}^d  \Big| \b{E}\big(\vec{\cE_t} \big)_j \big(\vec{\cE_t} \big)_l \Big| \\
    &=
    O_P\Bigg( \max_{k\in[K]} \Bigg\{ d_k^{2\delta_{k,1} -\delta_{k,r_k}} \Bigg( \frac{1}{T\dmk d_k^{2- \delta_{k,1}}} +\frac{1}{ d_k^{2+ \delta_{k, r_k}}} \Bigg) \frac{d^3}{g_s^2 g_w} \Bigg\} \Bigg) ,
\end{split}
\end{equation}
where the last equality used \eqref{eqn: Q_hat_rate_row}, \eqref{eqn: Q_kron_row_induction} and the first result on Lemma~\ref{lemma: correlation_Et_Ft}.1.

Lastly, consider $\c{I}_{f,3}$. With Assumptions (L1) and (L2), we have
\begin{equation}
\label{eqn: I_f3}
\begin{split}
    \|\c{I}_{f,3} \|^2 &= \Big\| \big(\Q_{K}\H_K' \otimes \dots \otimes \Q_{1}\H_1' \big)' \vec{\cE_t} \Big\|^2
    = O_P\Big( \big\| \big(\Q_{K}' \otimes \dots \otimes \Q_{1}' \big) \vec{\cE_t} \big\|^2 \Big) \\
    &=
    O_P\Big(\big\| \Z_\otimes^{-1/2} \big\|_F^2 \cdot \big\| \big(\A_{K}' \otimes \dots \otimes \A_{1}' \big) \vec{\cE_t} \big\|^2 \Big) \\
    &=
    O_P\Big\{ g_w^{-1} \cdot \sum_{j=1}^d \big\| \big(\vec{\cE_t} \big)_j \big(\A_{K} \otimes \dots \otimes \A_{1} \big)_{j\cdot} \big\|^2 \Big\} = O_P(d /g_w) ,
\end{split}
\end{equation}
where the last equality used (L1) and (L2) and the first result on Lemma~\ref{lemma: correlation_Et_Ft}.1 that
\begin{align*}
    &\hspace{12pt}
    \b{E}\Big\{ \sum_{j=1}^d \big\| \big(\vec{\cE_t} \big)_j \big(\A_{K} \otimes \dots \otimes \A_{1} \big)_{j\cdot} \big\|^2 \Big\} \\
    &\leq
    \max_{j\in[d]} \big\|\big(\A_{K} \otimes \dots \otimes \A_{1} \big)_{j\cdot} \big\|^2 \sum_{j,l=1}^d  \Big| \b{E}\big(\vec{\cE_t} \big)_j \big(\vec{\cE_t} \big)_l \Big| = O(d).
\end{align*}

Combining \eqref{eqn: cFt_hat_decomp}, \eqref{eqn: I_f1}, \eqref{eqn: I_f2} and \eqref{eqn: I_f3}, we obtain \eqref{eqn: hat_rate_cF}.

It remains to show \eqref{eqn: hat_rate_cC}. To this end, from \eqref{eqn: test_fm_w_Kron} and \eqref{eqn: hat_cCt} (where $\wh\cC_{m,t}$ is simplified as $\wh\cC_{t}$, explained in the summary of proofs), we have
\begin{align*}
    &\hspace{12pt}
    \vec{\wh\cC_{t}} - \vec{\cC_{t}} = (\wh\Q_{K} \otimes \dots \otimes \wh\Q_{1} ) \vec{\wh\cF_{t}} - \vec{\cF_t \times_{k=1}^K \A_k} \\
    &=
    (\wh\Q_{K} \otimes \dots \otimes \wh\Q_{1} ) \vec{\wh\cF_{t}} - (\Q_{K}\H_K' \otimes \dots \otimes \Q_{1}\H_1' ) (\H_\otimes')^{-1} \Z_\otimes^{1/2} \vec{\cF_t} \\
    &=
    (\wh\Q_{K} \otimes \dots \otimes \wh\Q_{1} ) \Big\{ \vec{\wh\cF_{t}} -(\H_\otimes')^{-1} \Z_\otimes^{1/2} \vec{\cF_t} \Big\} \\
    &\hspace{12pt}
    + \Big\{ (\wh\Q_{K} \otimes \dots \otimes \wh\Q_{1} ) - (\Q_{K}\H_K' \otimes \dots \otimes \Q_{1}\H_1' ) \Big\} (\H_\otimes')^{-1} \Z_\otimes^{1/2} \vec{\cF_t} ,
\end{align*}
which implies that
\begin{align*}
    &\hspace{12pt}
    \Big\{ \big(\vec{\wh\cC_{t}}\big)_i - \big( \vec{\cC_{t}} \big)_i \Big\}^2 \\
    &=
    O_P\Bigg( \prod_{k=1}^K \max_{j\in[d_k]} \norm{\wh\Q_{k,j\cdot}}^2 \cdot \big\|\vec{\wh\cF_{t}} - (\H_\otimes')^{-1} \Z_\otimes^{1/2} \vec{\cF_t} \big\|^2 \\
    &\hspace{36pt}
    + \max_{\ell\in[d]} \Big\| \big[\big( \wh\Q_{K} \otimes \dots \otimes \wh\Q_{1} \big) - \big(\Q_{K}\H_K' \otimes \dots \otimes \Q_{1}\H_1' \big) \big]_{\ell \cdot} \Big\|^2 \cdot \big\|\Z_\otimes^{1/2}\big\|_F^2 \Bigg) \\
    &=
    O_P\Bigg( \max_{k\in[K]}\Bigg\{ d_k^{2(\delta_{k,1} - \delta_{k,r_k})} \Bigg( \frac{1}{T\dmk d_k^{1- \delta_{k,1}}} +\frac{1}{ d_k^{1+ \delta_{k, r_k}}} \Bigg) \frac{d^2}{g_s g_w} \Bigg\} +\frac{d}{g_w^2} \Bigg) ,
\end{align*}
where the last equality used \eqref{eqn: hat_rate_cF}, \eqref{eqn: Q_hat_rate_row}, \eqref{eqn: Q_kron_row_induction}, Assumptions (L1) and (L2). This completes the proof of Lemma~\ref{lemma: hat_rate_cFcC}. $\square$

\begin{lemma}\label{lemma: Q_tilde_rate}
(Consistency of $\{\wt\Q_j\}_{j\in[v-1]}$, $\wt\Q_V$ and $\wt\cC_{\textnormal{reshape},t}$)
Let Assumptions (F1), (L1), (E1), (E2) and (R1) hold, and consider the model \eqref{eqn: test_fm_wo_Kron}. For $j\in[v-1]$, define $\wt\D_j$ as the $r_j \times r_j$ diagonal matrix with the first largest $r_j$ eigenvalues of
\[
\frac{1}{T} \sum_{t=1}^T \Reshape(\cY_t, \cA)_{(j)} \Reshape(\cY_t, \cA)_{(j)}'
\]
on the main diagonal, such that $\wt\D_j = \wt\Q_j' \big(T^{-1} \sum_{t=1}^T \Reshape(\cY_t, \cA)_{(j)} \Reshape(\cY_t, \cA)_{(j)}' \big) \wt\Q_j$. Similarly, $\wt\D_V$ denotes the $r_V \times r_V$ diagonal matrix with the first largest $r_V$ eigenvalues of
\[
\frac{1}{T} \sum_{t=1}^T \Reshape(\cY_t, \cA)_{(v)} \Reshape(\cY_t, \cA)_{(v)}' .
\]
Correspondingly, define an $r_j\times r_j$ matrix and an $r_V\times r_V$ matrix that
\begin{align*}
    \wt\H_j &:= T^{-1} \wt\D_j^{-1} \wt\Q_j' \Q_j \sum_{t=1}^T \F_{\textnormal{reshape},Z,t,(j)} \F_{\textnormal{reshape},Z,t,(j)}' , \\
    \wt\H_V &:= T^{-1} \wt\D_V^{-1} \wt\Q_V' \Q_V \sum_{t=1}^T \F_{\textnormal{reshape},Z,t,(v)} \F_{\textnormal{reshape},Z,t,(v)}' ,
\end{align*}
where $\cF_{\textnormal{reshape},Z,t} := \cF_{\textnormal{reshape},t} \times_{j=1}^{v-1} \Z_j^{1/2} \times_v \Z_V^{1/2}$. As $T,d_1,\dots, d_{v-1}, d_V \to\infty$ we have $\{\wt\H_j\}_{j\in[v-1]}$, $\wt\H_V$ are invertible, and for $j\in[v-1]$ that $\|\wt\H_j\|_F =O_P(1)$ and $\|\wt\H_V\|_F=O_P(1)$. For each $j\in[v-1]$,
\begin{align}
    \big\| \wt\Q_j - \Q_j\wt\H_j' \big\|_F^2 &=
    O_P\Bigg\{ d_j^{2(\delta_{j,1} - \delta_{j,r_j})} \Bigg( \frac{1}{T\dmk d_j^{1- \delta_{j,1}}} +\frac{1}{ d_j^{1+ \delta_{j, r_j}}} \Bigg) \frac{d^2}{\gamma_s^2} \Bigg\} ,
    \label{eqn: Qj_tilde_rate} \\
    \big\| \wt\Q_V - \Q_V\wt\H_V' \big\|_F^2 &=
    O_P\Bigg\{ d_V^{2(\delta_{V,1} - \delta_{V,r_V})} \Bigg( \frac{1}{Td d_V^{- \delta_{V,1}}} +\frac{1}{ d_V^{1+ \delta_{V, r_V}}} \Bigg) \frac{d^2}{\gamma_s^2} \Bigg\} ,
    \label{eqn: QV_tilde_rate}
\end{align}
where $\gamma_s$ is defined in Assumption (R1). Lastly, the $i$-th entry of the vectorised common component in \eqref{eqn: test_fm_wo_Kron} is also consistent such that
\begin{equation}
\label{eqn: cC_tilde_reshape_rate}
\begin{split}
    &\hspace{12pt}
    \Big\{ \big(\vec{\wt\cC_{\textnormal{reshape},t}}\big)_i - \big( \vec{\cC_{\textnormal{reshape},t}} \big)_i \Big\}^2 \\
    &=
    O_P\Bigg( \max_{j\in[v-1]}\Bigg\{ d_j^{2(\delta_{j,1} - \delta_{j,r_j})} \Bigg( \frac{1}{T\dmk d_j^{1- \delta_{j,1}}} +\frac{1}{ d_j^{1+ \delta_{j, r_j}}} \Bigg) \frac{d^2}{\gamma_s \gamma_w} \Bigg\} \\
    &\hspace{57pt}
    + d_V^{2(\delta_{V,1} - \delta_{V,r_V})} \Bigg( \frac{1}{Td d_V^{- \delta_{V,1}}} +\frac{1}{ d_V^{1+ \delta_{V, r_V}}} \Bigg) \frac{d^2}{\gamma_s \gamma_w}
    + \frac{d}{\gamma_w^2} \Bigg) ,
\end{split}
\end{equation}
where $\gamma_w$ is defined in Assumption (R2).
\end{lemma}

\textbf{\textit{Proof of Lemma~\ref{lemma: Q_tilde_rate}.}}
Consider \eqref{eqn: cEt_reshape}, we have by Assumptions (E1) and (E2) that
\begin{align*}
    \Reshape(\c{F}_{e,t}, \cA) = \sum_{q\geq 0} a_{e,q} \Reshape(\c{X}_{e,t-q}, \cA) , \;\;\;
    \Reshape(\bepsilon_t, \cA) = \sum_{q\geq 0} a_{\epsilon,q} \Reshape(\c{X}_{\epsilon, t-q}, \cA) ,
\end{align*}
which implies that the structure depicted in (E1) and (E2) for the noise $\cE_t$ in \eqref{eqn: test_fm_w_Kron} holds for $\cE_{\textnormal{reshape},t}$ in \eqref{eqn: test_fm_wo_Kron}. Read $\Reshape(\cY_t, \cA)$ as an order-$v$ tensor and consider the factor model \eqref{eqn: test_fm_wo_Kron}. Statements \eqref{eqn: Qj_tilde_rate} and \eqref{eqn: QV_tilde_rate} can be shown in exactly the same way as Lemma~\ref{lemma: Q_hat_rate} (without (L2) now; details omitted here), given that all rate conditions used in the proof of Lemma~\ref{lemma: Q_hat_rate} are fulfilled for \eqref{eqn: test_fm_wo_Kron}. In other words, it remains to show the rate conditions similar to the last equality of \eqref{eqn: Q_norm_rate1} are satisfied. Note that for $j\in[v-1]$, \eqref{eqn: Qj_tilde_rate} requires the same rate condition as Lemma~\ref{lemma: Q_hat_rate}. Hence we are left with the rate conditions for \eqref{eqn: QV_tilde_rate}, i.e.,
\[
d \gamma_s^{-2} T^{-1} d_V^{2(\delta_{V,1} - \delta_{V,r_V})+1} = o(1), \;\;\;
d \gamma_s^{-1} d_V^{\delta_{V,1} -\delta_{V,r_V}-1/2} = o(1) ,
\]
which are included in Assumption (R1). With \eqref{eqn: Qj_tilde_rate} and \eqref{eqn: QV_tilde_rate} shown, \eqref{eqn: cC_tilde_reshape_rate} follows similarly as Lemma~\ref{lemma: hat_rate_cFcC} (omitted). The proof of Lemma~\ref{lemma: Q_tilde_rate} is then complete. $\square$

\subsection{Proof of theorems}

\textbf{\textit{Proof of Theorem~\ref{thm: reshape}.}}
Let $\cA=\{a_1, \dots, a_\ell\}$ be given. We first show that the Tucker-decomposition tensor factor model \eqref{eqn: fm_w_Kron} implies \eqref{eqn: fm_wo_Kron} with $\A_{\textnormal{reshape}, K-\ell+1} \in \cK_{d_{a_1} \times \dots \times d_{a_\ell}}$. Suppose $\cY_t = \cF_t \times_{k=1}^K \A_k +\cE_t$.

Consider first $\ell=1$. With any $\cA=\{a_1\}$ (hence the corresponding $\cA^\ast = [K]\setminus \cA$ as defined in Theorem~\ref{thm: reshape}), it is direct that
\begin{align*}
    &\hspace{12pt}
    \Reshape\big( \cY_t ,\cA \big) = \Reshape\big( \cF_t \times_{k=1}^K \A_k +\cE_t, \cA \big) = \Reshape\big( \cF_t \times_{k=1}^K \A_k, \{a_1\} \big) + \Reshape\big( \cE_t, \cA \big) \\
    &=
    \Fold_K\big(\mat{a_1}{\cF_t \times_{k=1}^K \A_k}, \{d_1,\dots, d_{a_1-1}, d_{a_1+1}, \dots, d_K, d_{a_1}\} \big) + \Reshape\big( \cE_t, \cA \big) \\
    &=
    \Fold_K\big( \A_{a_1} \mat{a_1}{\cF_t} \A_{\text{-}a_1}', \{ d_1,\dots, d_{a_1-1}, d_{a_1+1}, \dots, d_K, d_{a_1}\} \big) + \Reshape\big( \cE_t, \cA \big), \\
    &=
    \Reshape\big( \cF_t ,\cA \big) \times_{i=1}^{K-1} \A_{\cA_i^\ast} \times_K \A_{a_1} + \Reshape\big( \cE_t, \cA \big),
\end{align*}
where the last equality used the fact that
\begin{equation}
\label{eqn: fold_reshape_equality}
\begin{split}
    &\hspace{12pt}
    \text{mat}_K\Big\{ \Reshape\big( \cF_t ,\cA \big) \times_{i=1}^{K-1} \A_{\cA_i^\ast} \times_K \A_{a_1} \Big\} = \A_{a_1} \text{mat}_K\Big\{ \Reshape\big( \cF_t ,\{a_1\} \big) \Big\} \A_{\text{-}a_1}' \\
    &=
    \A_{a_1} \text{mat}_K\Big\{ \Fold_K\big(\mat{a_1}{\cF_t}, \{d_1,\dots, d_{a_1-1}, d_{a_1+1}, \dots, d_K, d_{a_1}\} \big) \Big\} \A_{\text{-}a_1}'
    = \A_{a_1} \mat{a_1}{\cF_t} \A_{\text{-}a_1}' .
\end{split}
\end{equation}
Hence, $\cY_t$ follows \eqref{eqn: fm_wo_Kron} with variables defined according to Theorem~\ref{thm: reshape}, and the model has a Kronecker structure product since $\A_{a_1} \in\cK_{d_{a_1}}$.

Consider $\ell=2$ (implying at least $K=2$). In the following, we use the neater notation that mode-$k$ unfolding of some tensor $\cX$ is $\X_{(k)}$. Without loss of generality, let $\cA=\{a,b\}$ (with $a<b$) and the corresponding $\cA^\ast = [K]\setminus \cA$. Take the mode-$b$ unfolding on each $\cY_t$, we obtain
\begin{align*}
    \Y_{t,(b)} &= \A_b \F_{t,(b)} \big(\A_{K}\otimes \dots \otimes \A_{b+1} \otimes \A_{b-1} \otimes \dots\otimes \A_a \otimes \dots\otimes\A_1 \big)' + \E_{t,(b)} ,
\end{align*}
then for each row (as a column vector) of $\Y_{t,(b)}$, we fold it back to an order-$(K-1)$ tensor along the remaining dimensions, i.e., for any $i$-th row $\Y_{t,(b), i\cdot}$ with $i\in[d_b]$,
\begin{align*}
    &\hspace{12pt}
    \cY_{t,(b),i} := \Fold\big(\Y_{t,(b), i\cdot} , \{d_1, \dots, d_{b-1}, d_{b+1}, \dots, d_K\}\big) \\
    &=
    \Fold\big( \big(\A_{K}\otimes \dots \otimes \A_{b+1} \otimes \A_{b-1} \otimes \dots\otimes \A_a \otimes  \dots\otimes\A_1 \big) \F_{t,(b)}' \A_{b,i\cdot} , \{d_1, \dots, d_{b-1}, d_{b+1}, \dots, d_K\}\big) \\
    &\hspace{12pt}
    + \Fold\big(\E_{t,(b), i\cdot} , \{d_1, \dots, d_{b-1}, d_{b+1}, \dots, d_K\}\big) \\
    &=
    \Fold\big( \F_{t,(b)}' \A_{b,i\cdot} , \{r_1, \dots, r_{b-1}, r_{b+1}, \dots, r_K\}\big) \times_{k=1}^{b-1} \A_k \times_{h=b}^{K-1} \A_{h+1} \\
    &\hspace{12pt}
    + \Fold\big(\E_{t,(b), i\cdot} , \{d_1, \dots, d_{b-1}, d_{b+1}, \dots, d_K\}\big) .
\end{align*}
Define $\A_{\text{-}b,\text{-}a} := \A_{K}\otimes \dots \otimes \A_{b+1} \otimes \A_{b-1} \otimes \dots\otimes \A_{a+1} \otimes \A_{a-1} \otimes \dots\otimes\A_1$, where by convention $\A_{\text{-}b,\text{-}a} = 1$ if $K=2$. Take the mode-$a$ unfolding on $\cY_{t,(b),i}$, we have
\begin{align*}
    &\hspace{12pt}
    \big(\cY_{t,(b),i}\big)_{(a)} \\
    &=
    \A_a \big\{\Fold\big( \F_{t,(b)}' \A_{b,i\cdot} , \{r_1, \dots, r_{a-1}, r_{a+1}, \dots, r_K\} \big) \big\}_{(a)} \A_{\text{-}b,\text{-}a}' \\
    &\hspace{12pt}
    + \big\{ \Fold\big(\E_{t,(b), i\cdot} , \{d_1, \dots, d_{b-1}, d_{b+1}, \dots, d_K\} \big) \big\}_{(a)} \\
    &=
    \sum_{j=1}^{r_b} \A_a \big\{\Fold\big( A_{b,ij} \F_{t,(b), j\cdot} , \{r_1, \dots, r_{b-1}, r_{b+1}, \dots, r_K\}\big) \big\}_{(a)} \A_{\text{-}b,\text{-}a}' \\
    &\hspace{12pt}
    + \big\{ \Fold\big(\E_{t,(b), i\cdot} , \{d_1, \dots, d_{b-1}, d_{b+1}, \dots, d_K\}\big) \big\}_{(a)} \\
    &=
    \big(\A_{b,i\cdot}' \otimes \I_{d_a} \big) \big(\I_{r_b} \otimes \A_a \big)
    \begin{pmatrix}
    \big\{\Fold\big( \F_{t,(b),1\cdot} , \{r_1, \dots, r_{b-1}, r_{b+1}, \dots, r_K\}\big) \big\}_{(a)} \A_{\text{-}b,\text{-}a}' \\
    \ldots \\
    \big\{\Fold\big( \F_{t,(b),r_b\cdot} , \{r_1, \dots, r_{b-1}, r_{b+1}, \dots, r_K\} \big) \big\}_{(a)} \A_{\text{-}b,\text{-}a}'
    \end{pmatrix}  \\
    &\hspace{12pt}
    + \big\{ \Fold\big(\E_{t,(b), i\cdot} , \{d_1, \dots, d_{b-1}, d_{b+1}, \dots, d_K\}\big) \big\}_{(a)} \\
    &=
    \big(\A_{b,i\cdot}' \otimes \A_a \big)
    \begin{pmatrix}
    \big\{\Fold\big( \F_{t,(b),1\cdot} , \{r_1, \dots, r_{b-1}, r_{b+1}, \dots, r_K\}\big) \big\}_{(a)} \\
    \ldots \\
    \big\{\Fold\big( \F_{t,(b),r_b\cdot} , \{r_1, \dots, r_{b-1}, r_{b+1}, \dots, r_K\}\big) \big\}_{(a)}
    \end{pmatrix}
    \A_{\text{-}b,\text{-}a}' \\
    &\hspace{12pt}
    + \big\{ \Fold\big(\E_{t,(b), i\cdot} , \{d_1, \dots, d_{b-1}, d_{b+1}, \dots, d_K\}\big) \big\}_{(a)} .
\end{align*}
Therefore,
\begin{align*}
    \cY_{t, a\sim b} = \begin{pmatrix}
    \big(\cY_{t,(b),1}\big)_{(a)} \\
    \ldots \\
    \big(\cY_{t,(b),d_b}\big)_{(a)}
    \end{pmatrix}
    &=
    \big(\A_b \otimes \A_a \big)
    \begin{pmatrix}
    \big\{\Fold\big( \F_{t,(b),1\cdot} , \{r_1, \dots, r_{b-1}, r_{b+1}, \dots, r_K\}\big) \big\}_{(a)} \\
    \ldots \\
    \big\{\Fold\big( \F_{t,(b),r_b\cdot} , \{r_1, \dots, r_{b-1}, r_{b+1}, \dots, r_K\}\big) \big\}_{(a)}
    \end{pmatrix}
    \A_{\text{-}b,\text{-}a}' \\
    &\hspace{12pt}
    + \begin{pmatrix}
    \big\{ \Fold\big(\E_{t,(b), 1\cdot} , \{d_1, \dots, d_{b-1}, d_{b+1}, \dots, d_K\}\big) \big\}_{(a)} \\
    \ldots \\
    \big\{ \Fold\big(\E_{t,(b), d_b\cdot} , \{d_1, \dots, d_{b-1}, d_{b+1}, \dots, d_K\}\big) \big\}_{(a)}
    \end{pmatrix},
\end{align*}
so that by definition of the reshape operator with $\ell=2$,
\begin{align*}
    &\hspace{12pt}
    \Reshape\big( \cY_t , \{a, b\} \big) = \Fold_{K-1}\big( \cY_{t, a\sim b}, \{d_1, \dots, d_{a-1}, d_{a+1},\dots, d_{b-1}, d_{b+1}, \dots, d_{K}, d_{a}d_{b}\} \big) \\
    &=
    \Reshape\big( \cF_t ,\{a,b\} \big) \times_{i=1}^{K-2} \A_{\cA_i^\ast} \times_{K-1} \big(\A_b \otimes \A_a \big) + \Reshape\big( \cE_t, \{a,b\} \big),
\end{align*}
where the last line used similar arguments (omitted) as \eqref{eqn: fold_reshape_equality}. This implies $\cY_t$ follows \eqref{eqn: fm_wo_Kron} with a Kronecker structure product as $\big(\A_b \otimes \A_a \big) \in\cK_{d_a\times d_b}$.

Finally, consider any $\ell\geq 3$. We use an induction argument. With the definition of tensor reshape in Section~\ref{subsec: reshape} and the above for $\ell=2$, the initial case $\ell=3$ can be shown by
\begin{equation}
\label{eqn: reshape_proof_ell3}
\begin{split}
    &\hspace{12pt}
    \Reshape\big( \cY_t , \{a_1, a_2, a_3\} \big) =\Reshape\big[ \Reshape\big( \cY_t , \{a_2, a_3\} \big), \{a_1, K-1\} \big] \\
    &=
    \Reshape\Big\{ \Reshape\big( \cF_t ,\{a_2,a_3\} \big) \times_{i=1}^{K-2} \A_{[[K] \setminus \{a_2, a_3\} ]_i} \times_{K-1} \big(\A_{a_3} \otimes \A_{a_2} \big) \\
    &\hspace{12pt}
    + \Reshape\big( \cE_t, \{a_2,a_3\} \big) , \{a_1, K-1\}\Big\} \\
    &=
    \Reshape\Big\{ \Reshape\big( \cF_t ,\{a_2,a_3\} \big), \{a_1, K-1\}\Big\} \times_{i=1}^{K-3} \A_{[([K] \setminus \{a_2, a_3\} ) \setminus \{a_1\} ]_i} \times_{K-2} \big(\A_{a_3} \otimes \A_{a_2} \otimes \A_{a_1} \big) \\
    &\hspace{12pt}
    + \Reshape\Big\{ \Reshape\big( \cE_t, \{a_2,a_3\} \big) , \{a_1, K-1\}\Big\} \\
    &=
    \Reshape\big( \cF_t , \{a_1,a_2,a_3\} \big) \times_{i=1}^{K-3} \A_{[[K] \setminus \{a_1, a_2, a_3\}]_i} \times_{K-2} \big(\A_{a_3} \otimes \A_{a_2} \otimes \A_{a_1} \big) \\
    &\hspace{12pt}
    + \Reshape\big( \cE_t , \{a_1,a_2,a_3\} \big) ,
\end{split}
\end{equation}
where the last equality used again the definition of tensor shape, and note that $\big(\A_{a_3} \otimes \A_{a_2} \otimes \A_{a_1} \big) \in\cK_{d_{a_1}\times d_{a_2} \times d_{a_3}}$. Now if for all $\ell\in [L]$ with $L\geq 3$, \eqref{eqn: fm_w_Kron} implies \eqref{eqn: fm_wo_Kron} with $\A_{\textnormal{reshape}, K-\ell+1} \in \cK_{d_{a_1} \times \dots \times d_{a_\ell}}$, which is then also true for $\ell=L+1$ in a similar argument (omitted) as \eqref{eqn: reshape_proof_ell3}. This completes the induction.

Given any $\cA$, note also that if Assumption (F1) holds with $\cX_{\textnormal{reshape},t}$ and $\cF_{\textnormal{reshape},t}$ replaced by $\cX_t$ and $\cF_t$ respectively (with $\cX_t$ and $\cF_t$ from \eqref{eqn: F1_under_H0}), then it is immediate from the linearity of the reshape operator that
\begin{align*}
    \cF_{\textnormal{reshape},t} &= \Reshape\big( \cF_t , \cA \big) = \Reshape\Big( \sum_{w\geq 0} a_{f,w} \cX_{f,t-w} , \cA \Big) = \sum_{w\geq 0} a_{f,w} \Reshape\big( \cX_{f,t-w} , \cA \big),
\end{align*}
which implies $\cF_{\textnormal{reshape},t}$ follows Assumption (F1) by $\cX_{\textnormal{reshape},t} = \Reshape\big( \cX_{f,t} , \cA \big)$.

We have now proved that \eqref{eqn: fm_w_Kron} uniquely implies \eqref{eqn: fm_wo_Kron} with $\A_{\textnormal{reshape}, K-\ell+1} \in \cK_{d_{a_1} \times \dots \times d_{a_\ell}}$, and a version of Assumption (F1) on $\{\cF_t\}$ (from \eqref{eqn: F1_under_H0}) implies Assumption (F1) on $\cF_{\textnormal{reshape},t}$. It remains to show the other way around (for some $\cA$), but all the previous steps are reversible (note that in particular, the reshape operator is reversible as long as the dimension of the original tensor is known) and $\A_{\textnormal{reshape}, K-\ell+1} \in \cK_{d_{a_1} \times \dots \times d_{a_\ell}}$ ensures the existence of an appropriate set of low-rank matrices. Therefore, the proof for the theorem is completed. $\square$

\textbf{\textit{Proof of Theorem~\ref{thm: noise_aggregate_dist}.}}
Under $H_0$, consider \eqref{eqn: test_fm_w_Kron}. Since there exists $m\in [|\c{R}|]$ such that
\[
(\pi_{m,1}, \dots, \pi_{m,K-v+1}) = (r_v, \dots, r_K) ,
\]
we only consider such $m$ and simplify $\wh\cC_{m,t}$ as $\wh\cC_{t}$ and $\wh\cE_{m,t}$ as $\wh\cE_{t}$ (see the explanation at the beginning of Section~\ref{appendix: proof}). By \eqref{eqn: hat_cEt} and Assumption (E1),
\begin{align*}
    \wh\E_{t,(k^\ast)} &= \mat{k^\ast}{\wh\cE_{t}} = \mat{k^\ast}{(\cC_t - \wh\cC_{t}) + \c{F}_{e,t}\times_1 \A_{e,1}\times_2 \cdots \times_K\A_{e,K} + \bSigma_{\epsilon}\circ \bepsilon_t} \\
    &=
    (\cC_t - \wh\cC_{t})_{(k^\ast)} + \A_{e,k^\ast}\F_{e,t,(k^\ast)} \A_{e,\text{-}k^\ast}' + \bSigma_{\epsilon,(k^\ast)} \circ \bepsilon_{t,(k^\ast)} ,
\end{align*}
where $\A_{e,\text{-}k^\ast} := \A_{e,K} \otimes \dots \otimes \A_{e,k^\ast+1} \otimes \A_{e,k^\ast-1} \otimes \dots \otimes \A_{e,1}$. Hence, for any $t\in[T]$ and $j\in[d/d_k^\ast]$, we have
\begin{equation}
\label{eqn: residual_asymp_decomp}
\begin{split}
    \frac{1}{d_{k^\ast}} \sum_{i=1}^{d_{k^\ast}} \big(\wh{E}_{t,(k^\ast),ij}^2 - \Sigma_{\epsilon, (k^\ast), ij}^2 \big)
    &= \frac{1}{d_{k^\ast}} \big( \wh{\E}_{t,(k^\ast)}' \wh{\E}_{t,(k^\ast)} \big)_{jj} -\frac{1}{d_{k^\ast}} \sum_{i=1}^{d_{k^\ast}} \Sigma_{\epsilon, (k^\ast), ij}^2
    = \sum_{h=1}^6 \c{I}_{e,h}, \\
    \text{where} \;\;\;
    \c{I}_{e,1} &:=
    \frac{1}{d_{k^\ast}} \big\{ \big(\bSigma_{\epsilon,(k^\ast)} \circ \bepsilon_{t,(k^\ast)} \big)' \big(\bSigma_{\epsilon,(k^\ast)} \circ \bepsilon_{t,(k^\ast)} \big)\big\}_{jj} -\frac{1}{d_{k^\ast}} \sum_{i=1}^{d_{k^\ast}} \Sigma_{\epsilon, (k^\ast), ij}^2 , \\
    \c{I}_{e,2} &:=
    \frac{1}{d_{k^\ast}} \big\{ (\cC_t - \wh\cC_{t})_{(k^\ast)}' (\cC_t - \wh\cC_{t})_{(k^\ast)} \big\}_{jj} , \\
    \c{I}_{e,3} &:=
    \frac{1}{d_{k^\ast}} \big\{ \A_{e,\text{-}k^\ast} \F_{e,t,(k^\ast)}' \A_{e,k^\ast}' \A_{e,k^\ast} \F_{e,t,(k^\ast)} \A_{e,\text{-}k^\ast}'\big\}_{jj} , \\
    \c{I}_{e,4} &:=
    O_P\Big( d_{k^\ast}^{-1} \big\{ (\cC_t - \wh\cC_{t})_{(k^\ast)}' \A_{e,k^\ast} \F_{e,t,(k^\ast)} \A_{e,\text{-}k^\ast}' \big\}_{jj} \Big) , \\
    \c{I}_{e,5} &:=
    O_P\Big( d_{k^\ast}^{-1} \big\{ (\cC_t - \wh\cC_{t})_{(k^\ast)}' \big(\bSigma_{\epsilon,(k^\ast)} \circ \bepsilon_{t,(k^\ast)} \big) \big\}_{jj} \Big) , \\
    \c{I}_{e,6} &:=
    O_P\Big( d_{k^\ast}^{-1} \big\{ \A_{e,\text{-}k^\ast} \F_{e,t,(k^\ast)}' \A_{e,k^\ast}' \big(\bSigma_{\epsilon,(k^\ast)} \circ \bepsilon_{t,(k^\ast)} \big) \big\}_{jj} \Big) .
\end{split}
\end{equation}

Consider $\c{I}_{e,2}$. From Lemma~\ref{lemma: hat_rate_cFcC}, recall that
\begin{align*}
    &\hspace{12pt}
    \Big\{ \big(\vec{\wh\cC_{t}}\big)_i - \big( \vec{\cC_{t}} \big)_i \Big\}^2 \\
    &=
    O_P\Bigg( \max_{k\in[K]}\Bigg\{ d_k^{2(\delta_{k,1} - \delta_{k,r_k})} \Bigg( \frac{1}{T\dmk d_k^{1- \delta_{k,1}}} +\frac{1}{ d_k^{1+ \delta_{k, r_k}}} \Bigg) \frac{d^2}{g_s g_w} \Bigg\} +\frac{d}{g_w^2} \Bigg) ,
\end{align*}
which is the squared error for each entry of $\wh\cC_{t}$. With Assumption (R2), the above squared error rate is $o\big(1/\min_{k\in[K]} \{d_k\} \big) = o\big( d_{k^\ast}^{-1} \big)$. Hence $\c{I}_{e,2} =o_P\big( d_{k^\ast}^{-1} \big)$.

With Assumption (E1), $\|\A_{e,k^\ast}\|_F=O(1)$, $\|\A_{e,\text{-}k^\ast}\|_F=O(1)$ and $r_{e,k}$ for $k\in[K]$ are finite, so that $\c{I}_{e,3},\, \c{I}_{e,6} =O_P\big( d_{k^\ast}^{-1} \big)$. By the Cauchy--Schwarz inequality, immediately $\c{I}_{e,4} = O_P\big(\c{I}_{e,2}^{1/2}\cdot \c{I}_{e,3}^{1/2}\big) = o_P\big( d_{k^\ast}^{-1} \big)$.

Consider $\c{I}_{e,1}$, noting that
\[
\c{I}_{e,1} = \frac{1}{d_{k^\ast}} \sum_{i=1}^{d_{k^\ast}} \Sigma_{\epsilon,(k^\ast),ij}^2 \big( \epsilon_{t,(k^\ast),ij}^2 -1 \big) ,
\]
so that with Theorem 1 in \cite{AyvazyanUlyanov2023}, we have
\[
Z_{j} := \frac{d_{k^\ast}^{-1} \sum_{i=1}^{d_{k^\ast}} \Sigma_{\epsilon,(k^\ast),ij}^2 \big( \epsilon_{t,(k^\ast),ij}^2 -1 \big)} {\sqrt{d_{k^\ast}^{-2} \sum_{i=1}^{d_{k^\ast}} \text{Var}\big( \epsilon_{t,(k^\ast),ij}^2 \big) \Sigma_{\epsilon,(k^\ast),ij}^4}}
\xrightarrow{\c{D}} \c{N}(0,1) ,
\]
implying $\c{I}_{e,1}$ is of the rate $d_{k^\ast}^{-1/2}$ exactly. Note also that $Z_{j}$'s are independent of each other by Assumption (E1). It also follows that $\c{I}_{e,5} = O_P\big( \c{I}_{e,2}^{1/2}\cdot \c{I}_{e,1}^{1/2}\big) = o_P\big( d_{k^\ast}^{-3/4} \big)$. Finally, with \eqref{eqn: residual_asymp_decomp},
\begin{align*}
    \frac{\sum_{i=1}^{d_{k^\ast}} \big(\wh{E}_{t,(k^\ast),ij}^2 - \Sigma_{\epsilon, (k^\ast), ij}^2 \big)}{\sqrt{\sum_{i=1}^{d_{k^\ast}} \textnormal{Var}(\epsilon_{t, (k^\ast), ij}^2) \Sigma_{\epsilon, (k^\ast), ij}^4 }}
    &=
    \frac{d_{k^\ast}^{-1} \sum_{i=1}^{d_{k^\ast}} \big(\wh{E}_{t,(k^\ast),ij}^2 - \Sigma_{\epsilon, (k^\ast), ij}^2 \big)}{\sqrt{d_{k^\ast}^{-2} \sum_{i=1}^{d_{k^\ast}} \textnormal{Var}(\epsilon_{t, (k^\ast), ij}^2) \Sigma_{\epsilon, (k^\ast), ij}^4 }} \\
    &=
    Z_{j} (1+o_P(1)) \xrightarrow{p} Z_{j}
    \xrightarrow{\c{D}} \c{N}(0,1) .
\end{align*}
This shows the asymptotic result for $\wh\cE_t$ in Theorem~\ref{thm: noise_aggregate_dist} (i.e., the first asymptotic result). The result for $\wt\cE_t$ can be shown in the same manner, except that for $\c{I}_{e,2}$, we use (R2) on the squared error for each entry of $\wt\cC_{\textnormal{reshape},t}$ which is the following from Lemma~\ref{lemma: Q_tilde_rate},
\begin{align*}
    &\hspace{12pt}
    \Big\{ \big(\vec{\wt\cC_{\textnormal{reshape},t}}\big)_i - \big( \vec{\cC_{\textnormal{reshape},t}} \big)_i \Big\}^2 \\
    &=
    O_P\Bigg( \max_{j\in[v-1]}\Bigg\{ d_j^{2(\delta_{j,1} - \delta_{j,r_j})} \Bigg( \frac{1}{T\dmk d_j^{1- \delta_{j,1}}} +\frac{1}{ d_j^{1+ \delta_{j, r_j}}} \Bigg) \frac{d^2}{\gamma_s \gamma_w} \Bigg\} \\
    &\hspace{57pt}
    + d_V^{2(\delta_{V,1} - \delta_{V,r_V})} \Bigg( \frac{1}{Td d_V^{- \delta_{V,1}}} +\frac{1}{ d_V^{1+ \delta_{V, r_V}}} \Bigg) \frac{d^2}{\gamma_s \gamma_w}
    + \frac{d}{\gamma_w^2} \Bigg) .
\end{align*}
Under $H_0$, we arrive at the same conclusion with the same $Z_{j}$'s. Moreover, the asymptotic result for $\wh\cE_t$ holds true under $H_1$. This concludes the proof of the theorem. $\square$

\textbf{\textit{Proof of Theorem~\ref{thm: test}.}}
Under $H_0$, Theorem~\ref{thm: noise_aggregate_dist} implies that for each $t\in[T]$ and $j\in [d/d_{k^\ast}]$, there exists such $\sum_{i=1}^{d_{k^\ast}} \wh{E}_{m,t,(k^\ast),ij}^2$ distributed asymptotically the same as $\sum_{i=1}^{d_{k^\ast}} \wt{E}_{t,(k^\ast),ij}^2$. With the definition $\wh{q}_{y,m,j}(\alpha) := \inf\big\{ c\;|\;\b{F}_{y,m,j}(c)\geq 1-\alpha \big\}$, we arrive at the existence of such $\wh{q}_{y,m,j}(\alpha)$ asymptotically the same as $\wh{q}_{x,j}(\alpha)$. In this case, we have
\begin{align*}
\alpha \geq \b{P}_{y,m,j}\big( y_{m,j,t} \geq \wh{q}_{y,m,j}(\alpha) \big) \rightarrow \b{P}_{y,m,j}\big( y_{m,j,t} \geq \wh{q}_{x,j}(\alpha) \big).
\end{align*}
With the above, \eqref{eqn: min_test_probability} is direct from the definition of $\b{P}_{y,m,j}$. This completes the proof of the theorem. $\square$

\textbf{\textit{Proof of Theorem~\ref{thm: reshape2}.}}
Let $\cY_t \in\b{R}^{d_1\times \dots \times d_K}$ be an order-$K$ tensor and $\cA =\{a_1,\dots,a_\ell\}$. Then each $\Reshape(\cY_t, \cA)$ is an order-$(K-\ell+1)$ tensor. If $\big\{ \Reshape(\cY_t, \cA) \big\}$ has a Kronecker product structure, then Theorem~\ref{thm: reshape} allows us to write for $t\in[T]$,
\begin{equation}
\label{eqn: reshape2_proof_1}
\Reshape(\cY_t, \cA) = \cF_{\textnormal{reshape},t} \times_{j=1}^{K-\ell} \A_{\textnormal{reshape},j} \times_{K-\ell+1} \A_{\textnormal{reshape}, K-\ell+1} + \cE_{\textnormal{reshape},t} ,
\end{equation}
for some core factor $\{\cF_{\textnormal{reshape},t}\}$, loading matrices $\{\A_{\textnormal{reshape},j} \}_{j\in [K-\ell+1]}$, and noise $\cE_{\textnormal{reshape},t}$. Immediately by Definition~\ref{def: kron_structure}, if $\A_{\textnormal{reshape}, K-\ell+1} \in \cK_{d_{a_1} \times \dots\times d_{a_\ell}}$, then $\cY_t$ has a Kronecker product structure; otherwise, $\{\cY_t\}$ has no Kronecker product structure along $\cA$.

It remains to show that if $\{\cY_t\}$ either has a Kronecker product structure or has no Kronecker product structure along some set $\cA^\ast$, then $\big\{ \Reshape(\cY_t, \cA) \big\}$ with $\cA^\ast \subseteq \cA$ has a Kronecker product structure. The case where $\{\cY_t\}$ has a Kronecker product structure is trivial. We then only need to consider that $\{\cY_t\}$ has no Kronecker product structure along $\cA^\ast$, which implicitly assume a factor model of $\{\cY_t\}$ along $\cA^\ast$ by Definition~\ref{def: kron_structure}. Without loss of generality, let $\cA^\ast :=\{K-g+1, \dots, K\}$, otherwise redefine the mode indices of $\cY_t$. For the set $\cA$ with $\cA^\ast \subseteq \cA$, we now read $\cA=\{a_1, \dots, a_{\ell-g}, K-g+1, \dots, K\}$. Using the last property of tensor reshape in Section~\ref{subsec: reshape} (which can be easily seen by induction), we have
\begin{equation}
\label{eqn: reshape2_proof_2}
\Reshape(\cY_t, \cA) = \Reshape\big\{ \Reshape(\cY_t, \cA^\ast), \{a_1, \dots, a_{\ell-g}, K-g+1\}\big\} .
\end{equation}
Now that $\{\cY_t\}$ has no Kronecker product structure along $\cA^\ast$, similar to the form \eqref{eqn: reshape2_proof_1}, we have for $t\in[T]$,
\[
\Reshape(\cY_t, \cA^\ast) = \cF_{\textnormal{reshape},t}^\ast \times_{j=1}^{K-g} \A_{\textnormal{reshape},j}^\ast \times_{K-g+1} \A_{\textnormal{reshape}, K-g+1}^\ast + \cE_{\textnormal{reshape},t}^\ast ,
\]
which implies the time series $\{\Reshape(\cY_t, \cA^\ast)\}$ follows a Tucker-decomposition TFM. According to Theorem~\ref{thm: reshape}, $\{\Reshape(\cY_t, \cA^\ast)\}$ follows a factor model along any index set of $\{\Reshape(\cY_t, \cA^\ast)\}$ (with Kronecker product structure; but the factor model form is sufficient for our claim). In particular, with the index set $\{a_1, \dots, a_{\ell-g}, K-g+1\}$, we conclude from \eqref{eqn: reshape2_proof_2} that $\{\Reshape(\cY_t, \cA)\}$ has a Kronecker product structure. This also completes the proof of the theorem. $\square$

\newpage
\bibliographystyle{apalike}
\bibliography{ref}

\end{document}